\documentclass[11pt]{article}

\usepackage{amsmath,amssymb}

\newcommand{\epsln}{\varepsilon}

\newcommand{\rarrow}{\Longrightarrow}

\newcommand{\tr}{{\rm tr}}

\def\epsln{\varepsilon}
\def\bgeq{\begin{equation}}
\def\edeq{\end{equation}}
\def\bgar{\begin{array}}
\def\edar{\end{array}}

\title{Lie Group Action and Stability Analysis of Stationary Solutions for a
  Free Boundary Problem Modelling Tumor Growth}
\author{Shangbin Cui}
\date{\small Institute of Mathematics, Sun Yat-Sen University, Guangzhou,
  Guangdong 510275,\\ People's Republic of China. E-mail: cuisb3@yahoo.comc.cn}

\begin{document}

\maketitle

\begin{abstract}
  In this paper we study asymptotic behavior of solutions for a
  multidimensional free boundary problem modelling the growth of nonnecrotic
  tumors. We first establish a general result for differential equations in
  Banach spaces possessing a local Lie group action which maps a solution into
  new solutions. We prove that a center manifold exists under certain
  assumptions on the spectrum of the linearized operator without assuming that
  the space in which the equation is defined is of either $D_A(\theta)$ or
  $D_A(\theta,\infty)$ type. By using this general result and making delicate
  analysis of the spectrum of the linearization of the stationary free boundary
  problem, we prove that if the surface tension coefficient $\gamma$ is larger
  than a threshold value $\gamma^\ast$ then the unique stationary
  solution is asymptotically stable modulo translations, provided the constant
  $c$ representing the ratio between the nutrient diffusion time and the
  tumor-cell doubling time is sufficiently small, whereas if $\gamma<
  \gamma^\ast$ then this stationary solution is unstable.



\end{abstract}

\section{Introduction}
\setcounter{equation}{0}

  This paper aims at studying asymptotic behavior of solutions of the following
  free boundary problem:
$$
  c\partial_t\sigma=\Delta\sigma-f(\sigma),
  \quad  x\in \Omega(t),\;\; t>0,
\eqno{(1.1)}
$$
$$
  -\Delta p = g(\sigma),\quad x\in \Omega(t),\;\; t>0,
\eqno{(1.2)}
$$
$$
  \sigma = \bar{\sigma},\quad  x\in\partial\Omega(t),\;\; t>0,
\eqno{(1.3)}
$$
$$
   p = \gamma\kappa,\quad x\in\partial\Omega(t),\;\; t>0,
\eqno{(1.4)}
$$
$$
  {\mathbf V}=-\partial_{\mathbf n} p,\quad x\in\partial\Omega(t),\;\; t>0,
\eqno{(1.5)}
$$
$$
   \sigma(x,0)=\sigma_0(x),\quad x\in\Omega_0,
\eqno{(1.6)}
$$
$$
   \Omega(0)=\Omega_0.
\eqno{(1.7)}
$$
  Here $\sigma=\sigma(x,t)$ and $p=p(x,t)$ are unknown functions defined on
  the space-time manifold $\cup_{t\geq 0}(\overline{\Omega(t)}\times \{t\})$,
  and $\Omega(t)$ is an a priori unknown bounded time-dependent domain in
  ${\mathbb R}^n$, whose boundary $\partial\Omega(t)$ has to be determined
  together with the unknown functions $\sigma$ and $p$. Besides, $f$ and
  $g$ are given functions, $c$, $\bar{\sigma}$ are $\gamma$ are positive
  constants, $\kappa$, ${\mathbf V}$ and ${\mathbf n}$ are the mean curvature,
  the normal velocity and the unit outward normal vector of $\partial\Omega(t)$,
  respectively, and $\sigma_0$, $\Omega_0$ are given initial data of
  $\sigma=\sigma(\cdot,t)$ and $\Omega=\Omega(t)$, respectively. The sign of
  $\kappa$ is fixed on by the condition that $\kappa\geq 0$ at points where
  $\partial\Omega(t)$ is convex with regard to $\Omega(t)$.

  The above problem arises from recently developed subject of tumor growth
  modelling. It models the growth of tumors cultivated in laboratory or
  so-called {\em multicellular spheroids} (\cite{AdmBel}, \cite{Byrne},
  \cite{ByrCha}, \cite{Green}, \cite{KimStHa}, \cite{Klie}, \cite{Suther}). In
  this model $\Omega(t)$ represents
  the domain occupied by the tumor at time $t$, $\sigma$ and $p$ stand for the
  nutrient concentration and the tumor tissue pressure, respectively, and
  $f(\sigma)$, $g(\sigma)$ are the nutrient consumption rate and the tumor cell
  proliferation rate, respectively. It is assumed that all tumor cells are
  alive and dividable, and their density is constant, so that in $f$ and $g$ no
  cell density argument is involved. It is also assumed that the tumor is
  cultivated in a solution of nutrition materials whose concentration keeps
  constant in the process of cultivation, and $\bar{\sigma}$ reflects this
  constant nutrient supply to the tumor. The term $\gamma\kappa$ on the
  right-hand side of (1.4) stands for surface tension of the tumor. The
  equation (1.5) reflects the fact that the normal velocity of the tumor surface
  is equal to the normal component of the movement velocity of tumor cells
  adjacent to the surface. For more details of the modelling we refer the reader
  to see the references \cite{AdmBel}, \cite{Byrne}, \cite{ByrCha}, \cite{Cui1},
  \cite{Cui3}, \cite{Cui6}--\cite{CuiEsc2} and \cite{Green}. Here we point out
  that, by rescaling which we have pre-assumed and did not particularly mention,
  the constant $c$ represents the ratio between the nutrient diffusion time and
  the tumor-cell doubling time, so that $c\ll 1$, cf. \cite{AdmBel}, \cite{Byrne},
  and \cite{ByrCha}. Finally, we make the following assumptions on the functions
  $f$ and $g$:
\begin{flushright}
\begin{minipage}[t]{15cm}
\begin{flushleft}

  $(A1)$ $f\in C^\infty[0,\infty)$, $f'(\sigma)>0$ for $\sigma\geq 0$ and
  $f(0)=0$.

  $(A2)$ $g\in C^\infty[0,\infty)$, $g'(\sigma)>0$ for $\sigma\geq 0$ and there
  exists a number $\tilde{\sigma}>0$ such that\\ $\qquad$ $g(\tilde{\sigma})=0$
  ($\rarrow$ $g(\sigma)<0$ for $0\leq\sigma<\tilde{\sigma}$ and $g(\sigma)>0$
  for $\sigma>\tilde{\sigma}$).

  $(A3)$ $\tilde{\sigma}<\bar{\sigma}$.
\end{flushleft}
\end{minipage}
\end{flushright}

\noindent
  These assumptions are based on biological considerations, see  \cite{Cui3},
   \cite{CuiEsc1} and \cite{CuiEsc2}.

  Local well-posedness of the above problem has been recently established by
  the author in a more general framework in the reference \cite{Cui6} by using
  the analytic semigroup theory, which extends and modifies an earlier work of
  Escher \cite{Escher} for the special case that $f(\sigma)=f(\sigma)$ but
  $g(\sigma)=\mu(\sigma-\tilde{\sigma})$. In this paper we consider the more
  difficult topic of asymptotic behavior of the solution. More precisely, from
  \cite{Cui3} and \cite{CuiEsc1} we know that under the above assumptions
  $(A1)$--$(A3)$, the system (1.1)--(1.5) has a radially symmetric stationary
  solution $(\sigma_s,p_s,\Omega_s)$, which is unique up to translations and
  rotations of the coordinate of ${\mathbb R}^n$ and globally asymptotically
  stable under radially symmetric perturbations. This paper aims at studying
  the following question: Is $(\sigma_s,p_s,\Omega_s)$ also asymptotically
  stable under non-symmetric perturbations?

  We first make a short review to previous work on this topic. Rigorous
  analysis of free boundary problems of partial differential equations arising
  from tumor growth modelling has attracted a lot of attention during the past
  several years, and many interesting results have been systematically derived,
  cf. \cite{BazFri1}, \cite{BazFri2}, \cite{ChenCuiF}--\cite{CuiWei},
  \cite{Escher}, \cite{FriHu}--\cite{FriRei3}, and the references cited therein.
  As far as the problem (1.1)--(1.7) and its certain more specific
  forms are concerned, we cite the references \cite{BazFri1}, \cite{BazFri2},
  \cite{Cui1}, \cite{Cui3},  \cite{Cui6}--\cite{CuiEsc2}, \cite{Escher},
  \cite{FriHu}--\cite{FriRei2}. In particular, in \cite{FriRei1} Friedman and
  Reitich considered radially symmetric version of the problem (1.1)--(1.7) in
  the special case that $f(\sigma)=\lambda\sigma$ and $g(\sigma)=\mu(\sigma-
  \tilde{\sigma})$. Under the assumption $(A3)$, they proved the following
  results: (1) The problem is globally well-posed. (2) There exists a unique
  stationary solution. (3) For $c$ sufficiently small this stationary solution
  is globally asymptotically stable. (4) For $c$ large the stationary solution
  is unstable. The author of the present paper has recently extended the
  assertions (1), (2), (3) to the general case that $f$ and $g$ are general
  functions satisfying the conditions $(A1)$--$(A3)$, see \cite{Cui3}. The
  general non-symmetric version of (1.1)--(1.7) in the special case that
  $f(\sigma)=\lambda\sigma$ and $g(\sigma)=\mu(\sigma-\tilde{\sigma})$ has also
  been systematically studied by Friedman and his collaborators. Bazaliy and
  Friedman investigated local well-posedness of the time-dependent problem in
  the reference \cite{BazFri1}. In \cite{BazFri2} they studied asymptotic
  behavior of the solution starting from a neighborhood of the unique radially
  symmetric stationary solution ensured by the above assertion (2), and proved
  that, for $c=1$, $\lambda=1$, $\gamma=1$ and $\mu$ sufficiently small, the
  radially symmetric stationary solution is (locally) asymptotically stable
  under non-radial perturbations. This work was recently refined by Friedman
  and Hu \cite{FriHu}. They proved that, again for $c=1$, $\lambda=1$ and
  $\gamma=1$, there exists a threshold value $\mu^\ast>0$, such that for
  $0<\mu<\mu^\ast$ the radially symmetric stationary solution is (locally)
  asymptotically stable under non-radial perturbations, while for $\mu>
  \mu^\ast$ this stationary solution is unstable.

  In a recent work of the present author jointly with Escher \cite{CuiEsc2},
  the problem (1.1)--(1.7) with general functions $f$ and $g$ satisfying
  $(A1)$--$(A3)$ but $c=0$ was studied. We proved that there exists a
  threshold value $\gamma_\ast>0$, the supremum of all bifurcation points
  $\gamma_k$ ($k=2,3,\cdots$, see \cite{CuiEsc1}), such that
  if $\gamma>\gamma_\ast$ then the
  radially symmetric stationary solution $(\sigma_s,p_s,\Omega_s)$ is (locally)
  asymptotically stable {\em modulo translations}, i.e., any solution starting
  from a small neighborhood of $(\sigma_s,p_s,\Omega_s)$ is global and, as
  $t\to\infty$, it converges to either $(\sigma_s,p_s,\Omega_s)$ or an adjacent
  stationary solution $(\sigma_s',p_s',\Omega_s')$ obtained by translating
  $(\sigma_s,p_s,\Omega_s)$ (recall that any translation of $(\sigma_s,p_s,
  \Omega_s)$ is still a stationary solution), whereas if $\gamma<\gamma_\ast$
  then $(\sigma_s,p_s,\Omega_s)$ is unstable.

  In this paper we want to extend the above result of \cite{CuiEsc2} for the
  degenerate case $c=0$ to the more difficult non-degenerate case $c\neq
  0$, assuming that $c$ is sufficiently small. The main idea of analysis is
  the same with that of \cite{CuiEsc2}, namely, we shall first reduce the PDE
  problem into a differential equation in a Banach space and next use the
  abstract geometric theory for parabolic differential equations in Banach
  spaces to get the desired result. However, unlike in \cite{CuiEsc2} where
  we used the well-developed center manifold theorem by Da Prato and Lunardi
  \cite{DaPLun} to make the analysis, in this paper we shall have to first
  establish a new center manifold theorem, because the above-mentioned center
  manifold of Da Prato and Lunardi is not applicable to the case $c\neq 0$.
  The reason is as follows. Recall that the center manifold theorem of Da Prato
  and Lunardi requires the Banach space in which the differential equation is
  considered must be of the type either $D_A(\theta)$, the continuous
  interpolation space, or $D_A(\theta,\infty)$, the real interpolation space
  of the type $(\theta,\infty)$ ($0<\theta<1$). Such spaces cannot be reflexive
  (cf. \cite{Baillon}, \cite{Simon}). In the degenerate case $c=0$
  the reduced equation contains only the unknown function $\rho$ defining the
  free boundary $\partial\Omega(t)$, which is a quasi-linear parabolic
  pseudo-differential equation on a compact manifold, so that no boundary
  conditions appear and we can thus work on the little H\"{o}lder space
  $h^{m+\alpha}$ which is of the type $D_A(\theta)$. In the present
  non-degenerate case $c\neq 0$, however, since the reduced equation contains
  not only $\rho$ but also the unknown $\sigma$, the Dirichlet boundary
  condition for $\sigma$ renders it impossible for us to work on a space of
  the type either $D_A(\theta)$ or $D_A(\theta,\infty)$.

  To remedy this deficiency, in this paper we shall first establish a new
  center manifold theorem which removes this very
  restrictive assumption on the space $X$, but instead we shall assume that the
  equation admits a local Lie group action by which a solution is mapped into
  new solutions. We shall show that the phase diagram of a differential
  equation possessing such a Lie group action has a very nice structure and its
  center manifold can be very easily obtained. In particular, this new center
  manifold theorem does not make any additional assumption on the structure of
  the space $X$. Since the differential equation reduced from the the problem
  (1.1)--(1.7) naturally possesses a Lie group action induced by translations
  of the coordinate of ${\mathbb R}^n$, by using this new center manifold
  result we are able to make analysis in the framework of Sobolev and Besov
  spaces. Our final result says that similar assertions as for the case $c=0$
  also hold for the case that $c$ is non-vanishing but very small, and this
  result will be established in the space $W^{m-1,q}\times W^{m-3,q}\times
  B^{m-1/q}_{qq}$, where $W^{m-1,q}$ and $B^{m-1/q}_{qq}$ represent the Sobolev
  and Besov spaces, respectively.

  It should be noted that our center manifold theorem for differential equations
  in Banach spaces possessing Lie group action established in this paper not
  only works for the tumor model (1.1)--(1.7) as well as its special form of
  the case $c=0$, but also applies to other problems such as the Hele-Shaw
  problem. Thus, the center manifold theorem established in this paper has its
  own theoretic importance. More applications of this result will be given
  in our future work.

  To give a precise statement of our main result, let us first introduce some
  notation. Recall that the radially symmetric stationary
  solution $(\sigma_s,p_s,\Omega_s)$ of (1.1)--(1.5), where $\Omega_s=\{\, r<
  R_s\}$ with $r=|x|$, is the unique solution of the following free boundary
  problem:
$$
  \sigma_s''(r)+{n\!-\!1\over r}\sigma_s'(r)=f(\sigma_s(r)), \quad
  0<r<R_s,
\eqno{(1.8)}
$$
$$
  p_s''(r)+{n\!-\!1\over r}p_s'(r)=-g(\sigma_s(r)), \quad
  0<r<R_s,
\eqno{(1.9)}
$$
$$
  \sigma_s'(0)=0, \quad \sigma_s(R)=\bar{\sigma},
\eqno{(1.10)}
$$
$$
  p_s'(0)=0, \quad p_s(R_s)={\gamma\over R_s},
\eqno{(1.11)}
$$
$$
   p_s'(R_s)=0.
\eqno{(1.12)}
$$
  For $z\in {\mathbb R}^n$, we denote
$$
  \sigma_{s}^{z}(x)=\sigma_s(|x-z|), \quad p_{s}^{z}(x)=p_s(|x-z|),
 \quad
  \Omega_{s}^{z}=\{\; x\in {\mathbb R}^n:\;|x-z|<R_s\}.
$$
  Clearly, for any $z\in {\mathbb R}^n$ the triple $(\sigma_{s}^{z},p_{s}^{z},
  \Omega_{s}^{z})$ is a stationary solution of the system (1.1)--(1.5). If
  $|z|$ is sufficiently small then there exists a unique $\rho_{s}^{z}\in
  C^\infty({\mathbb S}^{n-1})$ which is sufficiently close to the constant
  function $R_s$, such that
$$
  \Omega_{s}^{z}=\{\, r<\rho_{s}^{z}
  (\omega),\;\omega\in {\mathbb S}^{n-1}\}.
$$
  Since we shall only consider solutions of (1.1)--(1.7) which are close to
  the stationary solution $(\sigma_s,p_s,\Omega_s)$, we can write $\Omega(t)$
  as $\Omega(t)=\{\,r<\rho(\omega,t),\;\omega\in {\mathbb S}^{n-1}\}$ for some
  $\rho(\cdot,t)\in C({\mathbb S}^{n-1})$ for every $t>0$, and, correspondingly,
  we write $\Omega_0$ as $\Omega_0=\{\, r<\rho_0(\omega),\;\omega\in {\mathbb
  S}^{n-1}\}$, where $\rho_0\in C({\mathbb S}^{n-1})$. Finally, from
  \cite{CuiEsc1} we know that the linearization of the stationary version of
  (1.1)--(1.5) has infinite many eigenvalues $\gamma_k$, $k=2,3,\cdots$, which
  are all positive and converge to zero as $k\to\infty$. As in \cite{CuiEsc2}
  we set
$$
  \gamma_\ast=\max\{\gamma_k,\; k=2,3,\cdots\}.
$$

  The main result of this paper is as follows:
\medskip

  {\bf Theorem 1.1}\ \ {\em If $\gamma>\gamma_\ast$ then there exists a
  corresponding $c_0>0$ such that for any $0<c<c_0$, the stationary solution
  $(\sigma_s,p_s,\Omega_s)$ of (1.1)--(1.5) is asymptotically stable modulo
  translations in the following sense: There exists $\epsln>0$ such that for
  any $\rho_0\in B^{m-1/q}_{qq}({\mathbb S}^{n-1})$ and $\sigma_0\in W^{m,q}
  (\Omega_0)$ $(m\in {\mathbb N}$, $m\geq 5$, $n/(m-4)<q<\infty)$ satisfying
$$
  \|\rho_0-R_s\|_{B^{m-1/q}_{qq}({\mathbb S}^{n-1})}<\epsln,\quad
  \|\sigma_0-\sigma_s\|_{W^{m,q}(\Omega_0)}<\epsln, \quad
  \sigma_0|_{\partial\Omega_0}=\bar{\sigma},
$$
  the problem (1.1)--(1.7) has a unique solution $(\sigma,p,\Omega)$ for all
  $t\geq 0$, and there exists $z\in {\mathbb R}^n$ uniquely determined by
  $\rho_0$ and $\sigma_0$ such that
$$
  \|\sigma(\cdot,t)-\sigma_{s}^{z}\|_{W^{m\!-\!1,q}(\Omega(t))} +
  \|p(\cdot,t)-p_{s}^{z}\|_{W^{m\!-\!3,q}(\Omega(t))}
  +\|\rho(\cdot,t)-\rho_{s}^{z}\|_{B^{m\!-\!1/q}_{qq}
  ({\mathbb S}^{n\!-\!1})}\leq Ce^{-\kappa t}
$$
  for some $C>0$, $\kappa>0$ and all $t\geq 0$.
  If $\gamma<\gamma_\ast$ then there also exists a corresponding $c_0>0$ such
  that for any $0<c<c_0$, $(\sigma_s,p_s,\Omega_s)$ is unstable.}
\medskip

  {\bf Remark 1.1}.\ \ By the assertion (4) of Friedman and Reitich reviewed
  before, we see that the condition $c<c_0$ cannot be removed. Besides, as we
  mentioned earlier, though we only consider solutions in $W^{m-1,q}\times
  W^{m-3,q}\times B^{m-1/q}_{qq}$, a similar result surely also holds for
  solutions in the space $C^{m+\alpha}\times C^{m-2+\alpha}\times C^{m+\alpha}$.
  In addition, the conditions $m\geq 5$ and $n/(m-4)<q<\infty$ can be weakened
  upto $m\geq 3$ and $n/(m-2)<q<\infty$. To achieve this improvement we need a
  modified version of Theorem 2.1 of the next section; see Remark 2.1 in the
  end of Section 2.

  The proof of the above theorem will be given in the last section of this
  paper, after step-by-step preparations in Sections 2--6. The layout of
  the rest part is as follows. In Section 2 we establish the general result for
  differential equations in Banach spaces mentioned earlier. In Section 2 we
  first use the so-called Hanzawa transformation to transform the problem
  (1.1)--(1.7) into an equivalent problem on the fixed domain $\Omega_s$, which
  for simplicity of notation will be assumed to be the unit sphere ${\mathbb
  B}^n$ later on, and next we further reduce the PDE problem into a
  differential equation in the Banach space $W^{m\!-\!3,q}({\mathbb B}^n)\times
  B^{m-3-1/q}_{qq}({\mathbb S}^{n-1})$ for the unknowns $(\sigma,\rho)$. In
  Section 4 we construct Lie group action for the reduced differential equation.
  In Section 5 we compute the linearization of the reduced equation. Section 6
  aims at studying the spectrum of the linearized problem. In the last section
  we complete the proof of Theorem 1.1.

\section{An abstract result}

  Let $X$ and $X_0$ be two Banach spaces such that $X_0\hookrightarrow X$.
  $X_0$ need not be dense in $X$. Let ${\mathcal O}$ be an open subset of
  $X_0$. Let $F\in C^{2-0}({\mathcal O},X)$, i.e. $F\in C^1({\mathcal O},X)$
  and $F'$ ($=DF=$ the Fr\'{e}chet derivative of $F$) is Lipschitz continuous.
  In this section we consider the initial value problem
$$
\left\{
\begin{array}{l}
  u'(t)=F(u(t)), \quad t>0,\\
  u(0)=u_0,
\end{array}
\right.
\eqno{(2.1)}
$$
  where $u_0\in {\mathcal O}$. By a {\em solution} of (2.1) we mean a solution
  of the class $u\in C([0,T),X)\cap C((0,T),{\mathcal O})\cap L^\infty((0,T),
  {\mathcal O})\cap C^1((0,T),X)$ defined in a maximal existence interval
  $I=[0,T)$ ($0<T\leq\infty$), which satisfies (2.1) in $[0,T)$ in usual sense
  and is not extendable. If $u$ satisfies the stronger condition $u\in C([0,T),
  {\mathcal O})\cap C^1([0,T),X)$ then we call it a {\em strict solution}.
  Later on we shall denote by $u(t,u_0)$ the solution of (2.1) when it exists
  and is unique. We always assume that for some $u_s\in {\mathcal O}$ there
  holds $F(u_s)=0$, so that $u(t)=u_s$, $t\geq 0$, is a stationary solution of
  the equation $u'=F(u)$. We want to study asymptotic stability of $u_s$.
  Our first assumption is as follows:
\medskip

\begin{minipage}[t]{14cm}
{\em
  $(B_1)$  $A=F'(u_s)$ is a sectorial operator in $X$ with domain $X_0$, and
  the graph norm of $A$ is equivalent to the norm of $X_0$: $\|u\|_{X_0}\sim
  \|u\|_X+\|Au\|_X$.
}
\end{minipage}
\medskip

\noindent
  Next, we consider some invariance property of $F$. Let $G$ be a {\em local
  Lie group of dimension $n$} in the sense of L. S. Pontrjagin \cite{Pontr}.
  Let ${\mathcal O}'$ be an open subset of $X$ such that ${\mathcal O}\subseteq
  {\mathcal O}'$. Let ${\mathcal O}_1$ be an open subset of $X_0$ contained in
  ${\mathcal O}$, and ${\mathcal O}'_1$ be an open subset of $X$ contained in
  ${\mathcal O}'$, such that $u_s\in {\mathcal O}_1\subseteq {\mathcal O}'_1$.
  We assume that there is a continuous mapping $p:G \times {\mathcal O}'_1\to
  {\mathcal O}'$, such that

\begin{enumerate}
\item[]  $(i)$ $p(G\times {\mathcal O}_1)\subseteq {\mathcal O}$;

\item[]  $(ii)$ $p(e,u)=u$ for every $u\in {\mathcal O}'_1$, where $e$ denotes
  the unit of $G$, and $p(\sigma,p(\tau,u))=p(\sigma\tau,u)$ for any $u\in
  {\mathcal O}'_1$ and $\sigma,\tau\in G$ such that $\sigma\tau$ is well-defined
  and $p(\tau,u)\in {\mathcal O}'_1$;

\item[]  $(iii)$ If $\sigma,\tau\in G$ such that $p(\sigma,u)=p(\tau,u)$ for
  some $u\in {\mathcal O}'_1$ then $\sigma=\tau$.

\item[]  $(iv)$ For any $\sigma\in G$, the mapping $u\to p(\sigma,u)$ from
  ${\mathcal O}'_1$ to ${\mathcal O}'$ is Fr\'{e}chet differentiable at every
  point in ${\mathcal O}_1$, and $[u\to D_u p(\sigma,u)]\in C({\mathcal O}_1,
  L(X,X))$.

\item[]  $(v)$ For any $u\in {\mathcal O}_1$, the mapping $\sigma\to p(\sigma,
  u)$ from $G$ to ${\mathcal O}$ is continuously Fr\'{e}chet differentiable
  when regarded as a mapping from $G$ to $X$ ($\Longrightarrow D_\sigma
  p(\sigma,u)\in L(T_\sigma(G),X)$, and $[\sigma\to p(\sigma,u)]\in C^1(G,X)$).
  Moreover, ${\rm rank}D_\sigma p(\sigma,u)=n$ for every $\sigma\in G$ and
  $u\in {\mathcal O}_1$.
\end{enumerate}

\noindent
  Later on we denote $S_\sigma(u)=p(\sigma,u)$ for $\sigma\in G$ and $u\in
  {\mathcal O}_1$. Our second assumption is as follows:
\medskip

\begin{minipage}[t]{14cm}
{\em
  $(B_2)$ There is a local Lie group $G$ satisfying the properties $(i)$--$(v)$,
  such that for any $u\in {\mathcal O}_1$ and $\sigma\in G$ there holds
}
\end{minipage}
$$
  F(S_\sigma(u))=DS_\sigma(u)F(u).
\eqno{(2.2)}
$$

\noindent
  This assumption has some obvious inferences. First, it implies that for any
  $u_0\in {\mathcal O}_1$ and $\sigma\in G$ there holds $u(t,S_\sigma(u_0))=
  S_\sigma(u(t,u_0))$, namely, if $t\to u(t)$ is a solution of the equation
  $u'=F(u)$ with initial value $u_0$, then $t\to S_\sigma(u(t))$ is also a
  solution of this equation, with initial value $S_\sigma(u_0)$. In particular,
  for any $\sigma\in G$, $S_\sigma(u_s)$ is a stationary solution of $u'=
  F(u)$. Next, if $u_s$ is more regular than $(v)$ in the sense that $[\sigma
  \to p(\sigma,u_s)]\in C^1(G,X_0)$ (so that $D_\sigma p(\sigma,u_s)\in
  L(T_\sigma(G),X_0)$ for any $\sigma\in G$), then by
  differentiating the relation $F(S_\sigma(u_s))=0$ in $\sigma$ at $\sigma=e$
  we see that $DF(u_s)D_\sigma p(e,u_s)\xi=0$ for any $\xi\in T_e(G)$, so that
  $A=DF(u_s)$ is degenerate, and $\dim{\rm Ker}A\geq n$. We now assume
  that
\medskip

\begin{minipage}[t]{14cm}
{\em
  $(B_3)$ $[\sigma\to p(\sigma,u_s)]\in C^1(G,X_0)$, $\dim{\rm Ker}A=n$, and
  the induced operator $\overline{A}:X_0/{\rm Ker}A\to X/{\rm Ker}A$ of $A$ is
  an isomorphism.
}
\end{minipage}
\medskip

\noindent
  Here and throughout this paper, by isomorphism from a Banach space $X_1$ to
  another Banach space $X_2$ we mean a linear mapping $T:X_1\to X_2$ such that
  it is an 1-1 correspondence, and both $T$ and $T^{-1}$ are continuous (i.e.,
  $T$ is not merely a linear isomorphism, but a topological homeomorphism as
  well). Finally, we assume that
\medskip

\begin{minipage}[t]{14cm}
{\em
  $(B_4)$ $\omega_-\equiv-\sup\{{\rm Re}\lambda:\lambda\in\sigma(A)\backslash
  \{0\}\}=-\sup\{{\rm Re}\lambda:\lambda\in\sigma(\overline{A})\}>0$.
}
\end{minipage}
\medskip

\noindent
  We point out that the condition $(B_3)$ is equivalent to the following
  condition:
\medskip

\begin{minipage}[t]{14cm}
{\em
  $(B_3')$  $\dim{\rm Ker}A=n$, ${\rm Range}A$ is closed in $X$, and
  $X={\rm Ker}A\oplus {\rm Range}A$.
}
\end{minipage}
\medskip

\noindent
  The proof of equivalence of $(B_3)$ with $(B_3')$ is simple, so that is
  omitted.
\medskip

  The main result of this section is as follows:
\medskip

  {\bf Theorem 2.1}\ \ {\em Let the assumptions $(B_1)$--$(B_4)$ be satisfied.
  Then there exists a neighborhood ${\mathcal O}_2$ of $u_s$, ${\mathcal O}_2
  \subseteq {\mathcal O}_1$, such that the following assertions hold:

  $(1)$\ \ For any $u_0\in {\mathcal O}_2$ the problem (2.1) has a unique
  solution $u(t,u_0)$ which exists for all $t\geq 0$, and if furthermore
  $F(u_0)\in\bar{X}_0$, then $u(t,u_0)$ is a strict solution.

  $(2)$\ \ The center manifold of the equation $u'=F(u)$ in ${\mathcal O}_2$ is
  given by ${\mathcal M}^c=\{S_\sigma(u_s):\sigma\in G\}\cap {\mathcal O}_2$,
  which is a $C^{2-0}$ manifold of dimension $n$ and consists of all stationary
  solutions of this equation in ${\mathcal O}_2$.

  $(3)$\ \ There exists a $C^{2-0}$ submanifold ${\mathcal M}^s\subseteq {\mathcal
  O}_2$ of codimension $n$ in $X_0$ passing $u_s$, such that for any $u_0\in
  {\mathcal M}^s$ there holds $\lim_{t\to\infty}u(t,u_0)=u_s$ and vice versa,
  i.e. ${\mathcal M}^s$ is the stable manifold of $u_s$ in ${\mathcal O}_2$.

  $(4)$\ \ For every $u_0\in {\mathcal O}_2$ there exist a unique $\sigma\in G$
  and a unique $v_0\in {\mathcal M}^s$ such that $u_0=S_\sigma(v_0)$, and we
  have
$$
  \lim_{t\to\infty}u(t,u_0)=S_\sigma(u_s).
\eqno{(2.3)}
$$
  Moreover, for any $0<\omega<\omega_-$ there exists corresponding $C=C(\omega)
  >0$ such that
$$
  \|u(t,u_0)-S_\sigma(u_s)\|_{X_0}\leq Ce^{-\omega t}\|u_0-S_\sigma(u_s)\|_{X_0}
  \quad \mbox{for all}\;\; t\geq 0.
\eqno{(2.4)}
$$
}

  To prove this theorem, we need a preliminary lemma. Let $X$ be a Banach space.
  Let $\alpha\in (0,1)$ and $T>0$. Recall that $C^\alpha_\alpha((0,T],X)$ is
  the Banach space of bounded mappings $u:(0,T]\to X$ such that $t^\alpha u(t)$
  is uniformly $\alpha$-H\"{o}lder continuous for $0<t\leq T$, with norm
$$
  \|u\|_{C^\alpha_\alpha((0,T],X)}=\sup_{0<t\leq T}\|u(t)\|_X+
  \sup_{0<s<t\leq T}\frac{\|t^\alpha u(t)-s^\alpha u(s)\|_X}{(t-s)^\alpha}.
$$
  For $\omega>0$, $C^\alpha([T,\infty),X,-\omega)$ is the Banach space of
  bounded mappings $u:[T,\infty)\to X$ such that $e^{\omega t}u(t)$ is uniformly
  $\alpha$-H\"{o}lder continuous for $t\geq T$, with norm
$$
  \|u\|_{C^\alpha([T,\infty),X,-\omega)}=\sup_{t\geq T}\|e^{\omega t} u(t)\|_X
  +\sup_{t>s\geq T}\frac{\|e^{\omega t}u(t) -e^{\omega s}u(s)\|_X}{(t-s)^\alpha}.
$$
\medskip

  {\bf Lemma 2.2}\ \ {\em Let $X$ and $X_0$ be two Banach spaces such that $X_0
  \hookrightarrow X$. Let $A$ be a sectorial operator
  in $X$ with domain $X_0$. Assume that $\omega_-=-\sup\{{\rm Re}\lambda:
  \lambda\in\sigma(A)\}>0$ and $f\in C^\alpha_\alpha((0,1],X)\cap C^\alpha([1,
  \infty),X,-\omega)$, where $\alpha\in (0,1)$ and $\omega\in (0,\omega_-)$.
  Let $u(t)=e^{tA}u_0+\displaystyle\int_0^t e^{(t-s)A}f(s)ds$, where $u_0\in
  X_0$. Then $u\in C^\alpha_\alpha((0,1],X_0)\cap C^\alpha([1,\infty),X_0,
  -\omega)$, and there exists constant $C=C(\alpha,\omega)>0$ independent of
  $f$ and $u_0$ such that
$$
  \|u\|_{C^\alpha_\alpha((0,1],X_0)}+\|u\|_{C^\alpha([1,\infty),X_0,-\omega)}
  \leq C(\|u_0\|_{X_0}+\|f\|_{C^\alpha_\alpha((0,1],X)}
  +\|f\|_{C^\alpha([1,\infty),X,-\omega)}).
\eqno{(2.5)}
$$
}

  {\em Proof}:\ \ By Theorem 4.3.5 and Corollary 4.3.6 (ii) of \cite{Lunar} we
  have $\|u\|_{C^\alpha_\alpha((0,1],X_0)}\leq C(\|u_0\|_{X_0}+
  \|f\|_{C^\alpha_\alpha((0,1],X)})$, and by Proposition 4.4.10 (i) of
  \cite{Lunar} we have $\|u\|_{C^\alpha([1,\infty),X_0,-\omega)}\leq
  C(\|u_0\|_X+\|f\|_{L^1([0,1],X)}+\|f\|_{C^\alpha([{1\over 2},\infty),X,
  -\omega)}$. Hence (2.5) holds. $\qquad\Box$
\medskip

  {\em Proof of Theorem 2.1}:\ \ Without loss of generality we assume that
  $u_s=0$. Since we are studying solutions of (2.1) in a neighborhood of $0$,
  by the assumption $(B_1)$ and a standard perturbation result, we may assume
  that $F'(u)$ is a sectorial operator for every $u\in {\mathcal O}$ (with
  domain $X_0$), and the graph norm of $F'(u)$ is equivalent to the norm of
  $X_0$. It follows by a standard result (cf. Theorem 8.1.1 of \cite{Lunar}
  and the remark in Lines 8--12 on Page 341 of \cite{Lunar}) that for any
  $u_0\in {\mathcal O}$, the problem (2.1) has a unique local solution $u\in
  C([0,T],X)\cap
  C((0,T],{\mathcal O})\cap L^\infty((0,T),{\mathcal O})\cap C^1((0,T],X)\cap
  C^\alpha_\alpha((0,T],X_0)$, and if further $F(u_0)\in\bar{X}_0$ then $u\in
  C([0,T],{\mathcal O})\cap C^1([0,T],X)\cap C^\alpha_\alpha((0,T],X_0)$,
  where $T>0$ depends on $u_0$ and $\alpha$ is an arbitrary number in $(0,1)$.
  Moreover, denoting by $T^\ast(u_0)$ the supreme of all such $T$, we know
  that there exists a constant $\epsln>0$ independent of $u_0$ such that if
  $\|u(t,u_0)\|_{X_0}<\epsln$ for all $t\in [0,T^\ast(u_0))$, then $T^\ast
  (u_0)=\infty$ (cf. Proposition 9.1.1 of \cite{Lunar}).

  Next we denote $\sigma_-(A)=\sigma(A)\backslash\{0\}$. Let $\Gamma$ be a
  closed smooth curve in the complex plane which encloses $0$ and separates it
  from $\sigma_-(A)$, and let $P$ be the projection operator in $X$ defined by
$$
  P={1\over 2\pi i}\int_\Gamma R(\lambda, A)d\lambda.
$$
  Since $X={\rm Ker}A\oplus {\rm Range}A$, we have $PX=PX_0={\rm Ker}A$, $(I-P)
  X={\rm Range}A$ (cf. Proposition A.2.2 of \cite{Lunar}), and $AP=0$. Let $A_-
  =(I-P)A|_{(I-P)X_0}:(I-P)X_0\to (I-P)X$. Then $\sigma(A_-)=\sigma(A)
  \backslash\{0\}$, so that $\sup\{{\rm Re}\lambda:\lambda\in\sigma(A_-)\}=-
  \omega_-<0$. Besides, by the assumption $(B_3)$ we see that $A_-:(I-P)X_0\to
  (I-P)X$ is an isomorphism.

  Let ${\mathcal M}^c=\{S_\sigma(0):\sigma\in G\}$. By $(v)$ in the assumption
  $(B_2)$ we see that ${\mathcal M}^c$ is a $C^1$ submanifold of $X$, and $\dim
  {\mathcal M}^c=n$. The equation $u=S_\sigma(0)$ ($\sigma\in G$) gives a
  parametrization of ${\mathcal M}^c$ by $G$. We can also give a parametrization
  of ${\mathcal
  M}^c$ by $PX$ as follows. For $u\in {\mathcal O}$ let $x=Pu$ and $y=(I-P)u$.
  Take two sufficiently small numbers $\delta>0$ and $\delta'>0$ such that $x
  \in B_1(0,\delta)$ and $y\in B_2(0,\delta')$ imply that $u=x+y\in {\mathcal
  O}_1$, where
$$
  B_1(0,\delta)=\{x\in PX:\|x\|_{X_0}<\delta\} \quad \mbox{and} \quad
  B_2(0,\delta')=\{y\in (I-P)X_0:\|y\|_{X_0}<\delta'\}.
$$
  For $(x,y)\in B_1(0,\delta)\times B_2(0,\delta')$ we denote ${\mathcal F}_1
  (x,y)=PF(x+y)$ and ${\mathcal F}_2(x,y)=(I-P)F(x+y)$. We have ${\mathcal F}_2
  \in C^{2-0}(B_1(0,\delta)\times B_2(0,\delta'),(I-P)X)$, ${\mathcal F}_2(0,0)
  =0$, and $D_y{\mathcal F}_2(0,0)=A_-$. Since $A_-:(I-P)X_0\to (I-P)X$ is an
  isomorphism, by the implicit function theorem we infer that if $\delta$ is
  sufficiently small then there exists $\varphi\in C^{2-0}(B_1(0,\delta),B_2(0,
  \delta'))$ such that $\varphi(0)=0$, ${\mathcal F}_2(x,\varphi(x))=0$ for
  every $x\in B_1(0,\delta)$, and for $(x,y)\in B_1(0,\delta)\times B_2(0,
  \delta')$, ${\mathcal F}_2(x,y)=0$ if and only if $y=\varphi(x)$. It follows
  that the equation ${\mathcal F}_2(x,y)=0$ defines a $C^{2-0}$ submanifold
  ${\mathcal M}_0$ of dimension $n$. Since $F(S_\sigma(0))=0$ for every $\sigma
  \in G$, which particularly implies that $(I-P)F(S_\sigma(0))=0$ for every
  $\sigma\in G$, we conclude that ${\mathcal M}^c\cap B_1(0,\delta)\times
  B_2(0,\delta')={\mathcal M}_0$. Hence, the equation $y=\varphi(x)$ gives a
  parametrization of ${\mathcal M}^c$ by $PX$. Furthermore, from this argument
  we also see that ${\mathcal F}_1(x,\varphi(x))=PF(S_\sigma(0))=0$ for every
  $x\in B_1(0,\delta)$.  Note that since $D_x{\mathcal F}_2(0,0)=(I-P)AP=0$, we
  have $\varphi'(0)=-[D_y{\mathcal F}_2(0,0)]^{-1}D_x{\mathcal F}_2(0,0)=0$.

  Let $N(u)=F(u)-Au$ (for $u\in {\mathcal M}_1$), ${\mathcal N}_1(x,y)=PN(x+y)$
  and ${\mathcal N}_2(x,y)=(I-P)N(x+y)$ (for $(x,y)\in B_1(0,\delta)\times B_2
  (0,\delta')$). Let $x_0=Pu_0$ and $y_0=(I-P)u_0$. Then (2.1) is equivalent
  to the following problem:
$$
\left\{
\begin{array}{l}
  x'={\mathcal N}_1(x,y),\;\;\; x(0)=x_0,\\
  y'=A_-y+{\mathcal N}_2(x,y),\;\;\;  y(0)=y_0.
\end{array}
\right.
\eqno{(2.6)}
$$
  Let $(x,y)=(x(t),y(t))$ be the solution of (2.6) defined in a maximal
  interval $[0,T_\ast)$ such that it exists for all $t\in[0,T_\ast)$ and lies
  in $B_1(0,\delta)\times B_2(0,\delta')$. Since $(x,y)=(0,0)$ is a solution
  defined for all $t\geq 0$, by continuous dependence of solutions on initial
  data, we see that there exists a neighborhood ${\mathcal O}_2$ of $0$
  contained in $B_1(0,\delta)\times B_2(0,\delta')$, such that for any $u_0\in
  {\mathcal O}_2$ there holds $T_\ast>1$. In the sequel we assume that $u_0\in
  {\mathcal O}_2$ so that $T_\ast>1$. Let $v(t)=y(t)-\varphi(x(t))$. Since
  $A_-\varphi(x)+{\mathcal N}_2(x,\varphi(x))={\mathcal F}_2(x,\varphi(x))=0$
  and ${\mathcal N}_1(x,\varphi(x))={\mathcal F}_1(x,\varphi(x))=0$ for all
  $x\in B_1(0,\delta)$, we have
$$
\begin{array}{rl}
  v'(t)=& A_-v(t)+[{\mathcal N}_2(x(t),y(t))-
  {\mathcal N}_2(x(t),\varphi(x(t)))]\\
  &-\varphi'(x(t))[{\mathcal N}_1(x(t),y(t))-
  {\mathcal N}_1(x(t),\varphi(x(t)))]
  \equiv A_-v(t)+{\mathcal G}(t),
\end{array}
$$
  so that
$$
  v(t)=e^{tA_-}v(0)+\int_0^t e^{(t-s)A_-}{\mathcal G}(s)ds.
$$
  It follows by Lemma 2.2 that for any $0<\alpha<1$ and $\omega\in
  (0,\omega_-)$ we have
$$
  \|v\|_{C^\alpha_\alpha((0,1],X_0)}+\|v\|_{C^\alpha([1,T_\ast),X_0,-\omega)}
  \leq C(\|v(0)\|_{X_0}+\|{\mathcal G}\|_{C^\alpha_\alpha((0,1],X)}
  +\|{\mathcal G}\|_{C^\alpha([1,T_\ast),X,-\omega)}),
\eqno{(2.7)}
$$
  where $C^\alpha([1,T_\ast),X,-\omega)$ is defined similarly as $C^\alpha([1,
  \infty),X,-\omega)$, with $\infty$ replaced with $T_\ast$. Note that all
  assertions in Lemma 2.2 clearly hold when $\infty$ is replaced by any $T_\ast
  \in (1,\infty]$. By a similar argument as in the proof of Theorem 9.1.2 (more
  precisely, as in Line 24, Page 342 through Line 10, Page 343) of \cite{Lunar},
  we have
$$
  \|{\mathcal G}\|_{C^\alpha_\alpha((0,1],X)}\leq
  C(\sup_{0<t\leq 1}\|u(t)\|_{X_0}+
  \sup_{0<t\leq 1}\|\widetilde{u}(t)\|_{X_0})
  \|v\|_{C^\alpha_\alpha((0,1],X)}
\eqno{(2.8)}
$$
  and
$$
  \|{\mathcal G}\|_{C^\alpha([1,T_\ast),X,-\omega)}\leq
  C(\sup_{0\leq t<T_\ast}\|u(t)\|_{X_0}+
  \sup_{0\leq t<T_\ast}\|\widetilde{u}(t)\|_{X_0})
  \|v\|_{C^\alpha([1,T_\ast),X,-\omega)},
\eqno{(2.9)}
$$
  where $\widetilde{u}(t)=x(t)+\varphi(x(t))$, and $C$ is a constant independent
  of $T_\ast$. Substituting (2.8), (2.9) into (2.7), we obtain
$$
  \|v\|_{C^\alpha_\alpha((0,1],X_0)}+\|v\|_{C^\alpha([1,T_\ast),X_0,-\omega)}
  \leq C[\|v(0)\|_{X_0}+(\delta+\delta')(\|v\|_{C^\alpha_\alpha((0,1],X_0)}
  +\|v\|_{C^\alpha([1,T_\ast),X_0,-\omega)})]
$$
  Thus, if $\delta$ and $\delta'$ are sufficiently small then we have
$$
  \|v\|_{C^\alpha_\alpha((0,1],X_0)}+\|v\|_{C^\alpha([1,T_\ast),X_0,-\omega)}
  \leq C\|v(0)\|_{X_0},
$$
  which implies, in particular, that
$$
  \|v(t)\|_{X_0}\leq Ce^{-\omega t}\|v(0)\|_{X_0} \quad \mbox{for}\;\;
  0\leq t<T_\ast,
\eqno{(2.10)}
$$
  where $C$ is independent of $T_\ast$. Next, since ${\mathcal N}_1(x,\varphi
  (x))=0$, we have
$$
  x'(t)={\mathcal N}_1(x(t),y(t))-{\mathcal N}_1(x(t),\varphi(x(t)))
  \equiv {\mathcal G}_1(t).
$$
  It can be easily shown that
$$
  \|{\mathcal G}_1(t)\|_X\leq C(\|u(t)\|_{X_0}+\|\widetilde{u}(t)\|_{X_0})
  \|v(t)\|_{X_0}.
$$
  Hence
$$
  \|x(t)\|_{X_0}\leq C\|x(t)\|_X\leq C\int_0^t\|{\mathcal G}_1(s)\|_Xds
  \leq C(\delta+\delta')\int_0^t\|v(s)\|_{X_0}ds\leq C\|v(0)\|_{X_0}.
\eqno{(2.11)}
$$
  Now, since $u(t)=x(t)+v(t)+\varphi(x(t))$ and $y(t)=v(t)+\varphi(x(t))$, by
  using (2.10) and (2.11) we can easily deduce that if ${\mathcal O}_2$ is
  sufficiently small then for any $u_0\in {\mathcal O}_2$ we have $T_\ast=
  T^\ast(u_0)=\infty$. This proves the assertion (1).

  Similarly as in the proof of (2.11), for any $s>t\geq 0$ we have
$$
  \|x(t)-x(s)\|_{X_0}\leq C\!\int_t^s\!\|{\mathcal G}_1(\tau)\|_X d\tau
  \leq C\!\int_t^s\!\|v(\tau)\|_{X_0}d\tau\leq
  C(e^{-\omega t}\!-\!e^{-\omega s})\|v(0)\|_{X_0}.
\eqno{(2.12)}
$$
  Hence $\lim_{t\to\infty}x(t)$ exists. Let $\bar{x}=\lim_{t\to\infty}x(t)$
  and $\bar{u}=\bar{x}+\varphi(\bar{x})$. Then $\bar{u}\in {\mathcal M}^c$, so
  that it is a stationary point of the equation $u'=F(u)$. Moreover, by the
  facts that $\lim_{t\to\infty}x(t)=\bar{x}$ and $\lim_{t\to\infty}v(t)=0$ in
  $(I-P)X_0$ we see that $\lim_{t\to\infty}u(t)=\bar{u}$ in $X_0$. Letting
  $s\to\infty$ in (2.12) we see that
$$
  \|x(t)-\bar{x}\|_{X_0}\leq Ce^{-\omega t}\|v(0)\|_{X_0}.
\eqno{(2.13)}
$$
  From (2.10) and (2.13) we obtain
$$
  \|u(t)-\bar{u}\|_{X_0}\leq Ce^{-\omega t}\|v(0)\|_{X_0}.
\eqno{(2.14)}
$$
  Hence ${\mathcal M}^c$ is the unique center manifold of the equation $u'=
  F(u)$ in a neighborhood of the origin. This proves the assertion (2).

  Next we note that the equation $u'=F(u)$ can be rewritten as $u'=Au+N(u)$.
  Besides, it is clear that $N\in C^{2-0}({\mathcal O}_1,X)$, $N(0)=0$, and
  $N'(0)=0$, so that $\|N(u)\|_X\leq C\|u\|_{X_0}^2$ and $\|N(u)-N(v)\|_X\leq
  (\|u\|_{X_0}+\|v\|_{X_0})\|u-v\|_{X_0}$. Given $y\in B_2(0,\delta')$, we
  consider the initial value problem
$$
  u'(t)=Au(t)+N(u(t))\quad \mbox{for}\;\; t>0, \quad \mbox{and} \quad
  (I-P)u(0)=y.
\eqno{(2.15)}
$$
  We assert that this problem has a unique solution defined for all $t\geq 0$
  and converging to $0$ as $t\to\infty$, provided $\delta'$ is sufficiently
  small. To prove existence let $\alpha$ and $\omega$ be as before, and for a
  positive number $R$ to be specified later we introduce a metric space
  $(M^\alpha_\omega(R),d)$ by letting
$$
  M^\alpha_\omega(R)=\{u\in C([0,\infty),X_0)\bigcap C^\alpha_\alpha((0,1],X_0)
  \bigcap C^\alpha([1,\infty),X_0,-\omega):|\|u\||\leq R\},
$$
  where
$$
  |\|u\||=\|u\|_{C^\alpha_\alpha((0,1],X_0)}
  +\|u\|_{C^\alpha([1,\infty),X_0,-\omega)},
$$
  and $d(u,v)=|\|u-v\||$. We define a mapping $\Psi_y:M^\alpha_\omega(R)\to
  C([0,\infty),X_0)$ by letting $\Psi_y(u)=v$ for every $u\in M^\alpha_\omega
  (R)$, where
$$
  v(t)=e^{tA_-}y+\int_0^t e^{(t-s)A_-}(I-P)N(u(s))ds
  -\int_t^\infty PN(u(s))ds.
$$
  Using Lemma 2.2, we can easily prove that for sufficiently small $R$,
  $\delta'$ and for any $y\in B_2(0,\delta')$,
  $\Psi_y$ is well-defined, maps $M^\alpha_\omega(R)$ into itself and is a
  contraction mapping. Hence, $\Psi_y$ has a unique fixed point in
  $M^\alpha_\omega(R)$ which we denote by $u_y$. Since $AP=0$ so that
  $e^{(t-s)A}P=P$, it is clear that $u_y$ is a solution of (2.15), and
  $\lim_{t\to\infty}\|u_y(t)\|_{X_0}=0$. This proves existence. To prove
  uniqueness, for any $(x,y)\in B_1(0,\delta)\times B_2(0,\delta')$ we denote
  by $u(t,x,y)$ the unique solution of the equation $u'=F(u)$ satisfying the
  initial conditions $Pu(0)=x$ and $(I-P)u(0)=y$. By Assertion (1) we know that
  $u(t,x,y)$ exists for all $t\geq 0$. Using the fact $AP=0$ we can easily
  deduce that $\lim_{t\to\infty}u(t,x,y)=0$ if and only if
$$
  x+\int_0^\infty PN(u(s,x,y))ds=0.
\eqno{(2.16)}
$$
  We introduce a mapping ${\mathcal F}:B(0,\delta)\times B(0,\delta')\to PX$ by
  letting
$$
  {\mathcal F}(x,y)=x+\int_0^\infty PN(u(s,x,y))ds
$$
  for  $(x,y)\in B_1(0,\delta)\times B_2(0,\delta')$. ${\mathcal F}$ is
  well-defined. Indeed, we know that for any $(x,y)\in B_1(0,\delta)\times
  B_2(0,\delta')$, $\bar{u}=\lim_{t\to\infty}u(t,x,y)$ exists and it belongs to
  ${\mathcal M}^c$. Let $\bar{x}=P\bar{u}$ and $\bar{y}=(I-P)\bar{u}$. Then
  $\bar{y}=\varphi(\bar{x})$, so that $PN(\bar{u})={\mathcal F}_1(\bar{x},
  \bar{y})={\mathcal F}_1(\bar{x},\varphi(\bar{x}))=0$. Thus we have
$$
  \|PN(u(s,x,y))\|_X=\|PN(u(s,x,y))-PN(\bar{u})\|_X\leq
  C\|u(s,x,y)-\bar{u}\|_{X_0}\leq C(x,y)e^{-\omega s}.
$$
  Hence, the integral in the definition of ${\mathcal F}$ is convergent. By a
  similar argument we can show that ${\mathcal F}\in C^{2-0}(B(0,\delta)\times
  B(0,\delta'),PX)$. Since $u(t,0,0)=0$ and $N(0)=N'(0)=0$, we have ${\mathcal
  F}(0,0)=0$ and $D_x{\mathcal F}(0,0)=id$. Thus, by the implicit function
  theorem we conclude that the solution of (2.16) is unique for fixed $y\in
  B_2(0,\delta')$, provided $\delta'$ is sufficiently small. This proves
  uniqueness.

  We now introduce a mapping $\psi: (I-P)X_0\to PX$ by define
$$
  \psi(y)=Pu_y(0)=-\int_0^\infty PN(u_y(s))ds \quad \mbox{for}\;\;
  y\in B_2(0,\delta').
$$
  Clearly, $x=\psi(y)$ is the implicit function solving the equation ${\mathcal
  F}(x,y)=0$, so that $\psi\in C^{2-0}(B_2(0,\delta'),PX)$. Letting ${\mathcal
  M}^s={\rm graph}\psi$, we see that all requirements of the assertion (3) are
  satisfied. This proves the assertion (3).

  Finally, given $u_0\in {\mathcal O}_3$ let $\bar{u}$ be as in (2.14). Since
  $\bar{u}\in {\mathcal M}^c$, there exists a unique $\sigma\in G$ such that
  $S_\sigma(0)=\bar{u}$. Let $v_0=S_{\sigma^{-1}}(u_0)=S_{\sigma}^{-1}(u_0)$.
  Then we have
$$
  \lim_{t\to\infty}u(t,v_0)=\lim_{t\to\infty}S_{\sigma}^{-1}(u(t,u_0))=
  S_{\sigma}^{-1}(S_\sigma(0))=0,
$$
  so that $v_0\in {\mathcal M}^s$. Noticing that $u_0=S_{\sigma}(v_0)$ and
  (2.4) is an immediate consequence of (2.14), we get the assertion (4). This
  completes the proof. $\qquad\Box$
\medskip

  {\bf Remark 2.1}.\ \ Checking the proof of Theorem 2.1, we see that the
  condition on the Lie group action $p$ can be weakened, that is, $p$ need not
  to act on the space $X$; an action on $X_0$ is sufficient.

\section{Reduction of the problem}

  In this section we shall reduce the problem (1.1)--(1.7) into an initial
  value problem of an abstract differential equation in some Banach space. The
  reduction will be fulfilled in two steps: First we use the Hanzawa
  transformation to convert the free boundary problem (1.1)--(1.7) into an
  initial-boundary value problem on the fixed domain $\Omega_s$. Next we solve
  the equations for $p$ in terms of $\sigma$ and $\rho$, the function defining
  the free boundary $\partial\Omega(t)$, to reduce this initial-boundary value
  problem into a purely evolutionary type and regard it as a differential
  equation in a suitable Banach space, which will be the desired abstract
  equation. For simplicity of notation, later on we always assume that $R_s=1$.
  Note that this assumption is reasonable because the general case can be
  reduced into this special case by making suitable rescaling. It follows that
$$
  \Omega_s={\mathbb B}^n=\{x\in {\mathbb R}^n:|x|<1\} \quad
  \mbox{and} \quad \partial\Omega_s=\partial{\mathbb B}^n={\mathbb S}^{n-1}.
$$
  Besides, throughout this paper we assume that the initial domain $\Omega_0$
  is a small perturbation of $\Omega_s={\mathbb B}^n$, so that $\partial\Omega_0$
  is contained in a small neighborhood of $\partial\Omega_s={\mathbb S}^{n-1}$.

  To perform the first step of reduction let us first consider the Hanzawa
  transformation.

  Fix a positive number $\delta$ such that $0<\delta<1$, and denote
$$
  {\mathcal O}_\delta ({\mathbb S}^{n-1})=\{\,\rho\in C^1({\mathbb S}^{n-1}):\;
  \|\rho\|_{C^1({\mathbb S}^{n-1})}<\delta\}.
$$
  Given $\rho\in {\mathcal O}_\delta ({\mathbb S}^{n-1})$, we define a mapping
  $\theta_\rho:\;{\mathbb S}^{n-1}\to {\mathbb R}^n$ by letting $\theta_\rho
  (\xi)=(1+\rho(\xi))\xi$ for $\xi\in {\mathbb S}^{n-1}$, and denote
$$
  \Gamma_\rho={\rm Im}(\theta_\rho)=\{\,x\in {\mathbb R}^n:\;\,
  x=(1+\rho(\xi))\xi,\;\,\xi\in {\mathbb S}^{n-1}\}.
$$
  Clearly, $\Gamma_\rho$ is a closed $C^1$-hypersurface diffeomorphic to
  ${\mathbb S}^{n-1}$, and $\theta_\rho$ is a $C^1$-diffeomorphism from
  ${\mathbb S}^{n-1}$ onto $\Gamma_\rho$. We denote by $\Omega_\rho$ the
  domain enclosed by $\Gamma_\rho$. In the following we always assume that
  $\partial\Omega_0$ is of $C^1$ class and is contained in the
  $\delta$-neighborhood of ${\mathbb S}^{n-1}$. More precisely, we assume
  that there exists $\rho_0\in {\mathcal O}_\delta({\mathbb S}^{n-1})$ such
  that $\partial\Omega_0=\Gamma_{\rho_0}$, and, accordingly, $\Omega_0=
  \Omega_{\rho_0}$.

  Let $m$ be an integer, $m\geq 2$, and let $n/(m-1)<q<\infty$. Then we
  have $B^{m-1/q}_{qq}({\mathbb S}^{n-1})\subseteq C^1({\mathbb S}^{n-1})$.
  The well-known trace theorem ensures that the trace
  operator $\tr(u)=u|_{{\mathbb S}^{n-1}}$ from $C^\infty(\overline{{\mathbb
  B}}^n)$ to $C^\infty({\mathbb S}^{n-1})$ can be extended to $W^{m,q}
  ({\mathbb B}^n)$ such that it maps $W^{m,q}({\mathbb B}^n)$ into
  $B^{m-1/q}_{qq}({\mathbb S}^{n-1})$ and is bounded and surjective. We
  introduce a right inverse $\Pi$ of this operator as follows: Given $\rho\in
  B^{m-1/q}_{qq}({\mathbb S}^{n-1})$, let $u\in W^{m,q}({\mathbb B}^n)$ be the
  unique solution of the boundary value problem
$$
  \Delta u=0 \quad \mbox{in}\;\;{\mathbb B}^n, \quad \mbox{and} \quad
  u=\rho \quad \mbox{on}\;\;{\mathbb S}^{n-1},
$$
  and define $\Pi(\rho)=u$. Then clearly $\tr(\Pi(\rho))=\rho$ for $\rho\in
  B^{m-1/q}_{qq}({\mathbb S}^{n-1})$, and the standard $L^p$ estimate and the
  maximum principle yield the following relations:
$$
  \|\Pi(\rho)\|_{W^{m,q}({\mathbb B}^n)}\leq C\|\rho\|_{B^{m-1/q}_{qq}
  ({\mathbb S}^{n-1})} \quad \mbox{and} \quad
  \sup_{x\in {\mathbb B}^n}|\Pi(\rho)(x)|
  =\max_{x\in {{\mathbb S}^{n-1}}}|\rho(x)|.
$$
  Note that since $W^{m,q}({\mathbb B}^n)\hookrightarrow C^1(\overline{{\mathbb
  B}^n})$, the first relation implies that
$$
  \|\Pi(\rho)\|_{C^1(\overline{{\mathbb B}^n})}\leq C_0\|\rho\|_{B^{m-1/q}_{qq}
  ({\mathbb S}^{n-1})}.
\eqno{(3.1)}
$$
  Here we use the special notation $C_0$ to denote the constant in (3.1)
  because later on this constant will play a special role. We now introduce
$$
  {\mathcal O}_\delta^{m,q}({\mathbb S}^{n-1})=\{\,\rho\in B^{m-1/q}_{qq}
  ({\mathbb S}^{n-1}):\; \|\rho\|_{B^{m-1/q}_{qq}({\mathbb S}^{n-1})}<\delta,
  \,\|\rho\|_{C^1({\mathbb S}^{n-1})}<\delta\}.
$$

  In the sequel we further assume that $\delta<\min\{1/5,(3C_0)^{-1}\}$. Take
  a function $\phi\in C^\infty({\mathbb R},[0,1])$ such that
$$
  \phi(\tau)=1\;\;\;\mbox{for}\;\; |\tau|\leq\delta, \quad
  \phi(\tau)=0\;\;\;\mbox{for}\;\; |\tau|\geq 3\delta,
  \quad \mbox{and} \quad \sup|\phi'|<{2\over 3}\delta^{-1}.
$$
  Given $\rho\in {\mathcal O}_\delta^{m,q}({\mathbb S}^{n-1})$, we define the
  {\it Hanzawa transformation} $\Theta_\rho:\;\overline{{\mathbb B}^n}\to
  \overline{\Omega}_\rho$ by
$$
  \Theta_\rho(x)=x+\phi(|x|-1)\Pi(\rho)(x)\omega(x)
  \;\;\;\mbox{for}\;\; x\in \overline{{\mathbb B}^n},
$$
  where $\omega(x)=x/|x|$ for $x\in {\mathbb R}^n\backslash\{0\}$, and $\omega
  (0)=0$. The choice of $\delta$ and the inequality (3.1) ensures that for
  fixed $\omega\in {\mathbb S}^{n-1}$, the function $r\to r+\phi(r-1)\Pi(\rho)
  (r\omega)$ is strictly monotone increasing for $0\leq r\leq 1$, so that
  $\Theta_\rho$ is a bijection from $\overline{{\mathbb B}^n}$ onto
  $\overline{\Omega}_\rho$. In fact, since the derivative of this function is
  strictly positive, it can be easily shown that $\Theta_\rho\in W^{m,q}
  ({\mathbb B}^n,{\Omega}_\rho)$ and $\Theta_\rho^{-1}\in W^{m,q}({\Omega}_\rho,
  {\mathbb B}^n)$. Besides, it is clear that $\Theta_\rho|_{{\mathbb S}^{n-1}}
  =\theta_\rho$. Since $W^{m,q}({\mathbb B}^n)\subseteq C^1(\overline{{\mathbb
  B}}^n)$ and $W^{m,q}({\Omega}_\rho)\subseteq C^1(\overline{{\Omega}}_\rho)$,
  we see that $\Theta_\rho$ is particularly a $C^1$-diffeomorphism from
  $\overline{{\mathbb B}^n}$ onto $\overline{\Omega}_\rho$.

  As usual we denote by $\Theta^\rho_\ast$ and $\Theta_\rho^\ast$ respectively
  the push-forward and pull-back operators induced by $\Theta_\rho$, i.e.,
  $\Theta^\rho_\ast u=u\circ\Theta_\rho^{-1}$ for $u\in C(\overline{{\mathbb
  B}^n})$, and $\Theta_\rho^\ast u=u\circ\Theta_\rho$ for $u\in
  C(\overline{\Omega}_\rho)$. Similarly, $\theta_\rho^\ast$ denotes the
  pull-back operator induced by $\theta_\rho$, i.e.,  $\theta_\rho^\ast u(\xi)
  =u(\theta_\rho(\xi))$ for $u\in C(\Gamma_\rho)$ and $\xi\in{\mathbb S}^{n-1}$.
  Later, we shall need the following result:
\medskip

  {\bf Lemma 3.1}\ \ {\em Let $m$ be an integer and $1\leq q<\infty$. Let
  $\Omega_1$ and $\Omega_2$ be two open subsets of ${\mathbb R}^n$. Let
  $\Phi$ be a diffeomorphism from ${\Omega}_1$ to ${\Omega}_2$ such that
  $\Phi\in W^{m,q}(\Omega_1,{\mathbb R}^n)$ and $\Phi^{-1}\in W^{m,q}
  (\Omega_2,{\mathbb R}^n)$. Assume that $m\geq 2$ and $q>n/(m-1)$. Then for
  any $0\leq k\leq m$ we have
$$
  \Phi_\ast\in L(W^{k,q}(\Omega_1),W^{k,q}(\Omega_2))
  \quad \mbox{and} \quad
  \Phi^\ast\in L(W^{k,q}(\Omega_2),W^{k,q}(\Omega_1)).
$$
  In particular, for any $\rho\in B^{m-1/q}_{qq}({\mathbb S}^{n-1})$ and
  $0\leq k\leq m$ we have}
$$
  \Theta^\rho_\ast\in L(W^{k,q}({\mathbb B}^n),W^{k,q}(\Omega_\rho))
  \quad \mbox{and} \quad
  \Theta_\rho^\ast\in L(W^{k,q}(\Omega_\rho),W^{k,q}({\mathbb B}^n)).
$$

  {\em Proof}:\ \ See the proof of Lemma 2.1 of \cite{Cui6} for the case $k=m$.
  Proofs for the rest cases $0\leq k\leq m-1$ are similar and simpler.
  $\quad\Box$
\medskip

  Next we introduce some notations.

  In the sequel we assume that $m\geq 2$ and $q>n/(m-1)$. As in \cite{Cui6},
  for $\rho\in {\mathcal O}_\delta^{m,q}({\mathbb S}^{n-1})$ we introduce a
  second-order partial differential operator ${\mathcal A}(\rho):W^{m,q}
  ({\mathbb B}^n)\to W^{m-2,q}({\mathbb B}^n)$ by
$$
  {\mathcal A}(\rho)u=\Theta_\rho^\ast\Delta(\Theta^\rho_\ast u) \quad
  \mbox{for} \quad u\in W^{m,q}({\mathbb B}^n).
$$
  By Lemma 3.1 we see that ${\mathcal A}(\rho)\in L(W^{m,q}({\mathbb B}^n),
  W^{m-2,q}({\mathbb B}^n))$. We also introduce nonlinear operators ${\mathcal
  F}$ and ${\mathcal G}$: $W^{m,q}({\mathbb B}^n)\to W^{m,q}({\mathbb B}^n)$
  respectively by
$$
  {\mathcal F}(u)=f\circ u,\quad {\mathcal G}(u)=g\circ u \quad
  \mbox{for} \quad u\in W^{m,q}({\mathbb B}^n).
$$
  Since the condition $q>n/(m-1)>n/m$ implies that $W^{m,q}({\mathbb B}^n)$ is
  an algebra, we see that these definitions make sense and we have ${\mathcal
  F},{\mathcal G}\in C^\infty(W^{m,q}({\mathbb B}^n), W^{m,q}({\mathbb B}^n))$.
  Given $\rho\in {\mathcal O}_\delta^{m,q}({\mathbb S}^{n-1})$ we denote
$$
  \psi_\rho(x)=|x|-1-\rho(\omega(x)) \quad \mbox{for}\;\; x\in
  {\mathcal R}\equiv\{x\in {\mathbb R}^n:\; 1-4\delta<|x|<1+4\delta\}.
$$
  Clearly, $\psi_\rho\in B^{m-1/q}_{qq}({\mathcal R})$.
  Since $\Gamma_\rho=\{x\in {\mathcal R}:\;\psi_\rho(x)
  =0\}$, we see that the unit outward normal field ${\mathbf n}$ on
  $\Gamma_\rho$ is given by ${\mathbf n}(x)=\nabla\psi_\rho(x)/|\nabla
  \psi_\rho(x)|$ for $x\in\Gamma_\rho$. We introduce a first-order trace
  operator ${\mathcal D}(\rho): W^{m,q}({\mathbb B}^n)\to B^{m-1-1/q}_{qq}
  ({\mathbb S}^{n-1})$ by
$$
  {\mathcal D}(\rho)u=\theta_\rho^\ast({\rm tr}_{\Gamma_\rho}
  (\nabla(\Theta^\rho_\ast u)\cdot\nabla\psi_\rho)) \quad
  \mbox{for} \quad u\in W^{m,q}({\mathbb B}^n),
$$
  where ${\rm tr}_{\Gamma_\rho}$ denotes the usual trace operator from
  $\overline{\Omega}_\rho\cap {\mathcal R}$ to $\Gamma_\rho$, i.e., ${\rm
  tr}_{\Gamma_\rho}(u)=u|_{\Gamma_\rho}$ for $u\in C(\overline{\Omega}_\rho\cap
  {\mathcal R})$. It can be easily seen that ${\mathcal D}(\rho)$ maps $W^{m,q}
  ({\mathbb B}^n)$ into $B^{m-1-1/q}_{qq}({\mathbb S}^{n-1})$, and
  ${\mathcal D}(\rho)\in L(W^{m,q}({\mathbb B}^n), B^{m-1-1/q}_{qq}
  ({\mathbb S}^{n-1}))$ for any $\rho\in {\mathcal O}_\delta^{m,q}
  ({\mathbb S}^{n-1})$. Similarly, given $(\rho,u)\in
  {\mathcal O}_\delta^{m,q}({\mathbb S}^{n-1})\times W^{m,q}({\mathbb B}^n)$,
  we introduce a first-order pseudo-differential operator ${\mathcal P}(\rho,
  u):\; W^{m,q}({\mathbb B}^n)\to W^{m-1,q}({\mathbb B}^n)$ as follows:
$$
  {\mathcal P}(\rho,u)v={\mathcal M}(\rho,u)\cdot\Pi({\mathcal D}(\rho)v)
  \quad \mbox{for} \quad v\in W^{m,q}({\mathbb B}^n).
$$
  Here we used the same notation $\Pi$ as before to denote the bounded right
  inverse of the trace operator $\tr:\; W^{m-1,q}({\mathbb B}^n)\to
  B^{m-1-1/q}_{qq}({\mathbb S}^{n-1})$ such that its restriction on
  $B^{m-1/q}_{qq}({\mathbb S}^{n-1})$ is equal to the previous $\Pi$, and
$$
   {\mathcal M}(\rho,u)(x)=\phi(|x|-1)\langle(\Theta_\rho^\ast
   \nabla\Theta^\rho_\ast u)(x),\omega(x)\rangle\quad \mbox{for}\;\;
   x\in {\mathbb B}^n,
$$
  where $\langle\cdot,\cdot\rangle$ denotes the inner product in ${\mathbb
  R}^n$. We note that ${\mathcal M}(\rho,u)\in W^{m-1,q}({\mathbb B}^n)$ and
  the mapping $u\to {\mathcal M}(\rho,u)$ is a first-order partial differential
  operator. Since $[v\to \Pi({\mathcal D}(\rho)v)]\in L(W^{m,q}({\mathbb B}^n),
  W^{m-1,q}({\mathbb B}^n))$ and the condition $q>n/(m-1)$ implies that
  $W^{m-1,q}({\mathbb B}^n)$ is an algebra, we see that ${\mathcal P}(\rho,u)
  \in L(W^{m,q}({\mathbb B}^n),W^{m-1,q}({\mathbb B}^n))$. Finally, we define
  the transformed mean curvature operator ${\mathcal K}$: $C^2({\mathbb
  S}^{n-1})\cap {\mathcal O}_\delta({\mathbb S}^{n-1})\to C({\mathbb
  S}^{n-1})$ by
$$
  {\mathcal K}(\rho)=\theta_\rho^\ast(\kappa_{\Gamma_\rho}),
$$
  where $\kappa_{\Gamma_\rho}$ denotes the mean curvature of the hypersurface
  $\Gamma_\rho$ (recall that $\kappa_{\Gamma_\rho}\in C(\Gamma_\rho, {\mathbb
  R})$ for $C^2$ class hypersurface $\Gamma_\rho$).
  Later we shall restrict ${\mathcal K}$ in ${\mathcal O}_\delta^{m,q}({\mathbb
  S}^{n-1})$ and shall see that ${\mathcal K}\in C^\infty({\mathcal
  O}_\delta^{m,q}({\mathbb S}^{n-1}),B^{m-2-1/q}_{qq}({\mathbb S}^{n-1}))$.

  Let $T$ be a given positive number and consider a function $\rho:\;[0,T]\to
  {\mathcal O}_\delta^{m,q}({\mathbb S}^{n-1})$. We assume that $\rho\in
  C([0,T],{\mathcal O}_\delta^{m,q}({\mathbb S}^{n-1})$. Given a such $\rho$,
  we denote
$$
  \Gamma_\rho(t)=\Gamma_{\rho(t)}, \quad
  \Omega_\rho(t)=\Omega_{\rho(t)} \quad (0\leq t\leq T).
$$
  Later on in case no confusion can be produced we shall occasionally abbreviate
  $\Gamma_\rho(t)$ and $\Omega_\rho(t)$ respectively as $\Gamma_\rho$ and
  $\Omega_\rho$. We shall briefly write the families of operators $t\to
  {\mathcal A}(\rho(t))$ and $t\to {\mathcal D}(\rho(t))$ ($0\leq t\leq T$)
  as ${\mathcal A}(\rho)$ and ${\mathcal D}(\rho)$, respectively, and for $u,v:
  [0,T]\to W^{m,q}({\mathbb B}^n)$, we briefly write the families of functions
  ${\mathcal F}(\rho(t),u(t))$, ${\mathcal G}(\rho(t),u(t))$ and ${\mathcal M}
  (\rho(t),u(t))v(t)$ ($0\leq t\leq T$) respectively as ${\mathcal F}(\rho,u)$,
  ${\mathcal G}(\rho,u)$ and ${\mathcal M}(\rho,u)v$. Besides, we shall identify
  a function $\rho: [0,T]\to C({\mathbb S}^{n-1})$ (resp. $u: [0,T]\to
  C(\overline{{\mathbb B}^n})$) with the corresponding function on ${\mathbb
  S}^{n-1}\times [0,T]$ (resp. $\overline{{\mathbb B}^n}\times [0,T]$) defined by
  $\rho(\xi,t)=\rho(t)(\xi)$ (resp. $u(x,t)=u(t)(x)$), where $t\in [0,T]$ and
  $\xi\in {\mathbb S}^{n-1}$ (resp. $x\in\overline{{\mathbb B}^n}$), and vice
  versa.

  With the above notations, it is not hard to verify that if we denote
$$
  u(x,t)=\sigma(\Theta_{\rho(t)}(x),t),\quad
  v(x,t)= p(\Theta_{\rho(t)}(x),t),
$$
  then the Hanzawa transformation transforms (1.1)--(1.7) into the following
  system of equations:
$$
  c\partial_t u-{\mathcal A}(\rho)u+c{\mathcal P}(\rho,u)v=-{\mathcal F}(u)
  \qquad\quad \mbox{in}\;\; {\mathbb B}^n\times (0,T],\qquad\quad
\eqno{(3.2)}
$$
$$
  \qquad \qquad\;\;\,\quad -{\mathcal A}(\rho)v={\mathcal G}(u),\quad
  \qquad\;\; \mbox{in}\;\; {\mathbb B}^n\times (0,T],\;\;\,
\eqno{(3.3)}
$$
$$
  \;\;\qquad\qquad\qquad\qquad u=\bar{\sigma}\qquad\qquad\quad\;\;\;
  \mbox{on}\;\; {\mathbb S}^{n-1}\times (0,T],\quad
\eqno{(3.4)}
$$
$$
  \qquad\qquad\qquad\qquad\;\; v=\gamma{\mathcal K}(\rho)\qquad\quad\quad\;\;
  \mbox{on}\;\; {\mathbb S}^{n-1}\times (0,T],\quad
\eqno{(3.5)}
$$
$$
  \;\;\;\qquad\partial_t\rho+{\mathcal D}(\rho)v=0 \qquad\qquad\quad
  \quad\; \mbox{on}\;\; {\mathbb S}^{n-1}\times (0,T],
\eqno{(3.6)}
$$
$$
  \qquad\qquad u(0)=u_0\qquad\qquad\quad\quad\;\;
  \mbox{on}\;\; {\mathbb B}^n,\qquad
\eqno{(3.7)}
$$
$$
  \qquad\qquad\rho(0)=\rho_0\qquad\qquad\quad\quad\;\;
  \mbox{on}\;\; {\mathbb S}^{n-1},\qquad
\eqno{(3.8)}
$$
  where $u_0=\Theta_{\rho_0}^\ast\sigma_0$. Indeed, it is immediate to see that
  (3.3), (3.4), (3.5), (3.7) and (3.8) are transformations of (1.2), (1.3),
  (1.4), (1.6) and (1.7), respectively. For the proof that the transformation
  of (1.5) is (3.6), we refer the reader to see the deduction of (2.19) in
  \cite{Cui6} and (2.8) in \cite{Escher}. Finally, (3.2) is obtained from
  transforming (1.1) and using (3.6).

  To establish properties of the operator ${\mathcal K}$, we need the following
  lemma:
\medskip

  {\bf Lemma 3.2}\ \ {\em $(i)$ Let $k$, $m$ be nonnegative integers, and $p,q
  \in [1,\infty]$. Let $\Omega$ be an open subset of ${\mathbb R}^n$ with a
  smooth boundary. Assume that $k\geq m$ and either $1\leq p\leq n/m$, $k>n/q$
  or $p>n/m$, $k-n/q\geq m-n/p$. Then we have
$$
  \|uv\|_{W^{m,p}(\Omega)}\leq C\|u\|_{W^{k,q}(\Omega)}\|v\|_{W^{m,p}(\Omega)}.
\eqno{(3.9)}
$$

  $(ii)$ Let $s,t>0$ and $p,q,r_1,r_2\in [1,\infty]$. Let $\Omega$
  be as before. Assume that $t\geq s$ and either $1\leq p\leq n/s$, $t>n/q$ or
  $p>n/s$, $t-n/q\geq s-n/p$. Then we have
$$
  \|uv\|_{B^s_{pr_1}(\Omega)}\leq C\|u\|_{B^t_{qr_2}(\Omega)}
  \|v\|_{B^s_{pr_1}(\Omega)}.
\eqno{(3.10)}
$$
  Here $r_1$, $r_2$ are arbitrary numbers in $[1,\infty]$ in case $t>s$, and
  $1\leq r_2\leq r_1\leq\infty$ if $t=s$.}
\medskip

  {\em Proof}:\ \ To prove (3.9), we first note that since $k>n/q$, we have
  $W^{k,q}(\Omega)\hookrightarrow L^\infty(\Omega)$, so that
$$
  \|uv\|_{L^p(\Omega)}\leq\|u\|_{L^\infty(\Omega)}\|v\|_{L^p(\Omega)}
  \leq C\|u\|_{W^{k,q}(\Omega)}\|v\|_{W^{m,p}(\Omega)}.
\eqno{(3.11)}
$$
  Next let $\alpha\in {\mathbb Z}_+^n$ be an arbitrary $n$-index of length $m$,
  i.e., $|\alpha|=m$. We write the Leibnitz formula:
$$
  \partial^\alpha(uv)=\sum_{\beta\leq\alpha}{\alpha!\over\beta!(\alpha-\beta)!}
  \partial^\beta u \;\partial^{\alpha-\beta} v.
$$
  For every $n$-index $\beta\leq\alpha$ we take $r_1,r_2\in [1,\infty]$ as
  follows:
$$
\left\{
\begin{array}{l}
  \displaystyle{1\over r_1}={1\over q}-{k\!-\!|\beta|\over n},\;
  {1\over r_2}={1\over p}-{1\over q}+{k\!-\!|\beta|\over n} \quad
  \mbox{if}\;\; |\beta|>k-{n\over q},\\
  \displaystyle r_1={1\over \epsln},\; {1\over r_2}={1\over p}-\epsln \quad
  \mbox{if}\;\; |\beta|=k-{n\over q},\\
  \displaystyle r_1=\infty,\; r_2=p \quad\mbox{if}\;\; |\beta|<k-{n\over q},
\end{array}
\right.
$$
  where $\epsln$ is a small positive number. Note that since $|\beta|\leq m\leq
  k$, we have ${1\over p}-{1\over q}+{k\!-\!|\beta|\over n}\geq 0$. Clearly
  ${1\over r_1}+{1\over r_2}={1\over p}$, $\|\partial^\beta u\|_{L^{r_1}
  (\Omega)}\leq C\|u\|_{W^{k,q}(\Omega)}$ and $\|\partial^{\alpha-\beta}
  v\|_{L^{r_2}(\Omega)}\leq C\|u\|_{W^{m,p}(\Omega)}$. Hence
$$
  \|\partial^\alpha(uv)\|_{L^p(\Omega)}\leq C\sum_{\beta\leq\alpha}
  \|\partial^\beta u\|_{L^{r_1}(\Omega)}
  \|\partial^{\alpha-\beta}v\|_{L^{r_2}(\Omega)}
  \leq C\|u\|_{W^{k,q}(\Omega)}\|v\|_{W^{m,p}(\Omega)}.
\eqno{(3.12)}
$$
  Combining (3.11) and (3.12), we get (3.9).

  Having proved (3.9), (3.10) easily follows by interpolation. $\qquad\Box$

\medskip

  {\bf Corollary 3.3}\ \ {\em Assume that $m\geq 2$ and either $q>n/(m\!-\!1)$,
  $0<s\leq m\!-\!1\!-\!1/q$ or $q>\max\{\,2n/(m\!+\!n\!-\!2),n/(m-1)\}$ and
  $-1/q\leq s\leq 0$. Then we have
$$
  \|uv\|_{B^s_{qq}({\mathbb S}^{n-1})}\leq
  C\|u\|_{B^{m-1-1/q}_{qq}({\mathbb S}^{n-1})}
  \|v\|_{B^s_{qq}({\mathbb S}^{n-1})}.
\eqno{(3.13)}
$$
}

  {\em Proof}:\ \ If $s>0$ then the desired assertion follows immediately from
  Lemma 3.2 $(ii)$, because we can easily verify that all conditions of Lemma
  3.2 $(ii)$ are satisfied when we replace $t$ with $m-1-1/q$, $p$ with $q$ and
  $n$ with $n-1$. Next we consider the case $-1/q\leq s\leq 0$. We can also
  easily verify that in this case all conditions of Lemma 3.2 $(ii)$ are
  satisfied when we replace $t$ with $m-1-1/q$, $s$ with $1/q$, $p$ with $q'$,
  and $n$ with $n-1$, so that
$$
  \|uv\|_{B^{1/q}_{q'q'}({\mathbb S}^{n-1})}\leq
  C\|u\|_{B^{m-1-1/q}_{qq}({\mathbb S}^{n-1})}
  \|v\|_{B^{1/q}_{q'q'}({\mathbb S}^{n-1})}.
$$
  By dual, this implies that
$$
  \|uv\|_{B^{-1/q}_{qq}({\mathbb S}^{n-1})}\leq
  C\|u\|_{B^{m-1-1/q}_{qq}({\mathbb S}^{n-1})}
  \|v\|_{B^{-1/q}_{qq}({\mathbb S}^{n-1})}.
$$
  Interpolating this inequality with (3.22) for $s>0$, we see that (3.22) also
  holds for $-1/q\leq s\leq 0$ under the prescribed conditions. $\qquad\Box$
\medskip

  {\bf Lemma 3.4}\ \ {\em Let $m\geq 2$ and $q>n/(m-1)$. Then for any $2\leq
  k\leq m$ we have the following assertions:
$$
   {\mathcal A}\in C^\infty({\mathcal O}_\delta^{m,q}({\mathbb S}^{n-1}),
   L(W^{k,q}({\mathbb B}^n),W^{k-2,q}({\mathbb B}^n))),
\eqno{(3.14)}
$$
$$
   {\mathcal D}\in C^\infty({\mathcal O}_\delta^{m,q}({\mathbb S}^{n-1}),
   L(W^{k,q}({\mathbb B}^n),B^{k-1-1/q}_{qq}({\mathbb S}^{n-1}))),
\eqno{(3.15)}
$$
$$
   {\mathcal P}\in C^\infty({\mathcal O}_\delta^{m,q}({\mathbb S}^{n-1})\times
   W^{m,q}({\mathbb B}^n),L(W^{k,q}({\mathbb B}^n), W^{k-1,q}({\mathbb B}^n))),
\eqno{(3.16)}
$$
  and for any $k>n/q$ we have
$$
   {\mathcal F},\;{\mathcal G}\in
   C^\infty(W^{k,q}({\mathbb B}^n),W^{k,q}({\mathbb B}^n)).
\eqno{(3.17)}
$$
}

  {\em Proof}:\ \ (3.14) is an immediate consequence of Lemma 3.1 and the
  fact that $\Delta\in L(W^{k,q}(\Omega_\rho)$, $W^{k-2,q}(\Omega_\rho))$ for
  any $k$. (3.15) is an immediate consequence of Lemmas 3.1, 3.2 and the fact
  that $\nabla\in L(W^{k,q}(\Omega_\rho),W^{k-1,q}(\Omega_\rho,{\mathbb R}^n))$
  for any $k$. (3.16) follows from similar reasons as for (3.15). Finally,
  (3.17) follows from the fact that $W^{k,q}(\Omega_\rho)$ is an algebra under
  the condition $k>n/q$, as we mentioned earlier. $\qquad \Box$
\medskip

  {\bf Lemma 3.5}\ \ {\em $(i)$ The mean curvature operator ${\mathcal K}(\rho)$
  has the following splitting:
$$
  {\mathcal K}(\rho)={\mathcal L}(\rho)\rho+{\mathcal K}_1(\rho),
\eqno{(3.18)}
$$
  where ${\mathcal L}(\rho)$ is a second-order elliptic linear partial
  differential operator on ${\mathbb S}^{n-1}$, with coefficients being
  functions of $\rho$ and its first-order derivatives, and ${\mathcal
  K}_1(\rho)$ is a first-order partial nonlinear differential operator on
  ${\mathbb S}^{n-1}$.

  $(ii)$ Assume that $m\geq 3$ and $q>\max\{\,2n/(m\!+\!n\!-\!2),n/(m-1)\}$.
  Then we have
$$
  {\mathcal L}\in C^\infty({\mathcal O}_\delta^{m,q}({\mathbb S}^{n-1}),
  L(B^{k-1/q}_{qq}({\mathbb S}^{n-1}),B^{k-2-1/q}_{qq}({\mathbb S}^{n-1}))),
  \quad 2\leq k\leq m,
\eqno{(3.19)}
$$
$$
  {\mathcal K}_1\in C^\infty({\mathcal O}_\delta^{m,q}({\mathbb S}^{n-1}),
  B^{m-1-1/q}_{qq}({\mathbb S}^{n-1}))),
\eqno{(3.20)}
$$
  so that
$$
  {\mathcal K}\in C^\infty({\mathcal O}_\delta^{m,q}({\mathbb S}^{n-1}),
  B^{m-2-1/q}_{qq}({\mathbb S}^{n-1})).
\eqno{(3.21)}
$$
}

  {\em Proof}:\ \ The Assertion $(i)$ is an immedaite consequence of the mean
  curvature formula, see \cite{Escher} and \cite{EscSim}. Next, since the
  condition $q>n/(m-1)$ implies that $B^{m-1-1/q}_{qq}({\mathbb S}^{n-1})$ is
  an algebra, (3.20) easily follows from the fact that ${\mathcal K}_1$ is a
  first-order nonlinear partial differential operator. Similarly, (3.19)
  follows from Corollary 3.3 and the facts that $B^{m-1-1/q}_{qq}({\mathbb
  S}^{n-1})$ is an algebra and ${\mathcal L}(\rho)$ is a second-order partial
  differential operator with coefficients being smooth functions of $\rho$ and
  its first-order partial derivatives. Finally, (3.21) follows readily from
  (3.18)--(3.20). $\qquad\Box$
\medskip

  In order to perform the second step of reduction, we need the following lemma:
\medskip

  {\bf Lemma 3.6}\ \ {\em Let $m\geq 2$, $q>n/(m-1)$ and $2\leq k\leq m$. Given
  $\rho\in {\mathcal O}_\delta^{m,q}({\mathbb S}^{n-1})$ and $(w,\eta)\in
  W^{k-2,q}({\mathbb B}^n)\times B^{k-1/q}_{qq}({\mathbb S}^{n-1})$, the
  problem
$$
\left\{
\begin{array}{rl}
   -{\mathcal A}(\rho)u=w \quad & \mbox{in}\;\; {\mathbb B}^n,\\
   u=\eta \quad & \mbox{on}\;\; {\mathbb S}^{n-1}
\end{array}
\right.
$$
  has a unique solution $u\in W^{k,q}({\mathbb B}^n)$, and it has the following
  expression:
$$
  u={\mathcal S}(\rho)w+{\mathcal T}(\rho)\eta,
$$
  where
$$
  {\mathcal S}\in C^\infty({\mathcal O}_\delta^{m,q}({\mathbb S}^{n-1}),
  L(W^{k-2,q}({\mathbb B}^n),W^{k,q}({\mathbb B}^n))),
\eqno{(3.22)}
$$
$$
   {\mathcal T}\in C^\infty({\mathcal O}_\delta^{m,q}({\mathbb S}^{n-1}),
   L(B^{k-1/q}_{qq}({\mathbb S}^{n-1}),W^{k,q}({\mathbb B}^n))).
\eqno{(3.23)}
$$
}

  {\em Proof}:\ \ All assertions easily follow from the standard theory of
  elliptic partial differential equations, cf. the proof of Lemma 3.1 of
  \cite{CuiEsc2}. $\qquad\Box$
\medskip

  In the sequel we perform the second step of reduction.

  By Lemmas 3.5 and 3.6 we see that given $u\in W^{m-1,q}({\mathbb B}^n)$ and
  $\rho\in B^{m-1/q}_{qq}({\mathbb S}^{n-1})$, the solution of Eq. (3.3)
  subject to the boundary value condition (3.5) is given by
$$
  v=\gamma {\mathcal T}(\rho){\mathcal L}(\rho)\rho+
  \gamma {\mathcal T}(\rho){\mathcal K}_1(\rho)+
  {\mathcal S}(\rho){\mathcal G}(u).
$$
  Substitute this expression into (3.2) and (3.6), we see that the problem
  (3.2)--(3.8) is reduced into the following problem:
$$
  \partial_t u-c^{-1}{\mathcal A}(\rho)u-{\mathcal Q}(\rho,u)\rho=
  {\mathcal F}_1(\rho,u)\;\;\quad \mbox{in}\;\;
  {\mathbb B}^n\times (0,\infty),
  \qquad\qquad
\eqno{(3.24)}
$$
$$
  \quad\partial_t\rho-{\mathcal B}(\rho)\rho={\mathcal G}_1(\rho,u)\;\;
  \quad \mbox{on}\;\; {\mathbb S}^{n-1}\times (0,\infty),
\eqno{(3.25)}
$$
$$
  \qquad\qquad\qquad\qquad u=\bar{\sigma}\qquad\quad\;\;
  \mbox{on}\;\; {\mathbb S}^{n-1}\times (0,\infty),\quad\quad\;
\eqno{(3.26)}
$$
$$
  \qquad\qquad u(0)=u_0\qquad\quad\;\;
  \mbox{on}\;\; {\mathbb B}^n,\qquad\;\;\;
\eqno{(3.27)}
$$
$$
  \qquad\qquad\rho(0)=\rho_0\qquad\quad\;\;
  \mbox{on}\;\; {\mathbb S}^{n-1},\qquad\;\;\;
\eqno{(3.28)}
$$
  where ${\mathcal A}(\rho)$ is as before, and
$$
\begin{array}{rl}
  {\mathcal B}(\rho)\zeta &=-\gamma{\mathcal D}(\rho){\mathcal T}(\rho)
  {\mathcal L}(\rho)\zeta,\\
  {\mathcal Q}(\rho,u)\zeta &={\mathcal M}(\rho,u)\cdot
  \Pi({\mathcal B}(\rho)\zeta),\\
  {\mathcal F}_1(\rho,u) &=-c^{-1}{\mathcal F}(u)-\gamma{\mathcal P}(\rho,u)
  {\mathcal T}(\rho){\mathcal K}_1(\rho)-{\mathcal P}(\rho,u){\mathcal S}(\rho)
  {\mathcal G}(u)\\
  &=-c^{-1}{\mathcal F}(u)-{\mathcal M}(\rho,u)\cdot
  \Pi({\mathcal G}_1(\rho,u)),\\
  {\mathcal G}_1(\rho,u) &=-\gamma{\mathcal D}(\rho)
  {\mathcal T}(\rho){\mathcal K}_1(\rho)-{\mathcal D}(\rho){\mathcal S}(\rho)
  {\mathcal G}(u).
\end{array}
$$

  To homogenize the boundary condition (3.26) we define
$$
  {\mathcal C}(\rho,u)={\mathcal Q}(\rho,u+\bar{\sigma}),\quad
  {\mathcal F}_2(\rho,u)={\mathcal F}_1(\rho,u+\bar{\sigma}),\quad
  {\mathcal G}_2(\rho,u)={\mathcal G}_1(\rho,u+\bar{\sigma}).
$$
  Replacing ${\mathcal Q}$, ${\mathcal F}_1$ and ${\mathcal G}_1$ in (3.24) and
  (3.25) with ${\mathcal C}$, ${\mathcal F}_2$ and ${\mathcal G}_2$,
  respectively, we see that the inhomogeneous boundary value condition (3.26)
  is replaced by the homogeneous boundary value condition
$$
  u=0 \qquad \mbox{on}\;\; {\mathbb S}^{n-1}\times (0,\infty).
\eqno{(3.29)}
$$

  We now denote
$$
  U=\left(
\begin{array}{c}
  u\\ \rho
\end{array}
\right), \quad
  {\mathbb A}(U)=\left(
\begin{array}{cc}
  c^{-1}{\mathcal A}(\rho)\;&\; {\mathcal C}(\rho,u)\\
  0 \; &\; {\mathcal B}(\rho)
\end{array}
\right), \quad
  {\mathbb F}_0(U)=\left(
\begin{array}{c}
  {\mathcal F}_2(\rho,u)\\
  {\mathcal G}_2(\rho,u)
\end{array}
\right), \quad
  U_0=\left(
\begin{array}{c}
  \sigma_0-\bar{\sigma}\\ \rho_0
\end{array}
\right),
$$
  and
$$
  {\mathbb F}(U)={\mathbb A}(U)U+{\mathbb F}_0(U).
$$
  We also denote
$$
  X=W^{m-3,q}({\mathbb B}^n)\times B^{m-3-1/q}_{qq}({\mathbb S}^{n-1}),\quad
  X_0=(W^{m-1,q}({\mathbb B}^n)\cap W^{1,q}_0({\mathbb B}^n))\times
  B^{m-1/q}_{qq}({\mathbb S}^{n-1}),
$$
$$
  Y=W^{m-2,q}({\mathbb B}^n)\times B^{m-2-1/q}_{qq}({\mathbb S}^{n-1}),
$$
  and
$$
  {\mathcal O}=(W^{m-1,q}({\mathbb B}^n)\cap W^{1,q}_0({\mathbb B}^n))\times
  {\mathcal O}_\delta^{m,q}({\mathbb S}^{n-1}).
$$
  Then the equations (3.24), (3.25) (with ${\mathcal Q}$, ${\mathcal F}_1$,
  ${\mathcal G}_1$ respectively replaced with ${\mathcal C}$, ${\mathcal F}_2$,
  ${\mathcal G}_2$) and (3.29) are reduced into the following abstract
  differential equation in the Banach space $X$:
$$
  {dU\over dt}={\mathbb F}(U),
\eqno{(3.30)}
$$
  and the problem (3.24)--(3.28) is reduced into the following initial value
  problem:
$$
\left\{
\begin{array}{l}
  U'(t)={\mathbb F}(U(t))\quad\mbox{for}\;\; t>0,\\
  U(0)=U_0
\end{array}
\right.
\eqno{(3.31)}
$$

  Clearly, $X$, $X_0$ and $Y$ are Banach spaces, $X_0\hookrightarrow X$, $Y$ is
  an intermediate space between $X$ and $X_0$, and ${\mathcal O}$ is an open
  subset of $X_0$. From (3.14)--(3.21) and (3.22)--(3.20) we see that
$$
   {\mathbb A}\in C^\infty({\mathcal O},L(X_0,X)),\quad
   {\mathbb F}_0\in C^\infty({\mathcal O},Y)\subseteq C^\infty({\mathcal O},X),
$$
  so that ${\mathbb F}\in C^\infty({\mathcal O},X)$. We note that for $m=3$,
  $X_0$ is dense in $X$, while for $m\geq 4$ the closure of $X_0$ in $X$ is
  given by
$$
  \bar{X}_0=(W^{m-3,q}({\mathbb B}^n)\cap W^{1,q}_0({\mathbb B}^n))\times
  B^{m-3-1/q}_{qq}({\mathbb S}^{n-1}).
$$

\section{The Lie group action}

  For $\epsln>0$ we denote by ${\mathbb B}^n_\epsln$ the ball in ${\mathbb
  R}^n$ centered at the origin with radius $\epsln$. Regarding ${\mathbb
  B}^n_\epsln$ as a neighborhood of the unit element $0$ of the commutative
  Lie group ${\mathbb R}^n$, we see that $G={\mathbb B}^n_\epsln$ is a local
  Lie group of dimension $n$. In this section we introduce an action ${\mathbf
  S}^\ast$ of this (local) Lie group $G$ to some open subset ${\mathcal O}'$
  of $X$, ${\mathcal O}'\cap X_0={\mathcal O}$, such that the relation
$$
  {\mathbb F}({\mathbf S}^\ast_z(u))=D{\mathbf S}^\ast_z(u) {\mathbb F}(u),
  \quad z\in G, \;\; u\in {\mathcal O}
\eqno{(4.1)}
$$
  is satisfied.

  Given $z\in {\mathbb R}^n$, we denote by $S_z$ the translation in ${\mathbb
  R}^n$ induced by $z$, i.e.,
$$
  S_z(x)=x+z \quad \mbox{for}\;\; x\in {\mathbb R}^n.
$$
  Let $\rho\in C^1({\mathbb S}^{n-1})$ such that $\|\rho\|_{C^1({\mathbb
  S}^{n-1})}$ is sufficiently small, say, $\|\rho\|_{C^1({\mathbb S}^{n-1})}
  <\delta$ for some small $\delta>0$. For any $z\in{\mathbb B}^n_\epsln$, where
  $\epsln$ is sufficiently small, consider the image of the
  hypersurface $r=1+\rho(\omega)$ under the translation $S_{z}$, which is still
  a hypersurface. This hypersurface has the equation $r=1+\tilde{\rho}(\omega)$
  with $\tilde{\rho}\in C^1({\mathbb S}^{n-1})$, and $\tilde{\rho}$ is uniquely
  determined by $\rho$ and $z$. We denote
$$
  \tilde{\rho}=S_z^*(\rho).
$$
  Let $r_0=|z|$ and $\omega_0=z/|z|$. Then the explicit expression of
  $\tilde{\rho}$ is as follows:
$$
  \tilde{\rho}(\omega')=\sqrt{[1+\rho(\omega)]^2+r_0^2+2r_0[1+\rho(\omega)]
  \omega\cdot\omega_0}-1,
\eqno{(4.2)}
$$
  where $\omega'\in {\mathbb S}^{n-1}$ and $\omega\in {\mathbb S}^{n-1}$ are
  connected by the following relation:
$$
  \omega'=\frac{[1+\rho(\omega)]\omega+r_0\omega_0}{\sqrt{[1+\rho(\omega)]^2
  +r_0^2+2r_0[1+\rho(\omega)]\omega\cdot\omega_0}}.
\eqno{(4.3)}
$$

  In the sequel, the notations ${\mathcal O}_\delta({\mathbb S}^{n-1})$ and
  ${\mathcal O}_\delta^{m,q}({\mathbb S}^{n-1})$ have same meaning as in the
  previous section.

\medskip
  {\bf Lemma 4.1}\ \ {\em If $\epsln$ and $\delta$ are sufficiently small then
  for any $z\in {\mathbb B}^n_\epsln$ and $\rho\in {\mathcal O}_\delta({\mathbb
  S}^{n-1})$, $S_z^*(\rho)$ is well-defined, and}
$$
  S_z^*\in C({\mathcal O}_\delta({\mathbb S}^{n-1}),C^1({\mathbb S}^{n-1}))
  \cap C^1({\mathcal O}_\delta({\mathbb S}^{n-1}),C({\mathbb S}^{n-1})).
$$

  {\em Proof}:\ \ Let $f_z(\rho,\omega)$ be the expression in the
  right-hand side of (4.3). We first prove that if $\epsln$ is sufficiently
  small then  for any $z\in {\mathbb B}^n_\epsln$ the mapping $\omega\to
  \omega'=f_z(\rho,\omega)$ from ${\mathbb S}^{n-1}$ to itself is an injection.
  Assume that $f_z(\rho,\omega_1)=f_z(\rho,\omega_2)$ for some $\omega_1,
  \omega_2\in {\mathbb S}^{n-1}$. Then there exists $\lambda>0$ such that
$$
  [1+\rho(\omega_2)]\omega_2+r_0\omega_0=
  \lambda\{[1+\rho(\omega_1)]\omega_1+r_0\omega_0\}.
\eqno{(4.4)}
$$
  Let $\lambda=1+\mu$, $\omega_2=\omega_1+\xi$ and $\rho(\omega_2)=\rho
  (\omega_1)+\eta$, where $\mu\in {\mathbb R}$, $\xi\in {\mathbb R}^n$ and
  $\eta\in {\mathbb R}$. Substituting these expressions into (4.4) we get
$$
  [1+\rho(\omega_1)]\xi+\omega_2\eta=
  \mu\{[1+\rho(\omega_1)]\omega_1+r_0\omega_0\},
$$
  which yields $\xi=\mu\omega_1+\zeta$, where $\zeta=(\mu r_0\omega_0-\omega_2
  \eta)/[1+\rho(\omega_1)]$. Since $|\rho(\omega_1)|<\delta$ and $|r_0|<\epsln$,
  from the expression of $\zeta$ we see that $|\zeta|\leq 2(\epsln|\mu|+|\eta|)$
  if $\delta\leq 1/2$. Since $\max_{\omega\in {\mathbb S}^{n-1}}|\nabla_\omega
  \rho(\omega)|<\delta$, by the mean value theorem we easily deduce that
  $|\eta|\leq\delta|\xi|$, so that
$$
  |\zeta|\leq 2(\epsln|\mu|+\delta|\xi|).
\eqno{(4.5)}
$$
  From the relation $\xi=\mu\omega_1+\zeta$ we have
$$
  |\xi|\leq |\mu|+|\zeta|.
\eqno{(4.6)}
$$
  Substituting the relation $\xi=\mu\omega_1+\zeta$ into $\omega_2=\omega_1+
  \xi$ we get $\omega_2=(1+\mu)\omega_1+\zeta$, or $(1+\mu)\omega_1=\omega_2
  -\zeta$. From this relation and the fact that $|\mu|<1$ (for $\epsln$ and
  $\delta$ sufficiently small) we obtain
$$
  |\mu|\leq |\zeta|.
\eqno{(4.7)}
$$
  From (4.5)--(4.7) we can easily deduce that $|\zeta|=|\xi|=|\mu|=0$ for
  sufficiently small $\epsln$ and $\delta$, which proves the desired assertion.

  Next we prove that if $\epsln$ and $\delta$ are sufficiently small then
  $D_\omega f_z(\rho,\omega): T_\omega({\mathbb S}^{n-1})\to T_\omega({\mathbb
  S}^{n-1})$ is non-degenerate for any $\omega\in {\mathbb S}^{n-1}$ and
  $\rho\in {\mathcal O}_\delta({\mathbb S}^{n-1})$. Note that since $\rho
  \in C^1({\mathbb S}^{n-1})$, we have $f_z(\rho,\cdot)\in C^1({\mathbb
  S}^{n-1},{\mathbb S}^{n-1})$. Let ${\mathbf a}=[1+\rho(\omega)]\omega+r_0
  \omega_0$ and ${\mathbf b}=[1+\rho(\omega)]\xi+[\nabla\rho(\omega)\cdot\xi]
  \omega$, where $\xi\in T_\omega({\mathbb S}^{n-1})$. Then a simple
  calculation shows that for any $\xi\in T_\omega({\mathbb S}^{n-1})$ we have
$$
  D_\omega f_z(\rho,\omega)\xi=\frac{|{\mathbf a}|^2 {\mathbf b}-
  ({\mathbf a}\cdot {\mathbf b}){\mathbf a}}{|{\mathbf a}|^3}.
$$
  Since ${\mathbf a}=\omega+O(\delta+\epsln)$, ${\mathbf b}=\xi+O(\delta)|\xi|$
  and $\omega\cdot\xi=0$, from the above expression we see that $D_\omega f_z
  (\rho,\omega)\xi=\xi+O(\delta+\epsln)|\xi|$, so that the desired assertion
  holds.

  It follows that for any $\rho\in {\mathcal O}_\delta({\mathbb S}^{n-1})$ and
  $z\in {\mathbb B}^n_\epsln$, the mapping $f_z(\rho,\cdot):{\mathbb S}^{n-1}
  \to {\mathbb S}^{n-1}$ is open. As a result, ${\rm Im}f_z(\rho,\cdot)$ is an
  open subset of ${\mathbb S}^{n-1}$. Since $f_z(\rho,\cdot)$ is continuous,
  ${\rm Im}\,f_z(\rho,\cdot)$ is also closed in ${\mathbb S}^{n-1}$. Thus,
  $f_z(\rho,\cdot):{\mathbb S}^{n-1}\to {\mathbb S}^{n-1}$ must be a surjection.

  Now let $g_z(\rho,\cdot)$ be the inverse of $f_z(\rho,\cdot)$. By the inverse
  function theorem we know that $g_z(\rho,\cdot)\in C^1({\mathbb S}^{n-1},
  {\mathbb S}^{n-1})$. Let $F_z(\rho,\omega)$ denote the right-hand side of
  (4.2). Substituting $\omega=g_z(\rho,\omega')$ into (4.2) we see that
$$
  \tilde{\rho}(\omega')=F_z(\rho,g_z(\rho,\omega')) \quad \mbox{for}\;\;
  \omega'\in {\mathbb S}^{n-1}.
\eqno{(4.8)}
$$
  Hence, the mapping $S_z^\ast$ is well-defined, and $S_z^\ast(\rho)=F_z(\rho,
  g_z(\rho,\cdot))$.

  Finally, it is clear that $F_z\in C^1({\mathcal O}_\delta({\mathbb S}^{n-1})
  \times {\mathbb S}^{n-1},{\mathbb R})$ and $f_z\in C^1({\mathcal O}_\delta
  ({\mathbb S}^{n-1})\times {\mathbb S}^{n-1},{\mathbb S}^{n-1})$. By the
  implicit function theorem, we also have $g_z\in C^1({\mathcal O}_\delta
  ({\mathbb S}^{n-1})\times {\mathbb S}^{n-1},{\mathbb S}^{n-1})$. Thus the
  mapping $(\rho,\omega)\to S_z^*(\rho)(\omega)$ from ${\mathcal O}_\delta
  ({\mathbb S}^{n-1})\times {\mathbb S}^{n-1}$ to ${\mathbb R}$ is of $C^1$
  class. Hence, we have $S_z^*\in C({\mathcal O}_\delta({\mathbb S}^{n-1}),
  C^1({\mathbb S}^{n-1}))\cap C^1({\mathcal O}_\delta({\mathbb S}^{n-1}),
  C({\mathbb S}^{n-1}))$. This completes the proof. $\quad\Box$

\medskip
  From the proof of Lemma 4.1 we can see that if $\rho\in {\mathcal
  O}_\delta^{m}({\mathbb S}^{n-1})=C^m({\mathbb S}^{n-1})\cap {\mathcal
  O}_\delta({\mathbb S}^{n-1})$ for some $m\geq 2$, then $F_z$ and $f_z$ are
  of $C^m$ class, which implies that $g_z$ and the mapping $(\rho,\omega)\to
  S_z^*(\rho)(\omega)=F_z(\rho,g_z(\rho,\omega))$ are of $C^m$ class, so that
  $S_z^*\in C^k({\mathcal O}_\delta^m({\mathbb S}^{n-1}), C^{m-k}({\mathbb
  S}^{n-1}))$ for any $0\leq k\leq m$. This particularly implies that
$$
  S_z^*\in C^\infty(C^\infty({\mathbb S}^{n-1})\cap {\mathcal
  O}_\delta({\mathbb S}^{n-1}), C^\infty({\mathbb S}^{n-1})).
\eqno{(4.9)}
$$
  Similarly, if $\rho\in {\mathcal O}_\delta^{m+\mu}({\mathbb S}^{n-1})=
  C^{m+\mu}({\mathbb S}^{n-1})\cap {\mathcal O}_\delta({\mathbb S}^{n-1})$ for
  some $m\geq 1$ and $0<\mu\leq 1$, then $S_z^*\in C^k({\mathcal
  O}_\delta^{m+\mu}({\mathbb S}^{n-1}),C^{m-k+\mu}({\mathbb S}^{n-1}))$ for
  any $0\leq k\leq m$. To establish a similar result for the space
  $B_{qq}^{m-{1/q}}({\mathbb S}^{n-1})$, we need the following lemma:

\medskip
  {\bf Lemma 4.2}\ \ {\em Let $\Omega_1$, $\Omega_2$ be two bounded smooth open
  subset of ${\mathbb R}^n$. Let $m\geq 2$ and $q>n/(m-1)$. Let $\Phi$ be a
  diffeomorphism from $\overline{\Omega}_1$ to $\overline{\Omega}_2$. Assume
  that $\Phi\in W^{m,q}(\Omega_1,{\mathbb R}^n)$. Then $\Phi^{-1}\in W^{m,q}
  (\Omega_2,{\mathbb R}^n)$. Moreover, given $\epsln>0$, the mapping $\Phi\to
  \Phi^{-1}$ from the set
$$
  \{\Phi\in W^{m,q}(\Omega_1,{\mathbb R}^n): |\det D\Phi(x)|\geq\epsln\;
  \mbox{for}\;\mbox{any}\;x\in\Omega_1\}
$$
  to $W^{m,q}(\Omega_2,{\mathbb R}^n)$
  is $C^\infty$, and there exists a continuous function $C_\epsln:[0,\infty)
  \to [0,\infty)$ such that
$$
  \|\Phi^{-1}\|_{W^{m,q}(\Omega_2,{\mathbb R}^n)}\leq
  C_\epsln(\|\Phi\|_{W^{m,q}(\Omega_1,{\mathbb R}^n)}).
\eqno{(4.10)}
$$
}

  {\em Proof}:\ \ We first note that the assumptions on $m$ and $q$ imply that
  $W^{m,q}(\Omega_1,{\mathbb R}^n)\hookrightarrow C^1(\overline{\Omega}_1,
  {\mathbb R}^n)$, and there exists constant $C>0$ such that $\|\Phi\|_{C^1
  (\overline{\Omega}_1,{\mathbb R}^n)}\leq C\|\Phi\|_{W^{m,q}(\Omega_1,{\mathbb
  R}^n)}$. In the sequel we denote by $x$ the variable in $\Omega_1$, and by
  $y$ the variable in $\Omega_2$. We also denote $\Psi=\Phi^{-1}$. Then we have
$$
  D\Psi(y)=[D\Phi(x)]^{-1}=[\det D\Phi(x)]^{-1}D^\ast\Phi(x),
$$
  where $D^\ast\Phi(x)$ denotes the co-matrix of the matrix $D\Phi(x)$. By this
  formula, the Leibnitz rule and the Gagliardo-Nirenberg inequality we can easily
  deduce that for any $\alpha\in {\mathbb Z}_+^n$ such that $0<|\alpha|\leq m$ and
  any $\epsln>0$ such that $|\det D\Phi(x)|\geq \epsln$
  for all $x\in\Omega_1$, we have
$$
  \|\partial^\alpha\Psi\|_{L^q(\Omega_2,{\mathbb R}^n)}
  \leq C\epsln^{-|\alpha|}\|D\Phi\|_{L^\infty(\Omega_1,{\mathbb R}^n)}^{n|\alpha|-2}
  \sum_{|\beta|=|\alpha|} \|\partial^{\beta}\Phi\|_{L^q(\Omega_1,{\mathbb R}^n)}
  \leq C\epsln^{-|\alpha|}\|\Phi\|_{W^{m,q}(\Omega_1,{\mathbb R}^n)}^{n|\alpha|-1}.
$$
  Hence (4.10) holds. The assertion that the mapping $\Phi\to\Phi^{-1}$ is
  smooth is an immediate consequence of the above argument. $\qquad\Box$

\medskip
  {\bf Lemma 4.3}\ \ {\em Let $m$ and $q$ be as in lemma 4.2. Then we have the
  following assertions:

  (i) For $\delta>0$ sufficiently small and for $z\in {\mathbb B}^n_\epsln$
  with $\epsln$ sufficiently small, we have
$$
  S_z^*\in C^k({\mathcal O}_\delta^{m,q}({\mathbb S}^{n-1}),
  B_{qq}^{m-k-{1/q}}({\mathbb S}^{n-1})),\quad 0\leq k\leq m-1.
$$
  In particular, $S_z^*\in C({\mathcal O}_\delta^{m,q}({\mathbb S}^{n-1}),
  B_{qq}^{m-{1/q}}({\mathbb S}^{n-1}))\cap C^1({\mathcal O}_\delta^{m,q}
  ({\mathbb S}^{n-1}),B_{qq}^{m-1-{1/q}}({\mathbb S}^{n-1}))$. Moreover, for
  any $\rho\in {\mathcal O}_\delta^{m,q}({\mathbb S}^{n-1})$ and $1\leq k\leq
  m-1$, the operator $DS_z^*(\rho)$ from $B_{qq}^{m-{1/q}}({\mathbb S}^{n-1})$
  to $B_{qq}^{m-1-{1/q}}({\mathbb S}^{n-1})$ can be uniquely extended to the
  space $B_{qq}^{m-k-{1/q}}({\mathbb S}^{n-1})$, such that
$$
  DS_z^*(\rho)\in L(B_{qq}^{m-k-{1/q}}({\mathbb S}^{n-1}),
  B_{qq}^{m-k-{1/q}}({\mathbb S}^{n-1})).
$$

  (ii) For any $z,w\in {\mathbb B}^n_\epsln$ with $\epsln$ sufficiently small,
  we have
$$
  S_z^*\circ S_w^*=S_{z+w}^*,\quad S_0^*=id, \quad \mbox{and}
  \quad (S_z^*)^{-1}=S_{-z}^*.
$$

  (iii) The mapping $S^*: z\to S_z^*$ from  ${\mathbb B}^n_\epsln$ to
  $C({\mathcal O}_\delta^{m,q}({\mathbb S}^{n-1})$, $B_{qq}^{m-{1/q}}
  ({\mathbb S}^{n-1}))$ is an injection, and
$$
  S^*\in C^k({\mathbb B}_\epsln^n,C^l({\mathcal O}_\delta^{m,q}({\mathbb S}^{n-1}),
  B_{qq}^{m-k-l-{1/q}}({\mathbb S}^{n-1}))),\quad k\geq 0, \;\; l\geq 0,\;\;
  k+l\leq m-1.
$$

  (iv) Finally assume that $2\leq k<m$, $q>n/(k-1)$ and define $p: {\mathbb
  B}^n_\epsln\times {\mathcal O}_\delta^{k,q}({\mathbb S}^{n-1})\to
  B_{qq}^{k-{1/q}}({\mathbb S}^{n-1})$ by $p(z,\rho)=S_z^*(\rho)$. Then for any
  $\rho\in {\mathcal O}_\delta^{m,q}
  ({\mathbb S}^{n-1})$ the mapping $z\to p(z,\rho)$ from ${\mathbb B}^n_\epsln$
  to $B_{qq}^{m-{1/q}}({\mathbb S}^{n-1})$ is Fr\'{e}chet differentiable when
  regarded as a mapping from ${\mathbb B}^n_\epsln$ to $B_{qq}^{k-{1/q}}
  ({\mathbb S}^{n-1})$, and we have ${\rm rank}\,D_z p(z,\rho)=n$ for any $z\in
  {\mathbb B}^n_\epsln$ and $\rho\in {\mathcal O}_\delta^{m,q}({\mathbb
  S}^{n-1})$. If furthermore $\rho\in C^\infty({\mathbb S}^{n-1})$ then
  $[z\to p(z,\rho)]\in C^\infty({\mathbb B}^n_\epsln,C^\infty({\mathbb
  S}^{n-1}))\subseteq C^\infty({\mathbb B}^n_\epsln,B_{qq}^{m-{1/q}}({\mathbb
  S}^{n-1}))$.}
\medskip

  {\em Proof}:\ \ We first note that the assumptions on $m$ and $q$ imply that
  $B_{qq}^{m-{1/q}}({\mathbb S}^{n-1})\hookrightarrow C^1({\mathbb S}^{n-1})$,
  so that by Lemma 4.1, $S_z^\ast(\rho)$ makes sense for $\rho\in
  {\mathcal O}_\delta^{m,q}({\mathbb S}^{n-1})$. Next, by (4.8) we see that
  $S_z^\ast(\rho)=F_z(\rho,g_z(\rho,\cdot))$. Considering (4.2) and (4.3), for
  given $z=r_0\omega_0\in {\mathbb B}^n_\epsln$ and any $u\in W^{m,q}({\mathbb
  B}^n)$ such that $\|u\|_{C^1({\mathbb B}^n)}<\delta$ we define
$$
  \widetilde{u}(x')=\sqrt{[1+u(x)]^2+r_0^2+2r_0[1+u(x)]x\cdot\omega_0}-1,
  \quad x\in \overline{{\mathbb B}}^n,
\eqno{(4.11)}
$$
  where $x'$ and $x$ are related by
$$
  x'={[1+u(x)]x+r_0\omega_0\over |[1+u(x)]x+r_0\omega_0|}|x|\phi(|x|-1)+
  [1-\phi(|x|-1)]x, \quad x\in \overline{{\mathbb B}}^n,
\eqno{(4.12)}
$$
  where $\phi$ is as in Section 3. As before we use the notation $F_z(u,x)$ to
  denote the expression on the right-hand side of (4.11). Since the assumptions
  on $m$ and $q$ imply that $W^{m,q}({\mathbb B}^n)$ is an algebra, it is clear
  that $F_z(u,\cdot)\in W^{m,q}({\mathbb B}^n)$, and the mapping $u\to F_z(u,
  \cdot)$ is $C^\infty$. We also use the same notation $f_z(u,x)$ as before to
  denote the expression on the right-hand side of (4.12), because if we
  particularly take $u=\Pi(\rho)$ and $x=\omega\in {\mathbb S}^{n-1}$ then we
  get $f_z(\rho,\omega)$ defined before. It can be easily shown that if
  $\epsln$ and $\delta$ are sufficiently small then the mapping $\Phi_u: x\to
  x'=f_z(u,x)$ is a diffeomorphism of $\overline{{\mathbb B}^n}$ to itself and
  $\det D\Phi_u(x)=1+O(\epsln+\delta)$. Moreover, since $W^{m,q}({\mathbb B}^n)$
  is an algebra, we have $\Phi_u\in W^{m,q}({\mathbb B}^n,{\mathbb R}^n)$ and
  it is clear that the mapping $u\to\Phi_u$ is $C^\infty$. By Lemma 4.2 we infer
  that $\Phi_u^{-1}\in W^{m,q}({\mathbb B}^n,{\mathbb R}^n)$, and the mapping
  $\Phi_u\to\Phi_u^{-1}$ is $C^\infty$. Substituting $x=\Phi_u^{-1}(x')$ into
  the right-hand side of (4.11) and using Lemma 3.1, we see that $\widetilde{u}
  =F_z(u,\Phi_u^{-1}(\cdot))\in W^{m,q}({\mathbb B}^n)$. Now, clearly if $u=
  \Pi(\rho)$ for some $\rho\in {\mathcal O}_\delta^{m,q}({\mathbb S}^{n-1})$
  then we have $\widetilde{u}|_{{\mathbb S}^{n-1}}=S_z^*(\rho)$, so that we
  have proved that $S_z^*(\rho)\in B_{qq}^{m-{1/q}}({\mathbb S}^{n-1})$ for any
   $\rho\in {\mathcal O}_\delta^{m,q}({\mathbb S}^{n-1})$. We note that though
  both the mappings $u\to F_z(u,\cdot)$ and $u\to\Phi_u^{-1}$ are $C^\infty$,
  the mapping $u\to\widetilde{u}=F_z(u,\Phi_u^{-1}(\cdot))$ is, however, not
  necessarily $C^\infty$, because $F_z(u,x)$ is generally not $C^\infty$ in
  $x$. Despite of this inconvenience, we still can ensure that the mapping
  $u\to\widetilde{u}=F_z(u,\Phi_u^{-1}(\cdot))$ from $W^{m,q}({\mathbb B}^n)
  \cap \{u\in C^1(\overline{{\mathbb B}^n}):\|u\|_{C^1(\overline{{\mathbb
  B}}^n)}<\delta\}$ to $W^{m,q}({\mathbb B}^n)$ is continuous, because both
  $(u,x)\to F_z(u,x)$ and $u\to\Phi_u^{-1}$ are continuous. Thus $S_z^*\in
  C({\mathcal O}_\delta^{m,q}({\mathbb S}^{n-1}),B_{qq}^{m-{1/q}}({\mathbb
  S}^{n-1}))$. Next, since $S_z^*(\rho)=\Gamma F_z(\Pi(\rho),
  \Phi_{\Pi(\rho)}^{-1}(\cdot))$, where $\Gamma$ denotes the trace operator,
  we have
$$
\begin{array}{rl}
  DS_z^*(\rho)\eta=&\Gamma D_1 F_z(\Pi(\rho),\Phi_{\Pi(\rho)}^{-1}(\cdot))
  \Pi(\eta)+\Gamma D_2 F_z(\Pi(\rho),\Phi_{\Pi(\rho)}^{-1}(\cdot))D_u
  \Phi_{\Pi(\rho)}^{-1}(\cdot)\Pi(\eta)\\
  \equiv & I(\rho)\eta+I\!I(\rho)\eta,
\end{array}
\eqno{(4.13)}
$$
  where $D_1$ and $D_2$ represent the Fr\'{e}chet derivatives in the first
  and the second arguments, respectively, and $D_u\Phi_{\Pi(\rho)}^{-1}=
  D_u\Phi_u^{-1}|_{u=\Pi(\rho)}$. By Lemma 3.2 $(i)$ it is obvious that
$$
\begin{array}{rl}
  I\in & C({\mathcal O}_\delta^{m,q}({\mathbb S}^{n-1}),
  L(B_{qq}^{m-{1/q}}({\mathbb S}^{n-1}), B_{qq}^{m-{1/q}}({\mathbb
  S}^{n-1})))\\
  &\bigcap C({\mathcal O}_\delta^{m,q}({\mathbb S}^{n-1}),L(B_{qq}^{m-k-{1/q}}
  ({\mathbb S}^{n-1}), B_{qq}^{m-k-{1/q}}({\mathbb S}^{n-1}))), \quad
  1\leq k\leq m-1.
\end{array}
$$
  To treat $I\!I$
  we denote $G_z(t,y)=\sqrt{(1+t)^2+r_0^2+2r_0(1+t)y\cdot\omega_0}-1$ for
  $t\in {\mathbb R}$ and $y\in {\mathbb B}^n$. Then $F_z(u,x)=G_z(u(x),x)$
  for $u\in W^{m,q}({\mathbb B}^n)$ and $x\in {\mathbb B}^n$, so that
$$
  D_2 F_z(u,x)=D_1 G_z(u(x),x)Du(x)+D_2 G_z(u(x),x).
$$
  Given $u\in W^{m,q}({\mathbb B}^n)$, from the above expression of $D_2 F_z
  (u,x)$ we see that $D_2 F_z(u,\cdot)=[x\to D_2 F_z(u,x)]\in W^{m-1,q}
  ({\mathbb B}^n,L({\mathbb R}^n,{\mathbb R}))$. Besides, since
$$
  D_u\Phi_u^{-1}\in L(W^{m,q}({\mathbb B}^n),W^{m,q}({\mathbb B}^n,
  {\mathbb R}^n))\cap L(W^{m-k,q}({\mathbb B}^n),W^{m-k,q}({\mathbb B}^n,
  {\mathbb R}^n))
$$
  ($1\leq k\leq m-1$), for any $\eta\in B^{m-1/q}_{qq}({\mathbb S}^{n-1})$ we
  have $D_u\Phi_u^{-1}(\cdot)\pi(\eta)=[x\to D_u\Phi_u^{-1}(x)\pi(\eta)]\in
  W^{m,q}({\mathbb B}^n,{\mathbb R}^n)$, and if $\eta\in B^{m-k-1/q}_{qq}
  ({\mathbb S}^{n-1})$ for some $1\leq k\leq m-1$ then we have $D_u\Phi_u^{-1}
  (\cdot)\pi(\eta)\in W^{m-k,q}({\mathbb B}^n,{\mathbb R}^n)$. Hence, given
  $\rho\in B^{m-1/q}_{qq}({\mathbb S}^{n-1})$, for any $\eta\in B^{m-1/q}_{qq}
  ({\mathbb S}^{n-1})$ we have
$$
   D_2 F_z(\pi(\rho),\Phi_{\pi(\rho)}^{-1}(\cdot))
   D_u\Phi_{\pi(\rho)}^{-1}(\cdot)\pi(\eta)
  =[x\to D_2 F_z(\pi(\rho),\Phi_{\pi(\rho)}^{-1}(x))
   D_u\Phi_{\pi(\rho)}^{-1}(x)\pi(\eta)]\in W^{m-1,q}({\mathbb B}^n),
$$
  and if $\eta\in B^{m-k-1/q}_{qq}({\mathbb S}^{n-1})$ for some $1\leq k\leq
  m-1$ then we have
$$
   D_2 F_z(\pi(\rho),\Phi_{\pi(\rho)}^{-1}(\cdot))
   D_u\Phi_{\pi(\rho)}^{-1}(\cdot)\pi(\eta)\in W^{m-k,q}({\mathbb B}^n).
$$
  This implies that for $\rho\in B^{m-1/q}_{qq}({\mathbb S}^{n-1})$, if $\eta
  \in B^{m-1/q}_{qq}({\mathbb S}^{n-1})$ then $I\!I(\rho)\eta\in
  B^{m-1-1/q}_{qq}({\mathbb S}^{n-1})$, whereas if $\eta\in B^{m-k-1/q}_{qq}
  ({\mathbb S}^{n-1})$ for some $1\leq k\leq m-1$ then $I\!I(\rho)\eta\in
  B^{m-k-1/q}_{qq}({\mathbb S}^{n-1})$. A similar analysis shows that
$$
\begin{array}{rl}
  I\!I\in & C({\mathcal O}_\delta^{m,q}({\mathbb S}^{n-1}),
  L(B_{qq}^{m-{1/q}}({\mathbb S}^{n-1}),B_{qq}^{m-1-{1/q}}({\mathbb
  S}^{n-1})))\\
  &\bigcap C({\mathcal O}_\delta^{m,q}({\mathbb S}^{n-1}),
  L(B_{qq}^{m-k-{1/q}}({\mathbb S}^{n-1}),B_{qq}^{m-k-{1/q}}({\mathbb
  S}^{n-1}))),\quad 1\leq k\leq m-1.
\end{array}
$$
  Hence, $S_z^*\in C^1({\mathcal O}_\delta^{m,q}({\mathbb S}^{n-1}),
  B_{qq}^{m-1-{1/q}}({\mathbb S}^{n-1}))$, and
$$
\begin{array}{rl}
  DS_z^*\in & C({\mathcal O}_\delta^{m,q}({\mathbb S}^{n-1}),
  L(B_{qq}^{m-{1/q}}({\mathbb S}^{n-1}),B_{qq}^{m-1-{1/q}}({\mathbb
  S}^{n-1})))\\
  &\bigcap C({\mathcal O}_\delta^{m,q}({\mathbb S}^{n-1}),
  L(B_{qq}^{m-k-{1/q}}({\mathbb S}^{n-1}),B_{qq}^{m-k-{1/q}}({\mathbb
  S}^{n-1}))),\quad 1\leq k\leq m-1.
\end{array}
$$
  Furthermore, by an induction argument we see that $S_z^*\in C^k({\mathcal
  O}_\delta^{m,q}({\mathbb S}^{n-1}),B_{qq}^{m-k-{1/q}}({\mathbb S}^{n-1}))$
  for any $0\leq k\leq m-1$. This proves Assertion $(i)$. Assertion $(ii)$
  is obvious. The first part of Assertion $(iii)$ is evident, and the second
  part follows by checking more carefully the argument in the proof of
  Assertion $(i)$, whcih we omit here. From the proof of Assertion
  $(i)$ we see that for any integers $2\leq k<m$ and $q>n/(k-1)$, the mapping
  $p:{\mathbb B}^n_\epsln\times {\mathcal O}_\delta^{k,q}({\mathbb S}^{n-1})\to
  B_{qq}^{k-{1/q}}({\mathbb S}^{n-1})$ defined by $p(z,\rho)=S_z^*(\rho)$ is
  continuously differentiable at any point $(z,\rho)\in {\mathbb B}^n_\epsln
  \times {\mathcal O}_\delta^{m,q}({\mathbb S}^{n-1})$. Moreover, a simple
  calculation shows that $D_1 p(0,0)z=z\cdot\omega$. Here $z\cdot\omega$
  represents the function $\omega\to z\cdot\omega$ on ${\mathbb S}^{n-1}$.
  This shows that ${\rm rank}D_1 p(0,0)=n$. By continuity, we infer that
  ${\rm rank}D_1 p(z,\rho)=n$ for any $(z,\rho)\in {\mathbb B}^n_\epsln\times
  {\mathcal O}_\delta^{m,q}({\mathbb S}^{n-1})$, provided $\epsln$ and $\delta$
  are sufficiently small. Finally, if  $\rho\in C^\infty({\mathbb S}^{n-1})$
  then from the construction of $S_z^*$ it is clear that $[z\to p(z,\rho)]\in
  C^\infty({\mathbb B}^n_\epsln,C^\infty({\mathbb S}^{n-1}))$. Hence, Assertion
  $(iv)$ follows. This completes the proof. $\quad\Box$
\medskip

  By Lemma 4.3 we see that the mapping $S^\ast$ provides an action of the local
  group $G={\mathbb B}^n_\epsln$ to some open subset of $B_{qq}^{m-{1/q}}
  ({\mathbb S}^{n-1})$. We note that if $c=0$ then by some similar arguments as
  in Section 3 we can reduce the problem (1.1)--(1.5) and (1.7) into a
  differential equation $\rho'(t)={\mathcal A}_\gamma (\rho(t))$ in the
  Banach space $B_{qq}^{m-3-{1/q}}({\mathbb S}^{n-1})$ in the unknown function
  $\rho=\rho(t)$ only, where ${\mathcal A}_\gamma$ is defined in some open
  subset of $B_{qq}^{m-{1/q}}({\mathbb S}^{n-1})$. It can be shown that in this
  case the reduced equation satisfies a similar relation as that in (4.1) under
  the above action of $G$ (cf. the proof of Lemma 4.6 below). For the equation
  (3.30), however, $G$ has to act on some open set in $X$. This is fulfilled in
  the following paragraph. In the sequel, the notations $X$, $X_0$ and
  ${\mathcal O}$ have the same meaning as introduced in the end of Section 3.

\medskip
  Given $z\in {\mathbb B}^n_\epsln$ and $\rho\in {\mathcal O}_\delta({\mathbb
  S}^{n-1})$, let $P_{z,\rho}:C({\mathbb B}^n)\to C({\mathbb B}^n)$ be the
  mapping
$$
   P_{z,\rho}(u)(x)=u(\Theta_{\rho}^{-1}(\Theta_{S_z^\ast(\rho)}(x)-z))\quad
   \mbox{for}\;\;u\in C({\mathbb B}^n).
$$
  Clearly, $P_{z,\rho}\in L(C({\mathbb B}^n),C({\mathbb B}^n))$. Moreover, if
  $\rho\in B_{qq}^{m-{1/q}}({\mathbb S}^{n-1})$ then $S_z^\ast(\rho)\in
  B_{qq}^{m-{1/q}}({\mathbb S}^{n-1})$, so that by Lemma 3.1 we have $P_{z,
  \rho}\in L(W^{m,q}({\mathbb B}^n)$, $W^{m,q}({\mathbb B}^n))$.
  For $u\in C({\mathbb B}^n)$, $\rho\in {\mathcal O}_\delta({\mathbb
  S}^{n-1})$ and $z\in {\mathbb B}^n_\epsln$ we denote
$$
  {\mathbf S}_z^\ast \left(
\begin{array}{c}
  u\\
  \rho
\end{array}
\right)=\left(
\begin{array}{c}
  P_{z,\rho}(u)\\
  S_z^\ast(\rho)
\end{array}
\right).
$$
  Note that ${\mathbf S}_0^\ast=id$.
\medskip

  {\bf Lemma 4.4}\ \ {\em Let $m\geq 5$ and $q>n/(m-4)$. Let
$$
  {\mathcal O}'=W^{m-3}({\mathbb B}^n)\times (B^{m-3-1/q}_{qq}({\mathbb
  S}^{n-1})\cap {\mathcal O}_\delta({\mathbb S}^{n-1})) \quad
  (\Longrightarrow {\mathcal O}=X_0\cap {\mathcal O}').
$$
  For sufficiently small $\epsln>0$ and $\delta>0$ we have the following
  assertions:

  $(i)$ For any $\epsln\in {\mathbb B}^n_\epsln$ we have ${\mathbf S}_z^\ast
  \in C({\mathcal O}',X)\cap C({\mathcal O},X_0)$. Moreover, regarded as a
  mapping from ${\mathcal O}'$ to $X$, ${\mathbf S}_z^\ast$ is Fr\'{e}chet
  differentiable at every point in ${\mathcal O}$, and $D{\mathbf S}_z^\ast
  \in C({\mathcal O},L(X,X))$.

  $(ii)$ For any $z,w\in {\mathbb B}^n_\epsln$ we have
$$
  {\mathbf S}_z^*\circ {\mathbf S}_w^*={\mathbf S}_{z+w}^*,\quad
  {\mathbf S}_0^*=id, \quad \mbox{and}
  \quad ({\mathbf S}_z^*)^{-1}={\mathbf S}_{-z}^*.
$$

  $(iii)$ The mapping ${\mathbf S}^*: z\to {\mathbf S}_z^*$ from ${\mathbb
  B}^n_\epsln$ to $C({\mathcal O}',X)$ is an injection, and
$$
  {\mathbf S}^*\in C^k({\mathbb B}^n_\epsln,C^l({\mathcal O},
  W^{m-k-l-1,q}({\mathbb B}^n)\times B_{qq}^{m-k-l-{1/q}}
  ({\mathbb S}^{n-1}))),\quad k\geq 0, \;\; l\geq 0,\;\; k+l\leq m-1.
$$

  (iv) Define $p:{\mathbb B}^n_\epsln\times {\mathcal O}'\to X$ by
  $p(z,U)={\mathbf S}_z^*(U)$. Then for any $U\in {\mathcal O}$ we have
  $p(\cdot,U)\in C^1({\mathbb B}^n_\epsln,X)$, and ${\rm rank}D_z p(z,U)=n$
  for every $z\in {\mathbb B}^n_\epsln$ and $U\in {\mathcal O}$. If
  furthermore $U\in X^\infty=C^\infty(\overline{{\mathbb B}^n})\times
  C^\infty({\mathbb S}^{n-1})$ then $p(\cdot,U)\in C^\infty({\mathbb
  B}^n_\epsln,X^\infty)$.}
\medskip

  {\em Proof}:\ \ All assertions of this lemma follow readily from the
  corresponding assertions in Lemma 4.3. $\quad\Box$
\medskip

  In the sequel, for $\rho=\rho(t)$, $u=u(x,t)$ and $U=\Big(^{u(x,t)}_{\;\;
  \rho(t)}\Big)$,  we denote by $P_{z,\rho}(u)$ the function $\widetilde{u}(x,
  t)=u(\Theta_{\rho(t)}^{-1}(\Theta_{S_z^\ast(\rho(t))}(x)-z),t)$, by
  $S_z^\ast(\rho)$ the function $\widetilde{\rho}(t)=S_z^\ast(\rho(t))$, and
  by ${\mathbf S}_z^\ast (U)$ the vector function $\Big(^{P_{z,\rho}(u)}_{\;\,
  S_z^\ast(\rho)}\Big)=\Big(^{\widetilde{u}(x,t)}_{\;\,\widetilde{\rho}(t)}
  \Big)$.
\medskip

  {\bf Lemma 4.5}\ \ {\em If $U=\Big(^u_\rho\Big)$ is a solution of the
  equation (3.30) such that $\|\rho\|_{C^1({\mathbb S}^{n-1})}$ is sufficiently
  small, then for any $z\in {\mathbb R}^n$ such that $|z|$ is sufficiently
  small, ${\mathbf S}_z^\ast (U)=\Big(^{P_{z,\rho}(u)}_{\;\,S_z^\ast(\rho)}
  \Big)$ is also a solution of (3.30).}
\medskip

  {\em Proof}:\ \ It is easy to see that if a triple $(\sigma,p,\Omega)$ is a
  solution of the problem (1.1)--(1.5), then for any $z\in {\mathbb R}^n$ the
  triple $(\widetilde{\sigma},\widetilde{p},\widetilde{\Omega})$ defined by
$$
  \widetilde{\sigma}(x,t)=\sigma(x-z,t), \quad \widetilde{p}(x,t)=p(x-z,t),
  \quad \widetilde{\Omega}(t)=\Omega(t)+z,
$$
  is also a solution of that problem. From this fact one can easily verify
  that if $U=\Big(^u_\rho\Big)$ is a solution of the equation (3.30) then
  $\widetilde{U}=\Big(^{\widetilde{u}}_{\widetilde{\rho}}\Big)$, where
$$
  \widetilde{u}(x,t)=u(\Theta_{\rho(t)}^{-1}
  (\Theta_{S_z^\ast(\rho(t))}(x)-z),t),
  \quad \widetilde{\rho}(t)=S_z^\ast(\rho(t)),
$$
  is also a solution of that equation, which is the desired assertion.
  $\qquad\Box$
\medskip

  {\bf Lemma 4.6}\ \ {\em The following relation holds for any $z\in {\mathbb
  B}^n_\epsln$ and any $U=\Big(^u_\rho\Big)\in {\mathcal O}$, provided $\epsln$
  and $\delta$ are sufficiently small:}
$$
  {\mathbb F}({\mathbf S}_z^\ast (U))=D{\mathbf S}_{z}^*(U){\mathbb F}(U).
\eqno{(4.14)}
$$

  {\em Proof}:\ \ By Theorem 1.1 of \cite{Cui6}, given any $U=\Big(^u_\rho\Big)
  \in X_0$ there exists $\delta>0$ such that the equation (3.30) has a unique
  solution $V=V(t)$ for $0\leq t\leq\delta$, which belongs to $C([0,\delta],X)
  \cap C((0,\delta],{\mathcal O})\cap L^\infty((0,\delta),X_0)\cap C^1((0,
  \delta],X)$ and satisfies the initial condition $V(0)=U$ (This result also
  follows from Corollary 5.3 in the next section and a standard existence
  theorem that we used in the proof of Theorem 2.1). Let $\widetilde{V}(t)=
  {\mathbf S}_z^\ast (V(t))$ for $0\leq t\leq\delta$. By Lemma 4.5,
  $\widetilde{V}$ is also a solution of (3.30), satisfying the initial
  condition $\widetilde{V}(0)={\mathbf S}_z^\ast(U)$. The fact that
  $\widetilde{V}$ is the solution of (3.30) implies that
$$
  {d\widetilde{V}(t)\over dt}={\mathbb F}(\widetilde{V}(t)) \quad
  \mbox{for}\;\; 0<t\leq\delta.
$$
  On the other hand, since $\widetilde{V}(t)={\mathbf S}_z^\ast (V(t))$, we
  have
$$
  {d\widetilde{V}(t)\over dt}=D{\mathbf S}_z^\ast (V(t)){dV(t)\over dt}
  =D{\mathbf S}_z^\ast (V(t)){\mathbb F}(V(t))
   \quad \mbox{for}\;\; 0<t\leq\delta.
$$
  Thus ${\mathbb F}(\widetilde{V}(t))=D{\mathbf S}_z^\ast (V(t)){\mathbb F}
  (V(t))$ for $0<t\leq\delta$. If $V(t)$ is a strict solution then clearly
  $\widetilde{V}(t)$ is also a strict solution, so that by directly letting
  $t\to\delta^+$ we get (4.14). If $V(t)$ is not a strict solution then we
  appeal to the quasi-linear structure of ${\mathbb F}(U)$ to prove (4.14):
  Since $V\in L^\infty((0,\delta),X_0)\cap C([0,\delta],X)$ and $V(0)=U$, we
  infer that $V(t)$ weakly converges to $U$ in $X_0$ as $t\to 0^+$. Similarly
  $\widetilde{V}(t)$ weakly converges to ${\mathbf S}_z^\ast(U)$ in $X_0$.
  Since ${\mathbb F}(U)={\mathbb A}(U)U+{\mathbb F}_0(U)$, we have
$$
\begin{array}{rl}
  {\mathbb F}(V(t))-{\mathbb F}(U)&=[{\mathbb A}(V(t))-{\mathbb A}(U)]V(t)
  +{\mathbb A}(U)[V(t)-U]+[{\mathbb F}_0(V(t))-{\mathbb F}(U)]\\
  &\equiv I(t)+I\!I(t)+I\!I\!I(t).
\end{array}
$$
  We have $\|I(t)\|_X\leq C\|{\mathbb A}(V(t))-{\mathbb A}(U)\|_{L(X_0,X)}$,
  so that $\lim_{t\to 0^+}\|I(t)\|_X=0$, because ${\mathbb A}$ maps $X_0$
  compactly into $L(X_0,X)$. We also have $\lim_{t\to 0^+}\|I\!I\!I(t)\|_X=0$
  by a similar reason. In addition, it is evident that $I\!I(t)$ weakly
  converges to $0$ in $X$ as $t\to 0^+$. Therefore, ${\mathbb F}(V(t))$
  weakly converges to ${\mathbb F}(U)$ in $X$. Similarly, ${\mathbb F}
  (\widetilde{V}(t))$ weakly converges to ${\mathbb F}({\mathbf S}_z^\ast(U))$
  in $X$. Finally, from the expression of $D{\mathbf S}_z^\ast$ (cf. (4.13))
  we can easily find that $D{\mathbf S}_z^\ast$ maps $X_0$ compactly into
  $L(X,X)$. Thus by a similar argument as above we infer that $D{\mathbf
  S}_z^\ast(V(t)){\mathbb F}(V(t))$ weakly converges to $D{\mathbf S}_z^\ast(U)
  {\mathbb F}(U)$ in $X$ as $t\to 0^+$. Hence (4.14) holds. $\quad\Box$
\medskip

  Lemma 4.6 has some obvious corollaries. First, let ${\mathbb F}_2$ be the
  second component of ${\mathbb F}$. Taking the second components of both
  sides of (4.14) we get
$$
  {\mathbb F}_2({\mathbf S}_z^\ast (U))=DS_{z}^*(\rho){\mathbb F}_2(U),
\eqno{(4.15)}
$$
  where $\rho$ is the second component of $U$. Next, let $u_s=\sigma_s-
  \bar{\sigma}$ and $U_s=\Big(^{u_s}_{\;\,0}\Big)$. $U_s$ is a stationary point
  of the equation (3.30), i.e., ${\mathbb F}(U_s)=0$. Taking $U=U_s$ in (4.14)
  we get ${\mathbb F}({\mathbf S}_z^\ast (U_s))=0$ for any $z\in {\mathbb
  B}^n_\epsln$. Since clearly $U_s\in X^\infty$, we have $[z\to {\mathbf
  S}_z^\ast (U_s)]\in C^\infty({\mathbb B}^n_\epsln,X^\infty)$. Thus,
  differentiating the equation ${\mathbb F}({\mathbf S}_z^\ast (U_s))=0$ in $z$
  at $z=0$, we obtain
$$
  {\mathbb F}'(U_s)W_j=0, \quad
  W_j=\left(\begin{array}{c} [\phi(r-1)r-1]\sigma_s'(r)\omega_j\\
  \omega_j\end{array}\right), \quad j=1,2,\cdots,n,
\eqno{(4.16)}
$$
  i.e., $0$ is an eigenvalue of ${\mathbb F}'(U_s)$ of (geometric) multiplicity
  $n$, and the corresponding linearly independent eigenvectors are $W_1$, $W_2$,
  $\cdots$, $W_n$.
\medskip

\section{Calculation of ${\mathbb F}'(U_s)$}

  In this section we calculate the Fr\'{e}chet derivative of ${\mathbb F}$ at
  the stationary point $U_s$. Since ${\mathbb F}(U)={\mathbb A}(U)U+{\mathbb
  F}_0(U)$, we have
$$
  {\mathbb F}'(U_s)V={\mathbb A}(U_s)V+[{\mathbb A}'(U_s)V]U_s
  +{\mathbb F}_0'(U_s)V \quad \mbox{for} \quad V\in X_0.
\eqno{(5.1)}
$$
  Recall that ${\mathbb A}\in C^\infty({\mathcal O},L(X_0,X))$, so that
  ${\mathbb A}'(U_s)\in L(X_0,L(X_0,X))$, and ${\mathbb A}'(U_s)V\in L(X_0,X)$
  for $V\in X_0$. Since $U_s=\Big(^{u_s}_{\;0}\Big)$, simple calculations
  show that for any $V=\Big(^v_\eta\Big)\in X_0$ we have
$$
  {\mathbb A}(U_s)V=\left(
\begin{array}{c}
  c^{-1}{\mathcal A}(0)v+{\mathcal C}(0,u_s)\eta\\
  {\mathcal B}(0)\eta
\end{array}
  \right),\qquad
  [{\mathbb A}'(U_s)V]U_s=\left(
\begin{array}{c}
  c^{-1}[{\mathcal A}'(0)\eta]u_s\\
  0
\end{array}
  \right),
\eqno{(5.2)}
$$
  and
$$
    {\mathbb F}_0'(U_s)V=\left(
\begin{array}{c}
  D_u{\mathcal F}_2(0,u_s)v+D_\rho{\mathcal F}_2(0,u_s)\eta\\
  D_u{\mathcal G}_2(0,u_s)v+D_\rho{\mathcal G}_2(0,u_s)\eta
\end{array}
  \right),
\eqno{(5.3)}
$$
  where $D_u{\mathcal F}_2$ and $D_\rho{\mathcal F}_2$ represent Fr\'{e}chet
  derivatives of ${\mathcal F}_2(\rho,u)$ in $u$ and $\rho$, respectively, and
  similarly for $D_u{\mathcal G}_2$ and $D_\rho{\mathcal G}_2$. Clearly,
$$
  {\mathcal A}(0)v=\Delta v,
\eqno{(5.4)}
$$
  and a simple computation shows that
$$
  {\mathcal C}(0,u_s)\eta=\phi(r-1)\sigma_s'(r)\Pi({\mathcal B}(0)\eta).
\eqno{(5.5)}
$$
  To compute ${\mathcal B}(0)\eta=-\gamma {\mathcal D}(0){\mathcal T}(0)
  {\mathcal L}(0)\eta$ we first note that, clearly,
$$
  {\mathcal D}(0)v={\partial v\over\partial r}\Big|_{r=1} \quad \mbox{and}
  \quad {\mathcal T}(0)\eta=\Pi(\eta).
$$
  Next, recall that
$$
  {\mathcal K}(\rho)={\mathcal L}(\rho)\rho+{\mathcal K}_1(\rho), \quad
  \mbox{so that} \quad
  {\mathcal K}'(0)\eta={\mathcal L}(0)\eta+{\mathcal K}_1'(0)\eta.
\eqno{(5.6)}
$$
  On the other hand, from \cite{FriRei3} we know that
$$
  {\mathcal K}(\epsln\eta)=1-\epsln[\eta(\omega)+
  {1\over n\!-\!1}\Delta_\omega \eta(\omega)]+o(\epsln),
$$
  which implies that ${\mathcal K}(0)=1$ and ${\mathcal K}'(0)\eta=-[\eta+
  {1\over n\!-\!1}\Delta_\omega\eta]$. Comparing these expressions with those
  in (5.6), we obtain
$$
  {\mathcal L}(0)\eta=-{1\over n\!-\!1}\Delta_\omega \eta, \quad
  {\mathcal K}_1(0)=1, \quad \mbox{and} \quad {\mathcal K}_1'(0)\eta=-\eta.
$$
  Hence we have
$$
  {\mathcal B}(0)\eta=-\gamma {\mathcal D}(0){\mathcal T}(0){\mathcal L}(0)\eta
  ={\gamma\over n\!-\!1}{\partial\over\partial r}
  \Pi(\Delta_\omega\eta)\Big|_{r=1}.
\eqno{(5.7)}
$$
  We denote $u^s_{\epsln,\eta}=\Theta^\ast_{\epsln\eta}\sigma_s-\bar{\sigma}$.
  Then we have
$$
  {\mathcal A}(\epsln\eta)u^s_{\epsln,\eta}=
  {\mathcal F}(u^s_{\epsln,\eta}+\bar{\sigma}),
$$
  so that
$$
  [{\mathcal A}(\epsln\eta)-{\mathcal A}(0)]u^s_{\epsln,\eta}
  +{\mathcal A}(0)(u^s_{\epsln,\eta}-u_s)
  ={\mathcal F}(u^s_{\epsln,\eta}+\bar{\sigma})
  -{\mathcal F}(u_s+\bar{\sigma}).
$$
  Dividing both sides with $\epsln$ and letting $\epsln\to 0$, we get
$$
  [{\mathcal A}'(0)\eta]u_s+{\mathcal A}(0)[{\mathcal M}(0,\sigma_s)
  \Pi(\eta)]={\mathcal F}'(\sigma_s)[{\mathcal M}(0,\sigma_s)\Pi(\eta)].
$$
  Here we used the fact that $\lim_{\epsln\to 0}\epsln^{-1}(u^s_{\epsln,\eta}-
  u_s)={\mathcal M}(0,\sigma_s)\Pi(\eta)$ ($=\phi(r-1)\sigma_s'(r)\Pi(\eta)$).
  Hence,
\setcounter{equation}{7}
\begin{eqnarray}
  [{\mathcal A}'(0)\eta]u_s & =&-{\mathcal A}(0)[{\mathcal M}(0,\sigma_s)
  \Pi(\eta)]+{\mathcal F}'(\sigma_s)[{\mathcal M}(0,\sigma_s)\Pi(\eta)]
\nonumber\\
  &=&-\Delta[\phi(r-1)\sigma_s'(r)\Pi(\eta)]+
  f'(\sigma_s(r))\phi(r-1)\sigma_s'(r)\Pi(\eta)
\nonumber\\
  &=&-[\Delta-f'(\sigma_s(r))]\big[\phi(r-1)\sigma_s'(r)\Pi(\eta)\big].
\end{eqnarray}
  To compute ${\mathbb F}_0'(U_s)$, we first note that since ${\mathcal P}(\rho,
  u)v={\mathcal M}(\rho,u)\Pi({\mathcal D}(\rho)v)$, we have
$$
\begin{array}{rcl}
  {\mathcal F}_2(\rho,u)&=&-c^{-1}{\mathcal F}(u+\bar{\sigma})-\gamma
  {\mathcal P}(\rho,u+\bar{\sigma}){\mathcal T}(\rho){\mathcal K}_1(\rho)
  -{\mathcal P}(\rho,u+\bar{\sigma}){\mathcal S}(\rho)
  {\mathcal G}(u+\bar{\sigma})\\
  &=&-c^{-1}{\mathcal F}(u+\bar{\sigma})-\gamma {\mathcal M}(\rho,u
  +\bar{\sigma})\Pi[{\mathcal D}(\rho){\mathcal T}(\rho){\mathcal K}_1(\rho)]
  \\
  &&-{\mathcal M}(\rho,u+\bar{\sigma})
  \Pi[{\mathcal D}(\rho){\mathcal S}(\rho){\mathcal G}(u+\bar{\sigma})].
\end{array}
$$
  Differentiating this expression in $u$ at $(\rho,u)=(0,u_s)$ yields
$$
\begin{array}{rcl}
  D_u{\mathcal F}_2(0,u_s)v &=&-c^{-1}f'(\sigma_s(r))v-\gamma
  {\mathcal M}(0,v)\Pi[{\mathcal D}(0){\mathcal T}(0){\mathcal K}_1(0)]\\
  &&-{\mathcal M}(0,v)\Pi[{\mathcal D}(0){\mathcal S}(0)g(\sigma_s(r))]
  -{\mathcal M}(0,u_s)\Pi[{\mathcal D}(0){\mathcal S}(0)g'(\sigma_s(r))v].
\end{array}
$$
  We have ${\mathcal D}(0){\mathcal T}(0){\mathcal K}_1(0)={\mathcal D}(0)
  {\mathcal T}(0)1={\mathcal D}(0)1=0$, and, by denoting $v_s(r)=p_s(r)-
  p_s(1)$, ${\mathcal D}(0){\mathcal S}(0)g(\sigma_s(r))={\mathcal D}(0)v_s
  =p_s'(1)=0$. Hence,
\setcounter{equation}{8}
\begin{eqnarray}
  D_u{\mathcal F}_2(0,u_s)v &=&-c^{-1}f'(\sigma_s(r))v-{\mathcal M}(0,u_s)
  \Pi[{\mathcal D}(0){\mathcal S}(0)g'(\sigma_s(r))v]
  \nonumber\\
  &=&-c^{-1}f'(\sigma_s(r))v-\phi(r-1)\sigma_s'(r)\Pi[{\mathcal D}(0)
  {\mathcal S}(0)g'(\sigma_s(r))v].
\end{eqnarray}
  In order to compute $D_\rho{\mathcal F}_2(0,u_s)$ we write
\setcounter{equation}{9}
\begin{eqnarray}
  &&{\mathcal D}(0){\mathcal T}(0){\mathcal K}_1'(0)\eta=
  -{\mathcal D}(0){\mathcal T}(0)\eta=-{\mathcal D}(0)\Pi(\eta)
  =-{\partial\Pi(\eta)\over\partial r}\Big|_{r=1},
  \nonumber\\
  &&{\mathcal D}(0)[{\mathcal T}'(0)\eta]{\mathcal K}_1(0)=
  {\mathcal D}(0)[{\mathcal T}'(0)\eta]1=0\;\;\;
  (\mbox{because}\;\;{\mathcal T}(\epsln\eta)1={\mathcal T}(0)1=1),
  \nonumber\\
  &&[{\mathcal D}'(0)\eta]{\mathcal T}(0){\mathcal K}_1(0)=
  [{\mathcal D}'(0)\eta]{\mathcal T}(0)1=[{\mathcal D}'(0)\eta]1=0
  \;\;\;(\mbox{because}\;\;{\mathcal D}(\epsln\eta)1={\mathcal D}(0)1=0),
  \nonumber\\
  &&[{\mathcal D}'(0)\eta]{\mathcal S}(0)g(\sigma_s(r))=
  [{\mathcal D}'(0)\eta]v_s=p_s'(1)\eta=0,
  \nonumber\\
  &&{\mathcal D}(0)[{\mathcal S}'(0)\eta]g(\sigma_s(r))={\mathcal D}(0)
  {\mathcal S}(0)[{\mathcal A}'(0)\eta]{\mathcal S}(0)g(\sigma_s(r))=
  {\mathcal D}(0){\mathcal S}(0)[{\mathcal A}'(0)\eta]v_s.
\end{eqnarray}
  In getting (5.10) we used the identity ${\mathcal S}'(0)\eta={\mathcal S}(0)
  [{\mathcal A}'(0)\eta]{\mathcal S}(0)$, which follows from the fact that
  ${\mathcal A}(\rho){\mathcal S}(\rho)=-id$ for any $\rho\in C^2({\mathbb
  S}^{n-1})$. By a similar argument as in the proof of (5.8) we see that
\begin{eqnarray}
  [{\mathcal A}'(0)\eta]v_s & =&-{\mathcal A}(0)[{\mathcal M}(0,p_s)
  \Pi(\eta)]-{\mathcal G}'(\sigma_s)[{\mathcal M}(0,\sigma_s)\Pi(\eta)]
  \nonumber\\
  &=&-\Delta[\phi(r-1)p_s'(r)\Pi(\eta)]-
  g'(\sigma_s(r))\phi(r-1)\sigma_s'(r)\Pi(\eta).
\end{eqnarray}
  Substituting (5.11) into (5.10) we get
$$
\begin{array}{rl}
  {\mathcal D}(0)[{\mathcal S}'(0)\eta]g(\sigma_s(r))&={\mathcal D}(0)
  [\phi(r-1)p_s'(r)\Pi(\eta)]-{\mathcal D}(0){\mathcal S}(0)[g'(\sigma_s(r))
  \phi(r-1)\sigma_s'(r)\Pi(\eta)]\\
   &=-g(\bar{\sigma})\eta-{\mathcal D}(0){\mathcal S}(0)[g'(\sigma_s(r))
  \phi(r-1)\sigma_s'(r)\Pi(\eta)].
\end{array}
$$
  Using these results and the relations ${\mathcal D}(0){\mathcal T}(0)
  {\mathcal K}_1(0)=0$ and ${\mathcal D}(0){\mathcal S}(0)g(\sigma_s(r))=0$,
  we see that
\begin{eqnarray}
  D_\rho{\mathcal F}_2(0,u_s)\eta &=&-\gamma D_\rho{\mathcal M}(0,\sigma_s)\eta
  \cdot\Pi\big({\mathcal D}(0){\mathcal T}(0){\mathcal K}_1(0)\big)
  -\gamma {\mathcal M}(0,\sigma_s)\Pi\big({\mathcal D}(0){\mathcal T}(0)
  {\mathcal K}_1'(0)\eta\big)
\nonumber\\
  &&-\gamma {\mathcal M}(0,\sigma_s)\Pi\big({\mathcal D}(0)
  [{\mathcal T}'(0)\eta]{\mathcal K}_1(0)\big)-\gamma {\mathcal M}(0,\sigma_s)
  \Pi\big([{\mathcal D}'(0)\eta]{\mathcal T}(0){\mathcal K}_1(0)\big)
\nonumber\\
  &&-D_\rho{\mathcal M}(0,\sigma_s)\eta\cdot\Pi\big({\mathcal D}(0){\mathcal S}(0)
  g(\sigma_s(r))\big)
  -{\mathcal M}(0,\sigma_s)\Pi\big({\mathcal D}(0)[{\mathcal S}'(0)\eta]
  g(\sigma_s(r))\big)
\nonumber\\
  &&-{\mathcal M}(0,\sigma_s)\Pi\big([{\mathcal D}'(0)\eta]{\mathcal S}(0)
  g(\sigma_s(r))\big)
\nonumber\\
  &=&\gamma\phi(r-1)\sigma_s'(r)\Pi\Big({\partial\Pi(\eta)\over\partial r}
  \Big|_{r=1}\Big)+\phi(r-1)\sigma_s'(r)\Pi\big({\mathcal D}(0){\mathcal S}(0)
  [g'(\sigma_s(r))
\nonumber\\
  &&\times\phi(r-1)\sigma_s'(r)\Pi(\eta)]\big)+
  g(\bar{\sigma})\phi(r-1)\sigma_s'(r)\Pi\big(\eta\big).
\end{eqnarray}
  Finally, differentiating ${\mathcal G}_2(\rho,u)=-\gamma{\mathcal D}(\rho)
  {\mathcal T}(\rho){\mathcal K}_1(\rho)-{\mathcal D}(\rho){\mathcal S}(\rho)
  {\mathcal G}(u+\bar{\sigma})$ in $u$ at $(\rho,u)=(0,u_s)$ yields
$$
  D_u{\mathcal G}_2(0,u_s)v=-{\mathcal D}(0){\mathcal S}(0)[g'(\sigma_s(r))v],
\eqno{(5.13)}
$$
  and differentiating in $\rho$ gives
\setcounter{equation}{13}
\begin{eqnarray}
  D_\rho{\mathcal G}_2(0,u_s)\eta&=&-\gamma{\mathcal D}(0){\mathcal T}(0)
  {\mathcal K}_1'(0)\eta-\gamma{\mathcal D}(0)[{\mathcal T}'(0)\eta]
  {\mathcal K}_1(0)-\gamma[{\mathcal D}'(0)\eta]{\mathcal T}(0)
  {\mathcal K}_1(0)
\nonumber\\
  &&-[{\mathcal D}'(0)\eta]{\mathcal S}(0)g(\sigma_s(r))
  -{\mathcal D}(0)[{\mathcal S}'(0)\eta]g(\sigma_s(r))
\nonumber\\
  &&=\gamma {\partial\Pi(\eta)\over\partial r}
  \Big|_{r=1}+{\mathcal D}(0){\mathcal S}(0)[g'(\sigma_s(r))
  \phi(r-1)\sigma_s'(r)\Pi(\eta)]+g(\bar{\sigma})\eta.
\end{eqnarray}
  From (5.1)--(5.5), (5.7)--(5.9) and (5.12)--(5.14) we obtain
$$
  {\mathbb F}'(U_s)=\left(
\begin{array}{cc}
  {\mathcal A}_{11} & \;\; {\mathcal A}_{12}\\
  {\mathcal A}_{21} & \;\; {\mathcal A}_{22}
\end{array}
  \right),
\eqno{(5.15)}
$$
  where, by denoting $m(r)=\phi(r-1)\sigma_s'(r)$,
\begin{eqnarray*}
  &{\mathcal A}_{11}v=&c^{-1}[\Delta-f'(\sigma_s(r))]v-m(r)
  \Pi\big[{\mathcal D}(0){\mathcal S}(0)[g'(\sigma_s(r))v]\big],\\[0.1cm]
  &{\mathcal A}_{12}\eta=& m(r)
  \Pi\Big[\gamma{\partial\over\partial r}\Pi\Big(\eta+{1\over n\!-\!1}
  \Delta_\omega\eta\Big)\Big|_{r=1}+g(\bar{\sigma})\eta\Big]
  -c^{-1}[\Delta-f'(\sigma_s(r))]
  \big[m(r)\Pi(\eta)\big]\\[0.1cm]
  &&+m(r)\Pi\big[{\mathcal D}(0){\mathcal S}(0)[g'(\sigma_s(r))m(r)
  \Pi(\eta)]\big],\\[0.1cm]
  &{\mathcal A}_{21}v=&-{\mathcal D}(0){\mathcal S}(0)[g'(\sigma_s(r))v],\\
  &{\mathcal A}_{22}\eta=&\gamma{\partial\over\partial r}
  \Pi\Big(\eta+{1\over n\!-\!1}\Delta_\omega\eta\Big)\Big|_{r=1}
  +g(\bar{\sigma})\eta
  +{\mathcal D}(0){\mathcal S}(0)[g'(\sigma_s(r))m(r)\Pi(\eta)].
\end{eqnarray*}
  We summarize the above result in the following lemma:
\medskip

  {\bf Lemma 5.1}\ \ {\em The Fr\'{e}chet derivative ${\mathbb F}'(U_s)$ is
  given by (5.15).}
\medskip

  In Section 4 we proved, by using the relation (4.14), that $W_j$ ($j=1,2,
  \cdots,n$) given in (4.16) are eigenvectors of ${\mathbb F}'(U_s)$
  corresponding to the eigenvalue $0$. We can easily reprove this result
  by using the expression (5.15) of ${\mathbb F}'(U_s)$.
\medskip

  {\bf Lemma 5.2}\ \ {\em The operator ${\mathbb F}'(U_s)$, regarded as an
  unbounded linear operator in $X$ with domain $X_0$, is a sectorial operator.}
\medskip

  {\em Proof}:\ \ Recalling (5.1), we see that ${\mathbb F}'(U_s)=A+B$, where
  $A={\mathbb A}(U_s)$ and $BV=[{\mathbb A}'(U_s)V]U_s+{\mathbb F}_0'(U_s)V$
  for $V\in X_0$. By Lemma 4.1 of \cite{Cui6} we know that $A$ is a sectorial
  operator in $X$ with domain $X_0$. Next we consider $B$. Since ${\mathbb F}_0
  \in C^\infty({\mathcal O},Y)$, we have ${\mathbb F}_0'(U_s)\in L(X_0,Y)$.
  Besides, from the second relation in (5.2) and the result obtained in (5.8)
  we easily see that the mapping $V\to [{\mathbb A}'(U_s)V]U_s$ also belong to
  $L(X_0,Y)$. Hence, in conclusion we have $B\in L(X_0,Y)$. Since $Y$ is
  clearly an intermediate space between $X_0$ and $X$, by a standard result we
  get the desired assertion. $\qquad\Box$
\medskip

  By a standard perturbation result, we have
\medskip

  {\bf Corollary 5.3}\ \ {\em If the neighborhood ${\mathcal O}$ of $U_s$ (in
  $X_0$) is sufficiently small, then for any $U\in {\mathcal O}$, ${\mathbb
  F}'(U)$ is a sectorial operator.}
\medskip

  Later on we shall assume that the number $\delta$ is so small that the open
  set ${\mathcal O}$ defined in the end of Section 3 satisfies the condition
  of the above corollary.

\section{The spectrum  of ${\mathbb F}'(U_s)$}

  Given a closed linear operator $B$ in a Banach space $X$, we denote by
  $\rho(B)$ and $\sigma(B)$ the resolvent set and the spectrum of $B$,
  respectively. In the sequel we study $\sigma({\mathbb F}'(U_s))$.

  We introduce the operator ${\mathcal A}_0:W^{m-1,q}({\mathbb B}^n)\to
  W^{m-3,q}({\mathbb B}^n)$ by
$$
  {\mathcal A}_0 v=[\Delta-f'(\sigma_s(r))]v \quad \mbox{for}\;\;
  v\in W^{m-1,q}({\mathbb B}^n),
$$
  the operator $Q: B^{m-1/q}_{qq}({\mathbb S}^{n-1})\to W^{m,q}({\mathbb B}^n)$
  by
$$
  Q\eta=m(r)\Pi(\eta)=\phi(r-1)\sigma_s'(r)\Pi(\eta) \quad \mbox{for}\;\;
   \eta\in B^{m-1/q}_{qq}({\mathbb S}^{n-1}),
$$
  and the operator ${\mathcal J}: W^{m-1,q}({\mathbb B}^n)\to B^{m-1/q}_{qq}
  ({\mathbb S}^{n-1})$ by
$$
  {\mathcal J} v=-{\mathcal D}(0){\mathcal S}(0)[g'(\sigma_s(r))v]
  ={\partial\over\partial r}\Big\{\Delta^{-1}[g'(\sigma_s(r))v]\Big\}
  \Big|_{r=1}\quad \mbox{for}\;\; v\in W^{m-1,q}({\mathbb B}^n)
$$
  Here $\Delta^{-1}$ denotes the inverse of the Laplacian under the homogeneous
  Dirichlet boundary condition. Let $\Pi_0:B^{m-1/q}_{qq}({\mathbb S}^{n-1})\to
  W^{m,q}({\mathbb B}^n)$ be the operator $\Pi_0(\eta)=v$, where for given
  $\eta\in B^{m-1/q}_{qq}({\mathbb S}^{n-1})$, $v\in W^{m,q}({\mathbb B}^n)$ is
  the solution of the boundary value problem
$$
  \Delta v-f'(\sigma_s(r))v=0\;\;\;\mbox{in}\;\;{\mathbb B}^n \quad
  \mbox{and} \quad v=\eta\;\;\;\mbox{on}\;\;{\mathbb S}^{n-1}.
$$
  Note that this definition implies that ${\mathcal A}_0\Pi_0=0$. We define
  ${\mathcal B}_\gamma: B^{m-1/q}_{qq}({\mathbb S}^{n-1})\to B^{m-3-1/q}_{qq}
  ({\mathbb S}^{n-1})$ by
$$
\begin{array}{rl}
  {\mathcal B}_\gamma\eta=&\gamma\displaystyle{\partial\over\partial r}
  \Big\{\Pi\Big(\eta+{1\over n\!-\!1}\Delta_\omega\eta\Big)\Big\}\Big|_{r=1}
  +g(\bar{\sigma})\eta-\sigma_s'(1){\mathcal J}\Pi_0(\eta)\\[0.2cm]
  =&\displaystyle{\partial\over\partial r}
  \Big\{\gamma\Pi\Big(\eta+{1\over n\!-\!1}\Delta_\omega\eta\Big)
  -\sigma_s'(1)\Delta^{-1}\Big(g'(\sigma_s(r))\Pi_0(\eta)\Big)\Big\}\Big|_{r=1}
  +g(\bar{\sigma})\eta,\\
  &\qquad\qquad\qquad\qquad\qquad\qquad\qquad
    \mbox{for}\;\;\eta\in B^{m-1/q}_{qq}({\mathbb S}^{n-1}).
\end{array}
$$
  Finally we define the operators ${\mathbb M}:X_0\to X$ and ${\mathbb T}:X\to
  X$ respectively by
$$
  {\mathbb M}=\left(
\begin{array}{cc}
  c^{-1}{\mathcal A}_0\!+\!\sigma_s'(1)\Pi_0{\mathcal J} & \;\;\;
  \sigma_s'(1)\Pi_0{\mathcal B}_\gamma\\
  {\mathcal J} & \;\;\; {\mathcal B}_\gamma
\end{array}
  \right) \quad \mbox{and} \quad
  {\mathbb T}=\left(
\begin{array}{cc}
  I & \;\;\; \sigma_s'(1)\Pi_0\!-\!Q\\
  0 & \;\;\; I
\end{array}
  \right).
$$
  Here the first $I$ in ${\mathbb T}$ represents the identity operator in
  $W^{m-3,q}({\mathbb B}^n)$, while the second $I$ in ${\mathbb T}$ represents
  the identity operator in $B^{m-3-1/q}_{qq}({\mathbb S}^{n-1})$. Note that
  $(\sigma_s'(1)\Pi_0\!-\!Q)\eta|_{{\mathbb S}^{n-1}}=0$ for any $\eta\in
  B^{m-1/q}_{qq}({\mathbb S}^{n-1})$, so that ${\mathbb T}$ maps $X_0$ to $X_0$.
\medskip

  {\bf Lemma 6.1}\ \ {\em For $V\in X_0$ and $\lambda\in {\mathbb C}$, the
  relation ${\mathbb F}'(U_s)V=\lambda V$ holds if and only if the relations
  ${\mathbb M}W=\lambda W$ and $W={\mathbb T}V$ hold.}
\medskip

  {\em Proof}:\ \ Clearly,
$$
\begin{array}{c}
  {\mathcal A}_{11}v=c^{-1}{\mathcal A}_0v+Q{\mathcal J}v,\quad
  {\mathcal A}_{12}\eta=-c^{-1}{\mathcal A}_0Q\eta
  +Q({\mathcal B}_\gamma+\sigma_s'(1){\mathcal J}\Pi_0-{\mathcal J}Q)\eta,\\
  {\mathcal A}_{21}v={\mathcal J}v,\quad
  {\mathcal A}_{22}\eta=({\mathcal B}_\gamma+\sigma_s'(1){\mathcal J}\Pi_0
  -{\mathcal J}Q)\eta.
\end{array}
$$
  Using these relations and the fact that ${\mathcal A}_0\Pi_0=0$ we can easily
  verify that
$$
  \left(
\begin{array}{cc}
  {\mathcal A}_{11} & \;\; {\mathcal A}_{12}\\
  {\mathcal A}_{21} & \;\; {\mathcal A}_{22}
\end{array}
  \right)=\left(
\begin{array}{cc}
  I & \;\;\; Q\!-\!\sigma_s'(1)\Pi_0\\
  0 & \;\;\; I
\end{array}
  \right)\left(
\begin{array}{cc}
  c^{-1}{\mathcal A}_0\!+\!\sigma_s'(1)\Pi_0{\mathcal J} & \;\;\;
  \sigma_s'(1)\Pi_0{\mathcal B}_\gamma\\
  {\mathcal J} & \;\;\; {\mathcal B}_\gamma
\end{array}
  \right)\left(
\begin{array}{cc}
  I & \;\;\; \sigma_s'(1)\Pi_0\!-\!Q\\
  0 & \;\;\; I
\end{array}
  \right),
$$
  or ${\mathbb F}'(U_s)={\mathbb T}^{-1}{\mathbb M}{\mathbb T}$.
  From this relation the desired assertion follows immediately. $\qquad\Box$
\medskip

  Since $X_0$ is clearly compactly embedded in $X$, by Lemma 5.2 we see that
  $\sigma({\mathbb F}'(U_s))$ consists entirely of eigenvalues. Hence, by Lemma
  6.1 we have
\medskip

  {\bf Corollary 6.2}\ \ {\em $\sigma({\mathbb F}'(U_s))=\sigma({\mathbb M})$.}
  $\qquad\Box$
\medskip

  We shall see that for sufficiently small $c$, $\sigma({\mathcal B}_\gamma)$
  plays a major role in determining $\sigma({\mathbb M})$. Hence, in the sequel
  we first compute $\sigma({\mathcal B}_\gamma)$. To this end we introduce some
  notation and recall some results of \cite{CuiEsc1}. For every nonnegative integer
  $k$, let $Y_{kl}(\omega)$, $l=1,2,\cdots,d_k$, be the normalized orthogonal
  basis of the space of all spherical harmonics of degree $k$, where $d_k$ is
  the dimension of this space, i.e.
$$
  d_0=1, \quad d_1=n, \quad d_k=\Big(\begin{array}{c} n\!+\!k\!-\!1\\[0.1cm]
   k\end{array}\Big)-\Big(\begin{array}{c} n\!+\!k\!-\!3\\[0.1cm]
   k\!-\!2\end{array}\Big)\;\;(k\geq 2).
$$
  It is well-known that
$$
  \Delta_\omega Y_{kl}(\omega)=-\lambda_k Y_{kl}(\omega), \quad
  \lambda_k=k^2+(n-2)k\;\;\;(k=0,1,2,\cdots),
$$
  and $\lambda_k$ ($k=0,1,2,\cdots$) are the all eigenvalues of $\Delta_\omega$.
  We denote
$$
  a_k=2k+n-1\geq n-1,
$$
  and denote by $\bar{u}_k(r)$ the solution of the initial value problem
$$
\left\{
\begin{array}{l}
  \bar{u}_k''(r)+\displaystyle\frac{a_k}{r}\bar{u}_k'(r)=
  f'(\sigma_s(r))\bar{u}_k(r), \\[0.2cm]
  \bar{u}_k(0)=1, \quad \bar{u}_k'(0)=0,
\end{array}
\right.
$$
  By using some ODE techniques we can show that this problem has a unique
  solution for all $r\in [0,R^\ast)$, where $[0,R^\ast)$ is the maximal
  existence interval of $\sigma_s(r)$. We also denote
$$
  \gamma_k={n-1\over (\lambda_k\!-\!n\!+\!1)k}
  \Big[g(\bar{\sigma})-{\sigma_0'(1)\over \bar{u}_k(1)}
  \int_0^1 g'(\sigma_0(\rho))\bar{u}_k(\rho)\rho^{a_k}d\rho\Big]
  \quad (k\geq 2).
$$
  From \cite{CuiEsc1} we know that $\gamma_k$'s and $\gamma_1=0$ are the all
  eigenvalues of the linearization of the stationary version of the system
  (1.1)--(1.5) at the radially symmetric stationary solution $(\sigma_s,p_s,
  \Omega_s)$, $\gamma_k>0$ for all $k\geq 2$, and $\lim_{k\to\infty}\gamma_k=0$.
  Next we denote
$$
  \alpha_{k,\gamma}=-\frac{(\lambda_k\!-\!n\!+\!1)k}{n-1}
  (\gamma-\gamma_k), \quad k=2,3,\cdots.
$$
  Note that $\alpha_{k,\gamma}\sim -\gamma k^3/(n\!-\!1)$ as $k\to\infty$.
  Finally, we denote $\alpha_{1,\gamma}=0$ and
$$
  \alpha_{0,\gamma}= g(\bar{\sigma})-{\sigma_0'(R_s)\over
  \bar{u}_0(R_s)R_s^{n\!-\!1}}\int_0^{R_s} g'(\sigma_0(r))
  \bar{u}_0(r)r^{n\!-\!1}dr.
$$
  From [16] we know that $\alpha_{0,\gamma}<0$ for all $\gamma>0$.
\medskip

  {\bf Lemma 6.3}\ \ {\em ${\mathcal B}_\gamma$ is a Fourier multiplication
  operator of the following form: For any $\eta(\omega)=\sum_{k=0}^\infty
  \sum_{l=1}^{d_k}b_{kl}Y_{kl}(\omega)\in C^\infty({\mathbb S}^{n-1})$,
$$
  {\mathcal B}_\gamma\eta(\omega)=\sum_{k=0}^\infty\sum_{l=1}^{d_k}
  \alpha_{k,\gamma}b_{kl}Y_{kl}(\omega).
\eqno{(6.1)}
$$
  As a result, we have $\sigma({\mathcal B}_\gamma)=\{\alpha_{k,\gamma}:
  k\in {\mathbb N,\,k\geq 2}\}\cup\{\,0,\alpha_{0,\gamma}\}$.}
\medskip

  {\em Proof}:\ \ It can be easily seen that ${\mathcal B}_\gamma$ has the
  same expression as that introduced in \cite{CuiEsc2} with the same notation
  (but notice that ${\mathcal B}_\gamma$ in \cite{CuiEsc2} is a mapping from
  $C^{m+\mu}({\mathbb S}^{n-1})$ to $C^{m-3+\mu}({\mathbb S}^{n-1})$ for some
  inter $m$ and $0<\mu<1$). Hence, by a similar calculation as in
  \cite{CuiEsc2} we get (6.1). $\qquad\Box$
\medskip

  {\bf Lemma 6.4}\ \ {\em For any $\gamma>0$ and $k\geq 2$, there exists
  corresponding $c_0>0$ such that for any $0<c\leq c_0$, ${\mathbb M}$ has an
  eigenvalue $\lambda_{k,\gamma}=\alpha_{k,\gamma}+c\mu_{k,\gamma}(c)$, where
  $\mu_{k,\gamma}(c)$ is a bounded continuous function in $0<c\leq c_0$.
  Moreover, the corresponding eigenvectors of ${\mathbb M}$ have the form
  $\Big(^{\;ca_{k,\gamma}(r,c)}_{\;\;\quad 1\quad}\Big)Y_{kl}(\omega)$
  $(l=1,2,\cdots,d_k)$, where $a_{k,\gamma}(r,c)$ is a smooth function in
  $r\in [0,1]$ and is bounded and continuous in  $0<c\leq c_0$.}
\medskip

  {\em Proof}:\ \ Let $U=\Big(^{\;ca_{k,\gamma}(r,c)}_{\;\;\quad 1\quad}\Big)
  Y_{kl}(\omega)$. Then the relation ${\mathbb M}U=\lambda_{k,\gamma}U$ holds
  if and only if the following relations hold:
\setcounter{equation}{1}
\begin{eqnarray}
  {\mathcal A}_0(a_{k,\gamma}Y_{kl})+c\sigma_s'(1)\Pi_0{\mathcal J}
  (a_{k,\gamma}Y_{kl})+\sigma_s'(1)\Pi_0{\mathcal B}_\gamma(Y_{kl}) &=&
  c\alpha_{k,\gamma}a_{k,\gamma}Y_{kl}+c^2\mu_{k,\gamma}a_{k,\gamma}Y_{kl},
  \\
  c{\mathcal J}(a_{k,\gamma}Y_{kl})+{\mathcal B}_\gamma(Y_{kl}) &=&
  \alpha_{k,\gamma}Y_{kl}+c\mu_{k,\gamma}Y_{kl}.
\end{eqnarray}
  Let ${\mathcal L}_k$ be the second-order differential operator ${\mathcal
  L}_k u(r)=u''(r)+{n-1\over r}u'(r)-{\lambda_k\over r^2}u(r)$, and ${\mathcal
  J}_k$ be the operator $u\to v_k'(1)$, where for a given continuous function
  $u=u(r)$ ($0\leq r\leq 1$), $v=v_k(r)$ is the solution of the boundary value
  problem:
$$
\left\{
\begin{array}{l}
  v''(r)+{n-1\over r}v'(r)-{\lambda_k\over r^2}v(r)=g'(\sigma_s(r))u(r),
  \quad 0<r<1,\\
  v'(0)=0,\quad v(1)=0.
\end{array}
\right.
$$
  Then we have ${\mathcal A}_0(a_{k,\gamma}Y_{kl})={\mathcal L}_k(a_{k,\gamma})
  Y_{kl}$ and ${\mathcal J}(a_{k,\gamma}Y_{kl})={\mathcal J}_k(a_{k,\gamma})
  Y_{kl}$. Besides, it can be easily seen that $\Pi_0(Y_{kl})=w_k(r)Y_{kl}$,
  where $w_k(r)$ ($0\leq r\leq 1$) is the solution of the boundary value
  problem:
$$
\left\{
\begin{array}{l}
  w_k''(r)+{n-1\over r}w_k'(r)-\Big({\lambda_k\over r^2}+
  f'(\sigma_s(r))\Big)w_k(r)=0,
  \quad 0<r<1,\\
  w_k'(0)=0,\quad w_k(1)=0.
\end{array}
\right.
$$
  Using these facts and the relation ${\mathcal B}_\gamma(Y_{kl})=
  \alpha_{k,\gamma}Y_{kl}$ (cf. (6.1)) we see that (6.2) and (6.3) reduce to
  the following system of equations:
\begin{eqnarray*}
  {\mathcal L}_k(a_{k,\gamma})+c\sigma_s'(1){\mathcal J}_k(a_{k,\gamma})
  w_k(r)+\sigma_s'(1)\alpha_{k,\gamma}w_k(r) &=&
  c\alpha_{k,\gamma}a_{k,\gamma}+c^2\mu_{k,\gamma}a_{k,\gamma},
  \\
  \mu_{k,\gamma}&=&{\mathcal J}_k(a_{k,\gamma}),
\end{eqnarray*}
  which can be further reduced to the following scaler equation in
  $a_{k,\gamma}$:
$$
  {\mathcal L}_k(a_{k,\gamma})=-c\sigma_s'(1){\mathcal J}_k(a_{k,\gamma})
  w_k(r)+c\alpha_{k,\gamma}a_{k,\gamma}
  +c^2a_{k,\gamma}{\mathcal J}_k(a_{k,\gamma})
  -\sigma_s'(1)\alpha_{k,\gamma}w_k(r).
$$
  By using a standard fixed point argument we can easily show that for $c$
  sufficiently small this equation complemented with the boundary value
  conditions ${\partial a_{k,\gamma}\over\partial r}\Big|_{r=0}=0$ and
  $a_{k,\gamma}\Big|_{r=1}=0$ has a unique solution. By this assertion, the
  desired result follows immediately. $\qquad\Box$

\medskip

  We denote
$$
  \gamma_\ast=\max_{k\geq 2}\gamma_k \quad \mbox{and} \quad
  \alpha^\ast_\gamma =\max_{k\geq 2}\alpha_{k,\gamma}.
$$
  Since $\gamma_k>0$, $\lim_{k\to\infty}\gamma_k=0$ and $\lim_{k\to\infty}
  \alpha_{k,\gamma}=-\infty$, $\gamma_\ast$ and $\alpha^\ast_\gamma $ are both
  well-defined. Clearly, we have $\alpha^\ast_\gamma <0$ for any $\gamma>
  \gamma_\ast$, while $\alpha^\ast_\gamma >0$ for any $0<\gamma<\gamma_\ast$.
\medskip

  {\bf Lemma 6.5}\ \ {\em Given $\gamma>\gamma_\ast$, there exists
  corresponding $c_0>0$ such that for any $0<c\leq c_0$ and any $\lambda\in
  {\mathbb C}\backslash\{\,0\}$ satisfying ${\rm Re}\lambda\geq {1\over 2}
  \alpha^\ast_\gamma $, there holds $\lambda\in\rho({\mathbb M})$, or
  equivalently,}
$$
  \sup\{{\rm Re}\lambda:\lambda\in\sigma({\mathbb M})\backslash\{\,0\}\}\leq
  {1\over 2}\alpha^\ast_\gamma<0.
$$

  {\em Proof}:\ \ We denote
$$
  {\mathbb M}_0=\left(
\begin{array}{cc}
  c^{-1}{\mathcal A}_{0}\;\;
  &\;\; 0\\
  {\mathcal J}\;\; &\;\; \mathcal B_\gamma
\end{array}
\right),
  \qquad
  {\mathbb N}=\left(
\begin{array}{cc}
  \sigma_s'(1)\Pi_0{\mathcal J}\;\;
  &\;\;\sigma_s'(1)\Pi_0\mathcal B_\gamma\\
  0\;\; &\;\; 0
\end{array}
\right).
$$
  Then ${\mathbb M}_0\in L({X}_0,{X})$, ${\mathbb N}\in L({X}_0,{X})$, and
  ${\mathbb M}={\mathbb M}_0+{\mathbb N}$. Since $f'(\sigma_s(r))\geq 0$, From
  the standard theory of elliptic partial differential equations of the
  second-order we know that all eigenvalues of ${\mathcal A}_{0}$ are negative
  and they make up a decreasing sequence tending to $-\infty$. Let $\nu_1$ be
  the largest eigenvalue of ${\mathcal A}_{0}$, and let $c_0=\nu_1/\alpha^\ast
  _\gamma$. Then for any $0<c\leq c_0$ and any $\lambda\in {\mathbb C}
  \backslash\{\,0\}$ such
  that ${\rm Re}\lambda\geq {1\over 2}\alpha^\ast_\gamma $ we have ${\rm Re}
  (c\lambda)\geq {1\over 2}\nu_1$, so that both $\lambda I-c^{-1}
  {\mathcal A}_{0}=c^{-1}(c\lambda I-{\mathcal A}_{0})$ and $\lambda I-
  {\mathcal B}_\gamma$ are invertible, which implies that $\lambda I-
  {\mathbb M}_0$ is invertible. In fact,
$$
  (\lambda I-{\mathbb M}_0)^{-1}=\left(
\begin{array}{cc}
  (\lambda I-c^{-1}{\mathcal A}_{0})^{-1}\;\;
  &\;\; 0\\
  (\lambda I-{\mathcal B}_\gamma)^{-1}
  {\mathcal J}(\lambda I-c^{-1}{\mathcal A}_{0})^{-1}\;\;
  &\;\; (\lambda I-{\mathcal B}_\gamma)^{-1}
\end{array}
\right).
$$
 Hence
$$
  \lambda I-{\mathbb M}=(\lambda I-{\mathbb M}_0)
  -{\mathbb N}=(\lambda I-{\mathbb M}_0)
  (I-c{\mathbb K}),
$$
  where
$$
\begin{array}{rl}
  {\mathbb K}=&c^{-1}(\lambda I-{\mathbb M}_0)^{-1}{\mathbb N}\\
  =&\left(
\begin{array}{cc}
  (c\lambda I-{\mathcal A}_{0})^{-1}\sigma_s'(1)\Pi_0{\mathcal J}\;\;
  &\;\; (c\lambda I-{\mathcal A}_{0})^{-1}\sigma_s'(1)\Pi_0{\mathcal B}_\gamma
  \\[0.2cm]
  (\lambda I-{\mathcal B}_\gamma)^{-1}{\mathcal J}
  (c\lambda I-{\mathcal A}_{0})^{-1}\sigma_s'(1)\Pi_0{\mathcal J}\;\;
  &\;\; (\lambda I-{\mathcal B}_\gamma)^{-1}{\mathcal J}
  (c\lambda I-{\mathcal A}_{0})^{-1}\sigma_s'(1)\Pi_0{\mathcal B}_\gamma
\end{array}
\right).
\end{array}
$$
  Since ${\mathcal A}_{0}$ is a self-adjoint sectorial operator and $\nu_1$
  is the maximal eigenvalue of ${\mathcal A}_{0}$, we have
$$
  \|(c\lambda I-{\mathcal A}_{0})^{-1}\|_{L(W^{m-3,q}({\mathbb B}^n),
  W^{m-3,q}({\mathbb B}^n))}\leq {C\over |c\lambda-\nu_1|}\leq 2C/\nu_1,
$$
  where $C$ is a constant independent of $c$ and $\lambda$. Using this fact,
  the identity
$$
  {\mathcal A}_{0}(c\lambda I-{\mathcal A}_{0})^{-1}=
  c\lambda (c\lambda I-{\mathcal A}_{0})^{-1}-I,
$$
  and the Agmon-Douglis-Nirenberg inequality, we obtain
$$
\begin{array}{rl}
  &\|(c\lambda I-{\mathcal A}_{0})^{-1}\|_{L(W^{m-3,q}({\mathbb B}^n),
  W^{m-1,q}({\mathbb B}^n)\cap W^{1,q}_0({\mathbb B}^n))}\\
  \leq & C[\|(c\lambda I-{\mathcal A}_{0})^{-1}\|_{L(W^{m-3,q}({\mathbb B}^n),
  W^{m-3,q}({\mathbb B}^n))}+
  \|{\mathcal A}_{0}(c\lambda I-{\mathcal A}_{0})^{-1}\|_{L(W^{m-3,q}({\mathbb B}^n),
  W^{m-3,q}({\mathbb B}^n))}]\\ [0.1cm]
  \leq &\displaystyle C+ {C|c\lambda|\over |c\lambda-\nu_1|}\leq C.
\end{array}
$$
  Similarly we have
$$
  \|(\lambda I-{\mathcal B}_\gamma)^{-1}\|_{L({B}^{m-3-1/q}_{qq}
  ({\mathbb S}^{n-1}),{B}^{m-1/q}_{qq}({\mathbb S}^{n-1}))}\leq C.
$$
  Using these estimates we can easily show that
$$
  \|{\mathbb K}\|_{L({X}_0,{X}_0)}\leq C
$$
  for any $0<c\leq c_0$ and any $\lambda\in {\mathbb C}$ such that ${\rm Re}
  \lambda\geq {1\over 2}\alpha^\ast_\gamma $. It follows that if we take
  $c_0$ further small such that $c_0 C<1$ then for $c$ and $\lambda$ in the
  set specified above, the operator $\lambda I-{\mathbb M}$ is invertible and
  the inverse is continuous. Hence, the desired assertion follows. $\quad\Box$
\medskip

\section{The proof of Theorem 1.1}

  {\em Proof of Theorem 1.1}:\ \ We first assume that $\gamma>\gamma_\ast$. By
  Lemma 5.2 we see that ${\mathbb F}'(U_s)$ is a sectorial operator in $X$ with
  domain $X_0$. In what follows we prove that the norm of $X_0$ coincides the
  graph norm of ${\mathbb F}'(U_s)$. From Section 6 we see that
  ${\mathbb F}'(U_s)={\mathbb T}^{-1}{\mathbb M}{\mathbb T}$. Clearly,
$$
  C\|U\|_X\leq\|{\mathbb T}U\|_X\leq C^{-1}\|U\|_X \quad \mbox{and} \quad
  C\|U\|_{X_0}\leq\|{\mathbb T}U\|_{X_0}\leq C^{-1}\|U\|_X
\eqno{(7.1)}
$$
  for some constants $C>0$. Thus the graph norm of ${\mathbb F}'(U_s)$ is
  equivalent to the graph norm of ${\mathbb M}$. Next, let
$$
  {\mathbb T}_0=\left(
\begin{array}{cc}
  I\;\;
  &\;\;\sigma_s'(1)\Pi_0\\
  0\;\; &\;\; I
\end{array}
\right).
$$
  Then we have ${\mathbb M}={\mathbb T}_0{\mathbb M}_0$. Clearly, all estimates
  in (7.1) still hold when ${\mathbb T}$ is replaced by ${\mathbb T}_0$. Hence
  the graph norm of ${\mathbb M}$ is equivalent to the graph norm of ${\mathbb
  M}_0$. Clearly, as an unbounded linear operator in $W^{m-3,q}({\mathbb B}^n)$
  with domain $W^{m-1,q}({\mathbb B}^n)$, the graph norm of ${\mathcal A}_0$ is
  equivalent to the norm of $W^{m-1,q}({\mathbb B}^n)$. Also, we know that as
  an unbounded linear operator in $B^{m-3-1/q}_{qq}({\mathbb S}^{n-1})$ with
  domain $B^{m-1/q}_{qq}({\mathbb S}^{n-1})$, the graph norm of ${\mathcal
  B}_\gamma$ is equivalent to the norm of $B^{m-1/q}_{qq}({\mathbb S}^{n-1})$
  (cf. \cite{CuiEsc2}). Besides, it is easy
  to see that ${\mathcal J}$ maps $W^{m-3,q}({\mathbb B}^n)$  continuously into
  $B^{m-2-1/q}_{qq}({\mathbb S}^{n-1})$, so that it is a compact operator from
  $W^{m-3,q}({\mathbb B}^n)$ to $B^{m-3-1/q}_{qq}({\mathbb S}^{n-1})$. From
  these facts, we can easily show that the graph norm of ${\mathbb M}_0$ is
  equivalent to the norm of $X_0$. Hence, the graph norm of ${\mathbb F}'(U_s)$
  is equivalent to the norm of $X_0$. This verifies that ${\mathbb F}'(U_s)$
  satisfies the condition $(B_1)$. By the results of Section 4 we see that
  ${\mathbb F}'(U_s)$ also satisfies the condition $(B_2)$. Next we consider
  the condition $(B_3)$. We first prove that ${\mathbb M}$ satisfies this
  condition. To this end we denote by ${\mathbb H}_1({\mathbb S}^{n-1})$ the
  linear space of all first-order spherical harmonics, and for every integer
  $k$ we introduce
$$
  \hat{B}^{k-1/q}_{qq}({\mathbb S}^{n-1})=
  \{\rho\in B^{k-1/q}_{qq}({\mathbb S}^{n-1}):\rho\;\;\mbox{is orthogonal to}
  \;\; {\mathbb H}_1({\mathbb S}^{n-1})\;\;
  \mbox{in}\;\; L^2({\mathbb S}^{n-1})\}.
$$
  We also denote $\hat{B}^\infty_{qq}({\mathbb S}^{n-1})=\bigcap_{k=1}^\infty
  \hat{B}^{k-1/q}_{qq}({\mathbb S}^{n-1})$.
  It can be easily shown that $\hat{B}^{k-1/q}_{qq}({\mathbb S}^{n-1})$ is a
  closed subspace of $B^{k-1/q}_{qq}({\mathbb S}^{n-1})$, and
$$
  B^{k-1/q}_{qq}({\mathbb S}^{n-1})=\hat{B}^{k-1/q}_{qq}({\mathbb S}^{n-1})
  \oplus {\mathbb H}_1({\mathbf S}^{n-1})
$$
  (for any integer $k$). By (6.1) we see that $\ker({\mathcal B}_\gamma)=
  {\mathbb H}_1({\mathbf S}^{n-1})$. We denote $\hat{{\mathcal B}}_\gamma=
  {\mathcal B}_\gamma|_{\hat{B}^{m-1/q}_{qq}({\mathbb S}^{n-1})}$, and split
  ${\mathcal J}$ into ${\mathcal J}_1+{\mathcal J}_2$ such that ${\mathcal
  J}_1 v\in\hat{B}^{m-3-1/q}_{qq}({\mathbb S}^{n-1})$ and ${\mathcal J}_2 v\in
  {\mathbb H}_1({\mathbb S}^{n-1})$ for any $v\in W^{m-1,q}({\mathbb B}^n)\cap
  W^{1,q}_0({\mathbb B}^n)$. We correspondingly split $X_0$ and $X$ into
  $(W^{m-1,q}({\mathbb B}^n)\cap W^{1,q}_0({\mathbb B}^n))\times
  \hat{B}^{m-1/q}_{qq}({\mathbb S}^{n-1})\times {\mathbb H}_1({\mathbb
  S}^{n-1})$ and $W^{m-3,q}({\mathbb B}^n))\times\hat{B}^{m-3-1/q}_{qq}
  ({\mathbb S}^{n-1})\times {\mathbb H}_1({\mathbb S}^{n-1})$, respectively.
  Then
$$
    {\mathbb M}=\left(
\begin{array}{ccc}
  c^{-1}{\mathcal A}_{0}+\sigma_s'(1)\Pi_0({\mathcal J}_1+{\mathcal J}_2)
  \;\;&\;\; \sigma_s'(1)\Pi_0\hat{{\mathcal B}}_\gamma\;\;&\;\; 0\\
  {\mathcal J}_1\;\; &\;\; \hat{{\mathcal B}}_\gamma\;\;&\;\; 0\\
  {\mathcal J}_2\;\; &\;\; 0\;\;&\;\; 0
\end{array}
\right)=\left(
\begin{array}{cc}
  \hat{{\mathbb M}}\;\;&\;\; 0\\
  \hat{{\mathcal J}}\;\; &\;\; 0
\end{array}
\right),
$$
  where
$$
\begin{array}{rl}
  \hat{{\mathbb M}}=&\left(
\begin{array}{cc}
  c^{-1}{\mathcal A}_{0}+\sigma_s'(1)\Pi_0({\mathcal J}_1+{\mathcal J}_2)
  \;\;&\;\; \sigma_s'(1)\Pi_0\hat{{\mathcal B}}_\gamma\\
  {\mathcal J}_1\;\; &\;\; \hat{{\mathcal B}}_\gamma
\end{array}
\right)\\[0.2cm]
  =&\left(
\begin{array}{cc}
  I\;\;&\;\; \sigma_s'(1)\Pi_0\\
  0\;\; &\;\; I
\end{array}
\right)\left(
\begin{array}{cc}
  c^{-1}{\mathcal A}_{0}+\sigma_s'(1)\Pi_0{\mathcal J}_2
  \;\;&\;\; 0\\
  {\mathcal J}_1\;\; &\;\; \hat{{\mathcal B}}_\gamma
\end{array}
\right)\equiv \hat{{\mathbb T}}_0\hat{{\mathbb M}}_1.
\end{array}
$$
  and $\hat{{\mathcal J}}=({\mathcal J}_2\;\;0)$.
  We claim that $\hat{{\mathcal B}}_\gamma$ is an isomorphism from
  $\hat{B}^{m-1/q}_{qq}({\mathbb S}^{n-1})$ to $\hat{B}^{m-3-1/q}_{qq}
  ({\mathbb S}^{n-1})$. Indeed, from (6.1) and the fact that ${\mathcal
  B}_\gamma$ maps $B^{m-1/q}_{qq}({\mathbb S}^{n-1})$ to
  $B^{m-3-1/q}_{qq}({\mathbb S}^{n-1})$ boundedly it is clear that
  $\hat{{\mathcal B}}_\gamma$ maps $\hat{B}^{m-1/q}_{qq}({\mathbb S}^{n-1})$ to
  $\hat{B}^{m-3-1/q}_{qq}({\mathbb S}^{n-1})$ boundedly and is an injection.
  Next, from (6.1) we see immediately that for any $\zeta\in \hat{B}^\infty_{qq}
  ({\mathbb S}^{n-1})$ there exists a unique $\eta\in\hat{B}^\infty_{qq}
  ({\mathbb S}^{n-1})$ such that ${\mathcal B}_\gamma\eta=\zeta$. Now assume
  that $\zeta\in\hat{B}^{m-3-1/q}_{qq}({\mathbb S}^{n-1})$. Let $\zeta_j\in
  \hat{B}^\infty_{qq}({\mathbb S}^{n-1})$ ($j=1,2,\cdots$) be such that $\zeta_j
  \to\zeta$ in $\hat{B}^{m-3-1/q}_{qq}({\mathbb S}^{n-1})$, and let $\eta_j\in
  \hat{B}^\infty_{qq}({\mathbb S}^{n-1})$ be the solution of the equation
  ${\mathcal B}_\gamma\eta_j=\zeta_j$ ($j=1,2,\cdots$). Take a real number $s$
  such that $s<m-3-1/q-(n\!-\!1)({1\over 2}-{1\over q})$. Then $B^{m-3-1/q}_{qq}
  ({\mathbb S}^{n-1})\hookrightarrow H^s({\mathbb S}^{n-1})$, where
  $H^s({\mathbb S}^{n-1})$ stands for the usual Sobolev space. Thus $\zeta_j\to
  \zeta$ in $H^s({\mathbb S}^{n-1})$. By (6.1) and the fact that $\alpha_{k,
  \gamma}\sim Ck^3$ we easily deduce that $\{\eta_j\}$ is a Cauchy sequence in
  $H^{s+3}({\mathbb S}^{n-1})$. Let $\eta\in H^{s+3}({\mathbb S}^{n-1})$ be the
  limit of $\{\eta_j\}$. By a standard argument we have
$$
  \|\rho\|_{B^{m-1/q}_{qq}({\mathbb S}^{n-1})}\leq
  C\big(\|\rho\|_{H^{s+3}({\mathbb S}^{n-1})}+
  \|{\mathcal B}_\gamma\rho\|_{B^{m-3-1/q}_{qq}({\mathbb S}^{n-1})}\big).
$$
  Applying this estimate to $\rho=\eta_j-\eta$, we conclude that $\eta_j\to
  \eta$ in $B^{m-1/q}_{qq}({\mathbb S}^{n-1})$. Since $\hat{B}^{m-1/q}_{qq}
  ({\mathbb S}^{n-1})$ is closed in $B^{m-1/q}_{qq}({\mathbb S}^{n-1})$, we
  get $\eta\in\hat{B}^{m-1/q}_{qq}({\mathbb S}^{n-1})$. This shows that
  $\hat{{\mathcal B}}_\gamma$ is a surjection. Hence, by the Banach inverse
  mapping theorem we see that $\hat{{\mathcal B}}_\gamma$ is an isomorphism
  from $\hat{B}^{m-1/q}_{qq}({\mathbb S}^{n-1})$ to $\hat{B}^{m-3-1/q}_{qq}
  ({\mathbb S}^{n-1})$, as desired. Next, since ${\mathcal A}_0$ is an
  isomorphism from $W^{m-1,q}({\mathbb B}^n)\cap W^{1,q}_0({\mathbb B}^n)$ to
  $W^{m-3,q}({\mathbb B}^n)$ and clearly $\sigma_s'(1)\Pi_0{\mathcal J}_2$ is
  a bounded operator from $W^{m-1,q}({\mathbb B}^n)\cap W^{1,q}_0({\mathbb
  B}^n)$ to $W^{m-3,q}({\mathbb B}^n)$ (actually a compact operator), it
  follows that for $c$ sufficiently small, $c^{-1}{\mathcal A}_{0}+\sigma_s'(1)
  \Pi_0{\mathcal J}_2$ is an isomorphism from $W^{m-1,q}({\mathbb B}^n)\cap
  W^{1,q}_0({\mathbb B}^n)$ to $W^{m-3,q}({\mathbb B}^n)$. By these results
  combined with the fact that ${\mathcal J}_1$ is a bounded operator from
  $W^{m-1,q}({\mathbb B}^n)\cap W^{1,q}_0({\mathbb B}^n)$ to $\hat{B}^{m-3-
  1/q}_{qq}({\mathbb S}^{n-1})$ (actually a compact operator), we immediately
  deduce that $\hat{{\mathbb M}}_1$ is an isomorphism from $(W^{m-1,q}
  ({\mathbb B}^n)\cap W^{1,q}_0({\mathbb B}^n))\times\hat{B}^{m-1/q}_{qq}
  ({\mathbb S}^{n-1})$ to $W^{m-3,q}({\mathbb B}^n)\times\hat{B}^{m-3-1/q}_{qq}
  ({\mathbb S}^{n-1})$. Since $\hat{{\mathbb T}}_0$ is clearly a
  self-isomorphism on $W^{m-3,q}({\mathbb B}^n)\times\hat{B}^{m-3-1/q}_{qq}
  ({\mathbb S}^{n-1})$, we conclude that $\hat{{\mathbb M}}$ is an isomorphism
  from $(W^{m-1,q}({\mathbb B}^n)\cap W^{1,q}_0({\mathbb B}^n))\times
  \hat{B}^{m-1/q}_{qq}({\mathbb S}^{n-1})$ to $W^{m-3,q}({\mathbb B}^n)\times
  \hat{B}^{m-3-1/q}_{qq}({\mathbb S}^{n-1})$. This easily implies that
  ${\mathbb M}$ satisfies the condition $(B_3)$. Now, since ${\mathbb F}'(U_s)
  ={\mathbb T}^{-1}{\mathbb M}{\mathbb T}$, it follows immediately that
  ${\mathbb F}'(U_s)$ also satisfies the condition $(B_3)$. Finally, by
  Corollary 6.2 and Lemma 6.5 we see that
$$
  \omega_-=-\sup\big\{{\rm Re}\lambda:\lambda\in
  \sigma({\mathbb F}'(U_s))\backslash\{\,0\}\big\}>0,
$$
  so that the condition $(B_4)$ is also satisfied by ${\mathbb F}'(U_s)$.
  Hence, by Theorem 2.1 we get the assertion $(i)$ of Theorem 1.1.

  Next we assume that $0<\gamma<\gamma_\ast$. Then there exists $k_0\geq 2$
  such that $\alpha_{k_0,\gamma}>0$. By Lemma 6.4 and Corollary 6.2, this
  implies that for sufficiently small $c$, ${\mathbb F}'(U_s)$ has a positive
  eigenvalue. Furthermore, if $\alpha_{k_1,\gamma}$, $\alpha_{k_2,\gamma}$,
  $\cdots$, $\alpha_{k_N,\gamma}$ are the all positive eigenvalues of ${\mathcal
  B}_\gamma$, then by Lemma 6.4 and a similar argument as in the proof of Lemma
  6.5 we see that for $c$ sufficiently small, $\lambda_{k_j,\gamma}=\alpha_{k_j,
  \gamma}+c\mu_{k_j,\gamma}(c)$ ($j=1,2,\cdots,N$) are the all positive
  eigenvalues of ${\mathbb F}'(U_s)$, and the following estimate holds:
$$
  \sup\{{\rm Re}\lambda:\lambda\in\sigma({\mathbb M})\backslash\{\,0,
  \lambda_{k_1,\gamma},\lambda_{k_2,\gamma},\cdots,\lambda_{k_N,\gamma}\}\}\leq
  {1\over 2}\max\{\alpha_k: k\geq 2, k\neq k_1,k_2,\cdots,k_N\}<0.
$$
  Thus by using Theorem 9.1.3 of \cite{Lunar}, we obtain the assertion $(ii)$
  of Theorem 1.1. This completes the proof of Theorem 1.1. $\qquad\Box$

\medskip
   {\bf Acknowledgement}.\hskip 1em This work is financially supported by the
   National Natural Science Foundation of China under the grant number 10471157.

\end{document}